# $p$-adic boundary values

P. Schneider, J. Teitelbaum

In this paper, we study in detail certain natural continuous representations of $G = GL_n(K)$ in locally convex vector spaces over a locally compact, non-archimedean field $K$ of characteristic zero. We construct boundary value maps, or integral transforms, between subquotients of the dual of a "holomorphic" representation coming from a $p$-adic symmetric space, and "principal series" representations constructed from locally analytic functions on $G$. We characterize the image of each of our integral transforms as a space of functions on $G$ having certain transformation properties and satisfying a system of partial differential equations of hypergeometric type.

This work generalizes earlier work of Morita, who studied this type of representation of the group $SL_2(K)$. It also extends the work of Schneider-Stuhler on the deRham cohomology of $p$-adic symmetric spaces. We view this work as part of a general program of developing the theory of such representations.

A major motivation for studying continuous representations of $p$-adic groups comes from the observation that, in traditional approaches to the representation theory of $p$-adic groups, one separates representations into two disjoint classes – the smooth representations (in the sense of Langlands theory) and the finite dimensional rational representations. Such a dichotomy does not exist for real Lie groups, where the finite dimensional representations are "smooth." The category of continuous representations which we study is broad enough to unify both smooth and rational representations, and one of the most interesting features of our results is the interaction between these two types of representations.

The principal tools of this paper are non-archimedean functional analysis, rigid geometry, and the "residue" theory developed in the paper [ST]. Indeed, the boundary value maps we study are derived from the residue map of [ST].

Before summarizing the structure of our paper and discussing our main results, we will review briefly some earlier, related results.

The pioneering work in this area is due to Morita ([Mo1-Mo6]). He intensively studied two types of representations of $SL_2(K)$. The first class of representations comes from the action of $SL_2(K)$ on sections of rigid line bundles on the one-dimensional rigid analytic space $\mathcal{X}$ obtained by deleting the $K$-rational points from $\mathbb{P}^1_{/K}$; this space is often called the $p$-adic upper half plane. The second class of representations is constructed from locally analytic functions on $SL_2(K)$ which transform by a locally analytic character under the right action by a Borel subgroup $P$ of $SL_2(K)$. This latter class make up what Morita called the ($p$-adic) principal series.

Morita showed that the duals of the "holomorphic" representations coming from the $p$-adic upper half plane occur as constituents of the principal series. The simplest example of this is Morita's pairing

$$(*) \qquad \Omega^1(\mathcal{X}) \times C^{\mathrm{an}}(\mathbb{P}^1(K), K)/K \to K$$

between the locally analytic functions on $\mathbb{P}^1(K)$ modulo constants (a "principal series" representation, obtained by induction from the trivial character) and the 1-forms on the one-dimensional symmetric space (a holomorphic representation.)

Morita's results illustrate how continuous representation theory extends the theory of smooth representations. Under the pairing (*), the locally constant functions on $\mathbb{P}^1(K)$ modulo constants (a smooth representation known as the Steinberg representation) are a $G$-invariant subspace which is orthogonal to the subspace of $\Omega^1(\mathcal{X})$ consisting of exact forms. In particular, this identifies the first deRham cohomology group of the $p$-adic upper half plane over $K$ with the $K$-linear dual of the Steinberg representation.

The two types of representations considered by Morita (holomorphic discrete series and principal series) have been generalized to $GL_n$.

The "holomorphic" representations defined in [Sch] use Drinfeld's $d$-dimensional $p$-adic symmetric space $\mathcal{X}$. The space $\mathcal{X}$ is the complement in $\mathbb{P}^d_{/K}$ of the $K$-rational hyperplanes. The action of the group $G := GL_{d+1}(K)$ on $\mathbb{P}^d$ preserves the missing hyperplanes, and induces an action of $G$ on $\mathcal{X}$ and a continuous action of $G$ on the infinite dimensional locally convex $K$-vector space $\mathcal{O}(\mathcal{X})$ of rigid functions on $\mathcal{X}$. The ($p$-adic) holomorphic discrete series representations are modelled on this example, and come from the action of $G$ on the global sections of homogeneous vector bundles on $\mathbb{P}^d$ restricted to $\mathcal{X}$. There is a close relationship between these holomorphic representations and classical automorphic forms, coming from the theory of $p$-adic uniformization of Shimura varieties ([RZ], [Var]).

The second type of representations we will study are the "locally analytic" representations. Such representations are developed systematically in a recent thesis of Feaux de Lacroix ([Fea]). He defines a class of representations (which he calls "weakly analytic") in locally convex vector spaces $V$ over $K$, relying on a general definition of a $V$-valued locally analytic function. Such a representation is a continuous linear action of $G$ on $V$ with the property that, for each $v$, the orbit maps $f_v(g) = g \cdot v$ are locally analytic $V$-valued functions on $G$. Notice that locally analytic representations include both smooth representations and rational ones.

Feaux de Lacroix's thesis develops some of the foundational properties of this type of representation. In particular, he establishes the basic properties of an induction functor (analytic coinduction). If we apply his induction to a one-dimensional locally analytic representation of a Borel subgroup of $G$, we obtain the $p$-adic principal series.

In this paper, we focus on one holomorphic representation and analyze it in terms of locally analytic principal series representations. Specifically, we study the representation of $G = GL_{d+1}(K)$ on the space $\Omega^d(\mathcal{X})$ of $d$-forms on the $d$-dimensional symmetric space $\mathcal{X}$. Our results generalize Morita, because we work in arbitrary dimensions, and Schneider-Stuhler, because we analyze all of $\Omega^d(\mathcal{X})$, not just its cohomology. Despite our narrow focus, we uncover new phenomena not apparent in either of the other works, and we believe that our results are representative of the general structure of holomorphic discrete series representations.

Our main results describe a $d$-step, $G$-invariant filtration on $\Omega^d(\mathcal{X})$ and a corresponding filtration on its continuous linear dual $\Omega^d(\mathcal{X})'$. We establish topological isomorphisms between the $d+1$ subquotients of the dual filtration and subquotients of members of the principal series. The $j$-th such isomorphism is given by a "boundary value map" $I^{[j]}$.

The filtration on $\Omega^d(\mathcal{X})$ comes from geometry and reflects the fact that $\mathcal{X}$ is a hyperplane complement. The first proper subspace $\Omega^d(\mathcal{X})^1$ in the filtration on $\Omega^d(\mathcal{X})$ is the space of exact forms, and the first subquotient is the $d$-th deRham cohomology group.

The principal series representation which occurs as the $j$-th subquotient of the dual of $\Omega^d(\mathcal{X})$ is a hybrid object blending rational representations, smooth representations, and differential equations. The construction of these principal series representations is a three step process. For each $j = 0, \ldots, d$, we first construct a representation $V_j$ of the maximal parabolic subgroup $P_{\underline{j}}$ of $G$ having a Levi subgroup of shape $GL_j(K) \times GL_{d+1-j}(K)$. The representation $V_j$ (which factors through this Levi subgroup) is the tensor product of a simple rational representation with the Steinberg representation of one of the Levi factors. In the second step, we apply analytic coinduction to $V_j$ to obtain a representation of $G$.

The third step is probably the most striking new aspect of our work. For each $j$, we describe a pairing between a generalized Verma module and the representation induced from $V_j$. We describe a submodule $\mathfrak{d}_{\underline{j}}$ of this Verma module such that $I^{[j]}$ is a topological isomorphism onto the subspace of the induced representation annihilated by $\mathfrak{d}_{\underline{j}}$:

$$I^{[j]} : [\Omega^d(\mathcal{X})^j/\Omega^d(\mathcal{X})^{j+1}]' \xrightarrow{\sim} C^{\mathrm{an}}(G, P_{\underline{j}}; V_j)^{\mathfrak{d}_{\underline{j}}=0}$$

The generators of the submodules $\mathfrak{d}_{\underline{j}}$ make up a system of partial differential equations. Interestingly, these differential equations are hypergeometric equations of the type studied by Gelfand and his collaborators (see [GKZ] for example). Specifically, the equations which arise here come from the adjoint action of the maximal torus of $G$ on the (transpose of) the unipotent radical of $P_{\underline{j}}$.

For the sake of comparison with earlier work, consider the two extreme cases when $j = 0$ and $j = d$. When $j = 0$, the group $P_{\underline{j}}$ is all of $G$, the representation $V_j$ is the Steinberg representation of $G$, and the induction is trivial. The submodule $\mathfrak{d}_{\underline{0}}$ is the augmentation ideal of $U(\mathfrak{g})$, which automatically kills $V_j$ because Steinberg is a smooth representation.

When $j = d$, $V_d$ is an one-dimensional rational representation of $P_{\underline{d}}$, and the module $\mathfrak{d}_{\underline{d}}$ is zero, so that there are no differential equations. In this case we obtain an isomorphism between the bottom step in the filtration and the locally analytic sections of an explicit homogeneous line bundle on the projective space $G/P_{\underline{d}}$. When $d = 1$, these two special cases ($j = 0$ and $j = 1$) together for $SL_2(K)$ are equivalent to Morita's theory applied to $\Omega^1(\mathcal{X})$.

We conclude this introduction with an outline of the sections of this paper. In sections one and two, we establish fundamental properties of $\Omega^d(\mathcal{X})$ as a topological vector space and as a $G$-representation. For example, we show that $\Omega^d(\mathcal{X})$ is a reflexive Fréchet space.

We introduce our first integral transform in section 2. Let $\xi$ be the logarithmic $d$-form on $\mathbb{P}^d$ with first order poles along the coordinate hyperplanes. We study the map

$$I : \Omega^d(\mathcal{X})' \longrightarrow C^{\mathrm{an}}(G, K)$$

$$\lambda \longmapsto [g \mapsto \lambda(g_*\xi)] .$$

We show that functions in the image of $I$ satisfy both discrete relations and differential equations, although we are unable to precisely characterize the image of the map $I$.

In section 3, we study the map $I$ in more detail. We make use of the kernel function introduced in [ST], and attempt to clarify the relationship between the transform $I$ and the results of that paper. Properties of the kernel function established in [ST], augmented by some new results, yield a map

$$I_{\mathrm{o}} : \Omega^d(\mathcal{X})' \to C(G/P, K)/C_{\mathrm{inv}}(G/P, K)$$

where $C(G/P, K)$ denotes the continuous functions on $G/P$ and $C_{\mathrm{inv}}(G/P, K)$ denotes the subspace generated by those continuous functions invariant by a larger parabolic subgroup. Using the "symmetrization map" of Borel and Serre, we show that the map $I_{\mathrm{o}}$ contains the same information as the original transform $I$. The map $I_{\mathrm{o}}$ has the advantage of targeting the possibly simpler space of functions on the compact space $G/P$. However, as was shown in [ST], the kernel function is locally analytic only on the big cell; it is continuous on all of $G/P$, but has complicated singularities at infinity. For this reason, the image of the map $I_{\mathrm{o}}$ does not lie inside the space of locally analytic functions. Introducing a notion of "analytic vectors" in a continuous representation, we prove that the image of $I_{\mathrm{o}}$ lies inside the subspace of analytic vectors, and so we can make sense of what it means for a function in the image of $I_{\mathrm{o}}$ to satisfy differential equations. However, as with $I$, we cannot completely describe the image of this "complete" integral transform, and to obtain precise results we must pass to subquotients of $\Omega^d(\mathcal{X})'$.

In the course of our analysis in section 3, we obtain the important result that the space of logarithmic forms (generated over $K$ by the $g_*\xi$) is dense in $\Omega^d(\mathcal{X})$, and consequently our maps $I$ and $I_\mathrm{o}$ are injective.

In section 4, we focus our attention on the differential equations satisfied by the functions in the image of the transform $I$. More precisely, let $\mathfrak{b}$ be the annihilator in $U(\mathfrak{g})$ of the special logarithmic form $\xi$. Any function in the image of $I$ is killed by $\mathfrak{b}$. The key result in this section is the fact that the left $U(\mathfrak{g})$-module $U(\mathfrak{g})/\mathfrak{b} = U(\mathfrak{g})\xi$ has one-dimensional weight spaces for each weight in the root lattice of $G$. In some weak sense, the $U(\mathfrak{g})$-module $U(\mathfrak{b})\xi$ plays the role of a Harish-Chandra $(\mathfrak{g}, K)$-module in our $p$-adic setting.

The filtration on $\Omega^d(\mathcal{X})$ is closely related to a descending filtration of $U(\mathfrak{g})$ by left ideals
$$U(\mathfrak{g}) = \mathfrak{b}_0 \supset \mathfrak{b}_1 \supset \ldots \supset \mathfrak{b}_{d+1} = \mathfrak{b}.$$

By combinatorial arguments using weights, we show that the subquotients of this filtration are finite direct sums of irreducible highest weight $U(\mathfrak{g})$-modules. Each of these modules has a presentation as a quotient of a generalized Verma module by a certain submodule. These submodules are the modules $\underline{\mathfrak{d}_j}$ which enter into the statement of the main theorem.

In section 5, we obtain a "local duality" result. Let $\Omega_b^d(U^0)$ be the Banach space of bounded differential forms on the admissible open set $U^0$ in $\mathcal{X}$ which is the inverse image, under the reduction map, of an open standard chamber in the Bruhat-Tits building of $G$. Let $B$ be the Iwahori group stabilizing this chamber, and let $\mathcal{O}(B)^{\mathfrak{b}=0}$ be the (globally) analytic functions on $B$ annihilated by the (left invariant) differential operators in $\mathfrak{b}$. We construct a pairing which induces a topological isomorphism between the dual space $(\mathcal{O}(B)^{\mathfrak{b}=0})'$ and $\Omega_b^d(U^0)$.

We go on in section 5 to study the filtration of $\mathcal{O}(B)^{\mathfrak{b}=0}$ whose terms are the subspaces killed by the successively larger ideals $\mathfrak{b}_i$. We compute the subquotients of this local filtration, and interpret them as spaces of functions satisfying systems of partial differential equations. These local computations are used in a crucial way in the proof of the main theorem.

In section 6, we return to global considerations and define our $G$-invariant filtration on $\Omega^d(\mathcal{X})$. We define this filtration first on the algebraic differential forms on $\mathcal{X}$. These are the rational $d$-forms having poles along an arbitrary arrangement of $K$-rational hyperplanes. The algebraic forms are dense in the rigid forms, and we define the filtration on the full space of rigid forms by taking closures. A "partial fractions" decomposition due to Gelfand-Varchenko ([GV]) plays a key role in the definition of the filtration and the proof of its main properties.

In section 7, we use rigid analysis to prove that the first step in the global filtration coincides with the space of exact forms; this implies in particular that the exact forms are closed in $\Omega^d(\mathcal{X})$. The desired results follow from a "convergent partial fractions" decomposition for global rigid forms on $\Omega^d(\mathcal{X})$. One major application of this characterization of the first stage of the filtration is that it allows us to relate the other stages with subspaces of forms coming by pull-back from lower dimensional $p$-adic symmetric spaces. Another consequence of the results of this section is an analytic proof of that part of the main theorem of [SS] describing $H_{DR}^d(\mathcal{X})$ in terms of the Steinberg representation.

In section 8, we prove the main theorem, identifying the subquotients of the filtration on the dual of $\Omega^d(\mathcal{X})$ with the subspaces of induced representations killed by the correct differential operators. All of the prior results are brought to bear on the problem. We show that the integral transform is bijective by showing that an element of the induced representation satisfying the differential equations can be written as a finite sum of $G$-translates of elements of a very special form, and then explicitly exhibiting an inverse image of such a special element. The fact that the map is a topological isomorphism follows from continuity and a careful application of an open-mapping theorem.

Part of this work was presented in a course at the Institut Henri Poincaré during the "$p$-adic semester" in 1997 . We are very grateful for this opportunity as well as for the stimulating atmosphere during this activity. The second author was supported by grants from the National Science Foundation.## 0. Notations and conventions

For the reader's convenience, we will begin by summarizing some of the notation we use in this paper. In general, we have followed the notational conventions of [ST].

Let $K$ denote a fixed, non-archimedean locally compact field of characteristic zero, residue characteristic $p > 0$ and ring of integers $o$. Let $|\cdot|$ be the absolute value on $K$, let $\omega : K \to \mathbb{Z}$ be the normalized additive valuation, and let $\pi$ be a uniformizing parameter. We will use $\mathbb{C}_p$ for the completion of an algebraic closure of $K$.
Fix an integer $d \geq 1$ and let $\mathbb{P}^d$ be the projective space over $K$ of dimension $d$. We let $G := GL_{d+1}(K)$, and adopt the convention that $G$ acts on $\mathbb{P}^d$ through the left action $g([q_0 : \cdots : q_d]) = [q_0 : \cdots : q_d]g^{-1}$. We let $T$ be the diagonal torus in $G$, and $\overline{T}$ the image of $T$ in $PGL_{d+1}(K)$. We use $\epsilon_0, \ldots, \epsilon_d$ for the characters of $T$, where, if $t = (t_{ii})_{i=0}^d$ is a diagonal matrix, then $\epsilon_i(t) = t_{ii}$. The character group $X^*(\overline{T})$ is the root lattice of $G$. It is spanned by the set $\Phi := \{\epsilon_i - \epsilon_j : 0 \leq i \neq j \leq d\}$ of roots of $G$. Let $\Xi_0, \ldots, \Xi_d$ be homogeneous coordinates for $\mathbb{P}^d$. Suppose that $\mu \in X^*(\overline{T})$, and write $\mu = \sum_{i=0}^d m_i \epsilon_i$. We let

$$\Xi_\mu = \prod_{i=0}^d \Xi_i^{m_i}.$$

Since $\mu$ belongs to the root lattice, we know that $\sum_{i=0}^{d} m_i = 0$, and therefore $\Xi_\mu$ is a well-defined rational function on $\mathbb{P}^d$. Certain choices of $\mu$ arise frequently and so we give them special names. For $i = 0, \ldots, d-1$ we let $\beta_i = \epsilon_i - \epsilon_d$ and $\beta = \beta_0 + \ldots + \beta_{d-1}$. We also let $\alpha_i = \epsilon_{i+1} - \epsilon_i$, for $i = 0, \ldots, d-1$. The set $\{\alpha_i\}_{i=0}^{d-1}$ is a set of simple roots. We also adopt the convention that $\alpha_d = \epsilon_0 - \epsilon_d$. Any weight $\mu$ in $X^*(\overline{T})$ may be written uniquely as a sum $\mu = \sum_{i=0}^{d} m_i \alpha_i$ with integers $m_i \geq 0$ of which at least one is equal to 0. If $\mu$ is written in this way, we let $\ell(\mu) := m_d$.

As mentioned in the introduction, we let $\mathcal{X}$ denote Drinfeld's $d$-dimensional $p$-adic symmetric space. The space $\mathcal{X}$ is the complement in $\mathbb{P}^d$ of the $K$-rational hyperplanes. The $G$-action on $\mathbb{P}^d$ preserves $\mathcal{X}$. The structure of $\mathcal{X}$ as a rigid analytic space comes from an admissible covering of $\mathcal{X}$ by an increasing family of open $K$-affinoid subvarieties $\mathcal{X}_n$. To define the subdomains $\mathcal{X}_n$, let $\mathcal{H}$ denote the set of hyperplanes in $\mathbb{P}^d$ which are defined over $K$. For any $H \in \mathcal{H}$ let $\ell_H$ be a unimodular linear form in $\Xi_0, \ldots, \Xi_d$ such that $H$ is the zero set of $\ell_H$. (Here, and throughout this paper, a linear form $\ell_H$ is called unimodular if it has coefficients in $o$ and at least one coefficient is a unit.) The set $\mathcal{X}_n$ consists of the set of points $q \in \mathbb{P}^d$ such that

$$\omega(\ell_H([q_0 : \cdots : q_d])) \leq n$$

for any $H \in \mathcal{H}$ whenever $[q_0 : q_1 : \cdots : q_d]$ is a unimodular representative for the homogeneous coordinates of $q$. We denote by $\mathcal{O}(\mathcal{X})$ the ring of global rigid analytic functions on $\mathcal{X}$, and by $\Omega^i(\mathcal{X})$ the global $i$-forms. By $H^*_{DR}(\mathcal{X})$ we mean the rigid-analytic deRham cohomology of $\mathcal{X}$.

The space $\mathcal{X}$ has a natural $G$-equivariant map (the reduction map) $r : \mathcal{X} \to \overline{X}$ to the Bruhat-Tits building $\overline{X}$ of $PGL_{d+1}(K)$. For the definition of this map, see Definition 2 of [ST].

The torus $\overline{T}$ stabilizes a standard apartment $\overline{A}$ in $\overline{X}$. The Iwahori group

$$B := \{g \in GL_{d+1}(o) : g \text{ is lower triangular mod } \pi\}$$

is the pointwise stabilizer of a certain closed chamber $\overline{C}$ in $\overline{A} \subset X$. Following the conventions of [ST], we mean by $(\overline{C}, 0)$ the chamber $\overline{C}$ together with the vertex 0 stabilized by $GL_{d+1}(o)$. We will frequently denote a random closed chamber in $X$ with the letter $\Delta$, while $\Delta^0$ will denote the interior of $\Delta$. The inverse image $U^0 = r^{-1}(\overline{C}^0)$ of the open standard chamber $\overline{C}^0$ under the reduction map is an admissible open subset in $\mathcal{X}$.

In addition to these conventions regarding roots and weights of $G$, we use the following letters for various objects associated with $G$:

$$\begin{array}{rcl} P & := & \text{the lower triangular Borel subgroup of } G \\ U & := & \text{the lower triangular unipotent group of } G \\ N & := & \text{the normalizer of } T \text{ in } G \\ W & := & \text{the Weyl group of } G \text{ (associated to diagonal torus } T) \\ w_{d+1} & := & \text{the longest element in } W \\ P_s & := & P \cup PsP \text{ for any simple reflection } s \in W \end{array}$$

For an element $g \in Uw_{d+1}P$ in the big cell we define $u_g \in U$ by the identity $g = u_g w_{d+1} h$ with $h \in P$.

Corresponding to a root $\alpha = \epsilon_i - \epsilon_j$ we have a homomorphism $\tilde{\alpha} : K^+ \to G$ sending $u \in K^+$ to the matrix $(u_{rs})$ with:

$$u_{rs} = \begin{cases} 1 & \text{if } r = s \\ u & \text{if } r = i \text{ and } s = j \\ 0 & \text{otherwise.} \end{cases}$$

The image $U_\alpha$ of $\tilde{\alpha}$ in $G$ is the root subgroup associated to $\alpha$. It is filtered by the subgroups $U_{\alpha,r} := \tilde{\alpha}(\{u \in K : \omega(u) \geq r\})$ for $r \in \mathbb{R}$. For a point $x \in \overline{A}$ we define $U_x$ to be the subgroup of $G$ generated by all $U_{\alpha, -\alpha(x)}$ for $\alpha \in \Phi$.

## 1. $\Omega^d(\mathcal{X})$ as a locally convex vector space

We begin by establishing two fundamental topological properties of $\Omega^d(\mathcal{X})$. We construct a family of norms on $\Omega^d(\mathcal{X})$, parameterized by chambers of the building $\overline{X}$, which defines the natural Fréchet topology (coming from its structure as a projective limit of Banach spaces) on $\Omega^d(\mathcal{X})$. We further show the fundamental result that $\Omega^d(\mathcal{X})$ is a reflexive Fréchet space.

We first look at the space $\mathcal{O}(\mathcal{X})$. For any open $K$-affinoid subvariety $\mathcal{Y} \subseteq \mathcal{X}$ its ring $\mathcal{O}(\mathcal{Y})$ of analytic functions is a $K$-Banach algebra with respect to the spectral norm. We equip $\mathcal{O}(\mathcal{X})$ with the initial topology with respect to the family of restriction maps $\mathcal{O}(\mathcal{X}) \longrightarrow \mathcal{O}(\mathcal{Y})$. Since the increasing family of open $K$-affinoid subvarieties $\mathcal{X}_n$ forms an admissible covering of $\mathcal{X}$ ([SS] Sect. 1) we have

$$\mathcal{O}(\mathcal{X}) = \varprojlim_n \mathcal{O}(\mathcal{X}_n)$$

in the sense of locally convex $K$-vector spaces. It follows in particular that $\mathcal{O}(\mathcal{X})$ is a Fréchet space. Using a basis $\eta_0$ of the free $\mathcal{O}(\mathcal{X})$-module $\Omega^d(\mathcal{X})$ of rank 1 we topologize $\Omega^d(\mathcal{X})$ by declaring the linear map

$$\mathcal{O}(\mathcal{X}) \xrightarrow{\sim} \Omega^d(\mathcal{X})$$

$$F \longmapsto F\eta_0$$

to be a topological isomorphism; the resulting topology is independent of the choice of $\eta_0$. In this way $\Omega^d(\mathcal{X})$ becomes a Fréchet space, too. Similarly each $\Omega^d(\mathcal{X}_n)$ becomes a Banach space. In the following we need a certain $G$-invariant family of continuous norms on $\Omega^d(\mathcal{X})$. First recall the definition of the weights

$$\beta_i := \varepsilon_i - \varepsilon_d \quad \text{for} \quad 0 \leq i \leq d-1 \ .$$

We have
$$\Omega^d(\mathcal{X}) = \mathcal{O}(\mathcal{X}) d\Xi_{\beta_0} \wedge \ldots \wedge d\Xi_{\beta_{d-1}} .$$
The torus $\overline{T}$ acts on the form $d\Xi_{\beta_0} \wedge \ldots \wedge d\Xi_{\beta_{d-1}}$ through the weight
$$\beta := \beta_0 + \ldots + \beta_{d-1} .$$
For any point $q \in \mathcal{X}$ such that $z := r(q) \in \overline{A}$ we define a continuous (additive) semi-norm $\gamma_q$ on $\Omega^d(\mathcal{X})$ by
$$\gamma_q(\eta) := \omega(F(q)) + \beta(z) \quad \text{if} \quad \eta = F d\Xi_{\beta_0} \wedge \ldots \wedge d\Xi_{\beta_{d-1}} .$$

**Lemma 1:**

Let $q \in \mathcal{X}$ such that $\overline{x} := r(q) \in \overline{A}$; we then have
$$\gamma_{gq} = \gamma_q \circ g^{-1} \quad \text{for any} \quad g \in N \cup U_x .$$

Proof: First let $g \in G$ be any element such that $g\overline{x} \in \overline{A}$. Using [ST] Cor. 4 and the characterizing property of the function $\mu(g^{-1}, .)$ ([ST] Def. 28) one easily computes
$$\gamma_{gq} - \gamma_q \circ g^{-1} = \omega\left( \frac{g_*^{-1}\Xi_0}{\Xi_0}(q) \cdot \ldots \cdot \frac{g_*^{-1}\Xi_d}{\Xi_d}(q) \right) + \omega(\det g) .$$
Obviously the right hand side vanishes if $g$ is a diagonal or permutation matrix and hence for any $g \in N$. It remains to consider a $g = \tilde{\alpha}(u) \in U_{\alpha, -\alpha(x)}$ for some root $\alpha \in \Phi$. Then the right hand side simplifies to $\omega(1 - u\Xi_\alpha(q)) = \omega(\Xi_\alpha(gq)) - \omega(\Xi_\alpha(q))$. According to [ST] Cor. 4 this is equal to $\alpha(r(gq)) - \alpha(r(q)) = \alpha(\overline{x}) - \alpha(\overline{x}) = 0$. □

This allows us to define, for any point $q \in \mathcal{X}$, a continuous semi-norm $\gamma_q$ on $\Omega^d(\mathcal{X})$ by
$$\gamma_q := \gamma_{gq} \circ g$$
where $g \in G$ is chosen in such a way that $r(gq) \in \overline{A}$. Moreover, for any chamber $\Delta$ in $\overline{X}$, we put
$$\gamma_\Delta := \inf_{r(q) \in \Delta^0} \gamma_q .$$
Since $r^{-1}(\Delta)$ is an affinoid ([ST] Prop. 13) this is a continuous semi-norm. To see that it actually is a norm let us look at the case of the standard chamber $\overline{C}$. Let $\eta = F \cdot d\Xi_{\alpha_{d-1}} \wedge \ldots \wedge d\Xi_{\alpha_0} \in \Omega^d(\mathcal{X})$. Since $F|U^0$ is bounded we have the expansion
$$F|U^0 = \sum_{\mu \in X^*(\overline{T})} a(\mu)\Xi_\mu$$

with $a(\mu) \in K$ and $\{\omega(a(\mu)) - l(\mu)\}_\mu$ bounded below. Since the restriction map $\Omega^d(\mathcal{X}) \longrightarrow \Omega^d(U^0)$ is injective we have the norm

$$\omega_C(\eta) := \inf_\mu \{\omega(a(\mu)) - l(\mu)\} = \inf_{q \in U^0} \omega(F(q))$$

on $\Omega^d(\mathcal{X})$.

**Lemma 2:**

$\omega_C \leq \gamma_{\overline{C}} \leq \omega_C + 1$.

Proof: Let $\eta := F \cdot d\Xi_{\alpha_{d-1}} \wedge \ldots \wedge d\Xi_{\alpha_0}$. The identity

$$d\Xi_{\alpha_{d-1}} \wedge \ldots \wedge d\Xi_{\alpha_0} = \pm \Xi_{-\beta-\alpha_d} d\Xi_{\beta_0} \wedge \ldots \wedge d\Xi_{\beta_{d-1}}$$

together with [ST] Cor.4 implies

$$\gamma_q(\eta) = \omega(F(q)) + \omega(\Xi_{-\beta-\alpha_d}(q)) + \beta(z) = \omega(F(q)) - \alpha_d(z)$$

for $r(q) = z \in \overline{C}^0$. Because of $-1 \leq \alpha_d|\overline{C} \leq 0$ we obtain

$$\omega(F(q)) \leq \gamma_q(\eta) \leq \omega(F(q)) + 1$$

for any $q \in U^0$. It remains to recall that $\omega_C(\eta) = \inf_{q \in U^0} \omega(F(q))$. □

This shows that all the $\gamma_\Delta$ are continuous norms on $\Omega^d(\mathcal{X})$. In fact the family of norms $\{\gamma_\Delta\}_\Delta$ defines the Fréchet topology of $\Omega^d(\mathcal{X})$. In order to see this it suffices to check that the additively written spectral norm $\omega_\Delta$ for the affinoid $r^{-1}(\Delta)$ satisfies

$$\omega_\Delta(F) = \inf_{r(q) \in \Delta^0} \omega(F(q)) \quad \text{for } F \in \mathcal{O}(\mathcal{X}) .$$

Let $\mathcal{X}_B$ denote Berkovich's version of the rigid analytic variety $\mathcal{X}$. Each point $q \in \mathcal{X}_B$ gives rise to the multiplicative semi-norm $F \mapsto \omega(F(q))$ on $\mathcal{O}(\mathcal{X})$. If one fixes $F \in \mathcal{O}(\mathcal{X})$ then the function $q \mapsto \omega(F(q))$ is continuous on $\mathcal{X}_B$. We need the following facts from [Be2]:

- The reduction map $r : \mathcal{X} \longrightarrow \overline{X}$ extends naturally to a continuous map $r_B : \mathcal{X}_B \longrightarrow \overline{X}$.

- The map $r_B$ has a natural continuous section $s_B : \overline{X} \longrightarrow \mathcal{X}_B$ such that $F \longmapsto \omega(F(s_B(z)))$, for $z \in r(\mathcal{X})$, is the spectral norm $\omega_{r^{-1}(z)}$ for the affinoid $r^{-1}(z)$.

In particular, for a fixed $F \in \mathcal{O}(\mathcal{X})$, the map $z \longmapsto \omega(F(s_B(z)))$ is continuous on $\overline{X}$. Since $r(\mathcal{X})$ is dense in $\overline{X}$ it follows that

$$\inf_{r(q) \in \Delta^0} \omega(F(q)) = \inf_{z \in r(\mathcal{X}) \cap \Delta^0} \omega_{r^{-1}(z)}(F) = \inf_{z \in r(\mathcal{X}) \cap \Delta} \omega_{r^{-1}(z)}(F)$$
$$= \inf_{r(q) \in \Delta} \omega(F(q)) = \omega_\Delta(F) .$$

**Lemma 3:**

*The $G$-action $G \times \Omega^d(\mathcal{X}) \longrightarrow \Omega^d(\mathcal{X})$ is continuous.*

Proof: Clearly each individual element $g \in G$ induces a continuous automorphism of $\Omega^d(\mathcal{X})$. As a Fréchet space $\Omega^d(\mathcal{X})$ is barrelled ([Tie] Thm. 3.15). Hence the Banach-Steinhaus theorem ([Tie] Thm. 4.1) holds for $\Omega^d(\mathcal{X})$ and we only have to check that the maps

$$\begin{aligned} G &\longrightarrow \Omega^d(\mathcal{X}) \quad \text{for } \eta \in \Omega^d(\mathcal{X}) \\ g &\longmapsto g\eta \end{aligned}$$

are continuous (compare the reasoning in [War] p. 219). By the universal property of the projective limit topology this is a consequence of the much stronger local analyticity property which we will establish in Prop. 1' of the next section.

**Proposition 4:**

*$\mathcal{O}(\mathcal{X})$ is reflexive and its strong dual $\mathcal{O}(\mathcal{X})'$ is the locally convex inductive limit*

$$\mathcal{O}(\mathcal{X})' = \varinjlim_n \mathcal{O}(\mathcal{X}_n)'$$

*of the dual Banach spaces $\mathcal{O}(\mathcal{X}_n)'$.*

The proof is based on the following concepts.

**Definition:**

*A homomorphism $\psi : \mathcal{A} \longrightarrow \mathcal{B}$ between $K$-Banach spaces is called compact if the image under $\psi$ of the unit ball $\{f \in \mathcal{A} : |f|_\mathcal{A} \leq 1\}$ in $\mathcal{A}$ is relatively compact in $\mathcal{B}$.*

We want to give a general criterion for a homomorphism of affinoid $K$-algebras to be compact. Recall that an affinoid $K$-algebra $\mathcal{A}$ is a Banach algebra with respect to the residue norm $|\ |_a$ induced by a presentation

$$a : K\langle T_1, \ldots, T_m \rangle \twoheadrightarrow \mathcal{A}$$

as a quotient of a Tate algebra. All these norms $|\ |_a$ are equivalent.

**Definition:** ([Ber] 2.5.1)

*A homomorphism $\psi : \mathcal{A} \longrightarrow \mathcal{B}$ of affinoid $K$-algebras is called inner if there is a presentation $a : K\langle T_1, \ldots, T_m \rangle \twoheadrightarrow \mathcal{A}$ such that*

$$\inf\{\omega(\psi a(T_i)(y)) : y \in \mathrm{Sp}(\mathcal{B}),\ 1 \leq i \leq m\} > 0 \ .$$

**Lemma 5:**

*Any inner homomorphism $\psi : \mathcal{A} \longrightarrow \mathcal{B}$ of affinoid $K$-algebras is compact.*

Proof: First of all we note that if the assertion holds for one residue norm on $\mathcal{A}$ then it holds for all of them. If $\psi$ is inner we find, according to [Ber] 2.5.2, a presentation $a : K\langle T_1, \ldots, T_m \rangle \twoheadrightarrow \mathcal{A}$ such that

$$\inf\{\omega(\psi a(T_i)(y)) \ :\ y \in \mathrm{Sp}(\mathcal{B}),\ 1 \leq i \leq m\} > 1 \ .$$

This means that we actually have a commutative diagram of affinoid $K$-algebras

$$\begin{array}{ccc} K\langle T_1, \ldots, T_m \rangle & \xrightarrow{\imath} & K\langle \pi^{-1}T_1, \ldots, \pi^{-1}T_m \rangle \\ a \downarrow & & \downarrow \\ \mathcal{A} & \xrightarrow{\psi} & \mathcal{B} \end{array}$$

where $\imath$ is the obvious inclusion of Tate algebras. Since the valuation of $K$ is discrete the unit ball in $K\langle T_1, \ldots, T_m \rangle$ (with respect to the Gauss norm) is mapped surjectively, by $a$, onto the unit ball in $\mathcal{A}$ (with respect to $|\ |_a$). Hence it suffices to prove that the inner monomorphism $\imath$ is compact. But this is a straightforward generalization of the argument in the proof of [Mo1] 3.5.

Proof of Prop. 4: In the proof of [SS] §1 Prop. 4 the following two facts are established:

- The restriction maps $\mathcal{O}(\mathcal{X}_{n+1}) \longrightarrow \mathcal{O}(\mathcal{X}_n)$ are inner;
- $\mathcal{X}_n$ is a Weierstraß domain in $\mathcal{X}_{n+1}$ for each $n$.

The second fact implies that the restriction map $\mathcal{O}(\mathcal{X}_{n+1}) \longrightarrow \mathcal{O}(\mathcal{X}_n)$ has a dense image. It then follows from Mittag-Leffler ([B-TG3] II §3.5 Thm. 1) that the restriction maps $\mathcal{O}(\mathcal{X}) \longrightarrow \mathcal{O}(\mathcal{X}_n)$ have dense images. Using Lemma 5 we see that the assumptions in [Mo1] 3.3(i) and 3.4(i) are satisfied for the sequence of Banach spaces $\mathcal{O}(\mathcal{X}_n)$. Our assertion results. □

Of course then also $\Omega^d(\mathcal{X})$ is reflexive with $\Omega^d(\mathcal{X})' = \varinjlim_n \Omega^d(\mathcal{X}_n)'$.

## 2. $\Omega^d(\mathcal{X})$ as a locally analytic $G$-representation

In this section, we study the $G$-action on $\Omega^d(\mathcal{X})$ and investigate in which sense it is locally analytic. Using this property of the $G$-action, we construct a continuous map $I$ from $\Omega^d(\mathcal{X})'$ to the space of locally analytic $K$-valued functions on $G$. It follows from the construction of this map that its image consists of functions annihilated by a certain ideal $\mathfrak{a}$ in the algebra of punctual distributions on $G$. In particular, this means that functions in the image of $I$ satisfy both discrete relations (meaning that their values at certain related points of $G$ cannot be independently specified) and differential equations. We will study these relations in more detail in later sections.

We will use the notion of a locally analytic map from a locally $K$-analytic manifold into a Hausdorff locally convex $K$-vector space as it is defined in [B-VAR] 5.3.1. But we add the attribute "locally" in order to make clearer the distinction from rigid analytic objects.

**Proposition 1:**

*For any $\eta \in \Omega^d(\mathcal{X})$ and any $\lambda \in \Omega^d(\mathcal{X})'$ the function $g \longmapsto \lambda(g_*\eta)$ on $G$ is locally analytic.*

Since, by Prop. 1.4, $\lambda$ comes from a continuous linear form on some $\Omega^d(\mathcal{X}_n)$ this is an immediate consequence of the following apparently stronger fact.

**Proposition 1':**

*Whenever $\Omega^d(\mathcal{X})$ is equipped with the coarser topology coming from the spectral norm on $\mathcal{X}_n$ for some fixed but arbitrary $n \in \mathbb{N}$ then the map $g \mapsto g_*\eta$, for any $\eta \in \Omega^d(\mathcal{X})$, is locally analytic.*

Proof: For the moment being we fix a natural number $n \in \mathbb{N}$. In the algebraic, and hence rigid analytic, $K$-group $GL_{d+1}$ we have the open $K$-affinoid subgroup

$$H_n := \{h \in GL_{d+1}(o_{\mathbb{C}_p}) : h \equiv g \mod \pi^{n+1} \text{ for some } g \in GL_{d+1}(o)\}$$

which contains the open $K$-affinoid subgroup

$$D_n := 1 + \pi^{n+1} M_{d+1}(o_{\mathbb{C}_p}) \ ;$$

here $o$, resp. $o_{\mathbb{C}_p}$, denotes the ring of integers in $K$, resp. $\mathbb{C}_p$. As a rigid variety over $K$ the latter group $D_n$ is a polydisk of dimension $r := (d+1)^2$. Since $H_n$ preserves the $K$-affinoid subdomain $\mathcal{X}_n$ of $\mathbb{P}^d$ the algebraic action of $GL_{d+1}$ on $\mathbb{P}^d$ restricts to a rigid analytic action $m : H_n \times \mathcal{X}_n \longrightarrow \mathcal{X}_n$ which corresponds to a homomorphism of $K$-affinoid algebras

$$\begin{aligned}\mathcal{O}(\mathcal{X}_n) &\longrightarrow \mathcal{O}(H_n \times \mathcal{X}_n) = \mathcal{O}(H_n) \hat{\otimes}_K \mathcal{O}(\mathcal{X}_n) \\ F &\longmapsto m^*F\ .\end{aligned}$$

For any $h \in H_n$ we clearly have

$$[(\text{evaluation in } h) \otimes id] \circ m^*F = hF\ .$$

For a fixed $g \in GL_{d+1}(o)$ we consider the rigid analytic "chart"

$$\begin{aligned}\imath_g : D_n &\longrightarrow H_n \\ h &\longmapsto gh\ .\end{aligned}$$

Fixing coordinates $T_1, \ldots, T_r$ on the polydisk $D_n$ we have $\mathcal{O}(D_n) \hat{\otimes}_K \mathcal{O}(\mathcal{X}_n) \cong \mathcal{O}(\mathcal{X}_n)\langle T_1, \ldots, T_r\rangle$. The power series

$$\mathcal{F}_g(T_1, \ldots, T_r) := (\imath_g^* \otimes id) m^*F \in \mathcal{O}(\mathcal{X}_n)\langle T_1, \ldots, T_r\rangle$$

has the property that $ghF = \mathcal{F}_g(T_1(h), \ldots, T_r(h))$ for any $h \in D_n$. This shows that, for any $F \in \mathcal{O}(\mathcal{X}_n)$, the map

$$\begin{aligned}GL_{d+1}(o) &\longrightarrow \mathcal{O}(\mathcal{X}_n) \\ g &\longmapsto gF\end{aligned}$$

is locally analytic.

This construction varies in an obvious way with the natural number $n$. In particular if we start with a function $F \in \mathcal{O}(\mathcal{X}) \subseteq \mathcal{O}(\mathcal{X}_n)$ then the coefficients of the power series $\mathcal{F}_g$ also lie in $\mathcal{O}(\mathcal{X})$. It follows that actually, for any $F \in \mathcal{O}(\mathcal{X})$, the map

$$\begin{aligned}GL_{d+1}(o) &\longrightarrow \mathcal{O}(\mathcal{X}) \\ g &\longmapsto gF\end{aligned}$$

is locally analytic provided the right hand side is equipped with the sup-norm on $\mathcal{X}_n$ for a fixed but arbitrary $n \in \mathbb{N}$. Since $F$ was arbitrary and $GL_{d+1}(o)$ is open in $G$ the full map

$$\begin{aligned}G &\longrightarrow \mathcal{O}(\mathcal{X}) \\ g &\longmapsto gF\end{aligned}$$

has to have the same local analyticity property.

This kind of reasoning extends readily to any $GL_{d+1}$-equivariant algebraic vector bundle $\mathbb{V}$ on $\mathbb{P}^d$. Then the space of rigid analytic sections $\mathbb{V}(\mathcal{X})$ is a Fréchet space as before on which $G$ acts continuously and such that the maps

$$\begin{aligned}G &\longrightarrow \mathbb{V}(\mathcal{X}) \\ g &\longmapsto gs\end{aligned}$$

for any $s \in \mathbb{W}(\mathcal{X})$ have the analogous local analyticity property. The reason is that the algebraic action induces a rigid analytic action

$$H_n \times \mathbb{W}_{/\mathcal{X}_n} \longrightarrow \mathbb{W}_{/\mathcal{X}_n}$$

which is compatible with the action of $H_n$ on $\mathcal{X}_n$ via $m$. But this amounts to the existence of a vector bundle isomorphism

$$m^*(\mathbb{W}_{/\mathcal{X}_n}) \xrightarrow{\cong} pr_2^*(\mathbb{W}_{/\mathcal{X}_n})$$

satisfying a certain cocycle condition (compare [Mum] 1.3); here $pr_2 : H_n \times \mathcal{X}_n \longrightarrow \mathcal{X}_n$ is the projection map. Hence similarly as above the $H_n$-action on the sections $\mathbb{W}(\mathcal{X}_n)$ is given by a homomorphism

$$\mathbb{W}(\mathcal{X}_n) \longrightarrow m^*(\mathbb{W}_{/\mathcal{X}_n})(H_n \times \mathcal{X}_n) \xrightarrow{\cong} pr_2^*(\mathbb{W}_{/\mathcal{X}_n})(H_n \times \mathcal{X}_n) = \mathcal{O}(\mathcal{X}_n) \hat{\otimes}_K \mathbb{W}(\mathcal{X}_n) .$$

The rest of the argument then is the same as above. □

That result has two important consequences for our further investigation. In the first place it allows us to introduce the basic map for our computation of the dual space $\Omega^d(\mathcal{X})'$. Let

$$C^{\mathrm{an}}(G, K) := \text{ space of locally } K\text{-analytic functions on } G .$$

We always consider this space as the locally convex inductive limit

$$C^{\mathrm{an}}(G, K) = \varinjlim_{\mathcal{U}} C^{\mathrm{an}}_{\mathcal{U}}(G, K) .$$

Here $\mathcal{U} = \{U_i\}_{i \in I}$ is a disjoint covering of the locally $K$-analytic manifold $G$ by closed balls (in the sense of charts) and

$$C^{\mathrm{an}}_{\mathcal{U}}(G, K) := \{f \in C^{\mathrm{an}}(G, K) : f|U_i \text{ is analytic for any } i \in I\}$$

is the direct product of the Banach spaces of analytic functions on each $U_i$ (where the Banach norm is the spectral norm on $U_i$). The group $G$ acts by left translations on $C^{\mathrm{an}}(G, K)$.

**Lemma 2:**

*The $G$-action $G \times C^{\mathrm{an}}(G, K) \longrightarrow C^{\mathrm{an}}(G, K)$ is continuous.*

Proof: Clearly each individual group element $g \in G$ acts continuously on $C^{\mathrm{an}}(G, K)$. Being the locally convex inductive limit of a direct product of Banach spaces $C^{\mathrm{an}}(G,K)$ is barrelled. Hence it suffices (as in the proof of Lemma 1.3) to check that the maps

$$\begin{aligned} G &\longrightarrow C^{\mathrm{an}}(G,K) \quad \text{for } f \in C^{\mathrm{an}}(G,K) \\ g &\longmapsto gf \end{aligned}$$

are continuous. But those maps actually are differentiable ([Fea] 3.3.4). $\square$

In all that follows, the $d$-form

$$\xi := \frac{d\Xi_{\beta_0} \wedge \ldots \wedge d\Xi_{\beta_{d-1}}}{\Xi_{\beta_0} \cdot \ldots \cdot \Xi_{\beta_{d-1}}}$$

on $\mathcal{X}$ is the basic object. Because of Prop. 1 we have the $G$-equivariant map

$$\begin{aligned} I : \Omega^d(\mathcal{X})' &\longrightarrow C^{\mathrm{an}}(G,K) \\ \lambda &\longmapsto [g \mapsto \lambda(g_*\xi)] \ . \end{aligned}$$

**Lemma 3:**

*The map $I$ is continuous.*

Proof: Since $\Omega^d(\mathcal{X})'$ is the locally convex inductive limit of the Banach spaces $\Omega^d(\mathcal{X}_n)'$ it suffices to establish the corresponding fact for $\Omega^d(\mathcal{X}_n)$. In the proof of Prop. 1' we have seen that the map

$$\begin{aligned} G &\longrightarrow \Omega^d(\mathcal{X}_n) \\ g &\longmapsto g_*\xi \end{aligned}$$

is analytic on the right cosets of $G \cap D_n$ in $G$. We obtain that, for $\lambda \in \Omega^d(\mathcal{X}_n)'$, the function $g \longmapsto \lambda(g_*\xi)$ lies in $C^{\mathrm{an}}_{\mathcal{U}}(G,K)$ with $\mathcal{U} := \{(G\cap D_n)g\}_{g\in G}$ and that on a fixed coset $(G\cap D_n)g$ the spectral norms satisfy the inequality

$$\|\lambda(._*\xi)\| \leq \|\lambda\| \cdot \|._*\xi\| \ . \qquad \square$$

We also have the right translation action of $G$ on $C^{\mathrm{an}}(G,K)$ which we write as

$$\delta_g f(h) \ = \ f(hg) \ .$$

In addition we have the action of the Lie algebra $\mathfrak{g}$ of $G$ by left invariant differential operators; for any $\mathfrak{x} \in \mathfrak{g}$ the corresponding operator on $C^{\mathrm{an}}(G,K)$ is given by

$$(\mathfrak{x}f)(g) := \frac{d}{dt} f(g\exp(t\mathfrak{x}))_{|t=0} \ .$$

here exp : $\mathfrak{g} --- > G$ denotes the exponential map which is defined locally around 0. This extends by the universal property to a left action of the universal enveloping algebra $U(\mathfrak{g})$ on $C^{\mathrm{an}}(G, K)$. For any $f \in C^{\mathrm{an}}(G, K)$, any $g \in G$, and any $\mathfrak{x} \in \mathfrak{g}$ sufficiently close to 0 (depending on $g$) we have Taylor's formula

$$f(g \exp(\mathfrak{x})) = \sum_{n=0}^{\infty} \frac{1}{n!} (\mathfrak{x}^n f)(g)$$

(compare, for example, the proof in [Hel] II.1.4 which goes through word for word for $p$-adic Lie groups). We actually find for any $h \in G$ a neighbourhood $N_0$ of $h$ in $G$ and a neighbourhood $\mathfrak{n}$ of 0 in $\mathfrak{g}$ such that the above formula holds for all $(g, \mathfrak{x}) \in N_0 \times \mathfrak{n}$.

The right translation action of $G$ and the $U(\mathfrak{g})$-action on $C^{\mathrm{an}}(G, K)$ combine into an action of the algebra $\mathcal{D}(G)$ of punctual distributions on $G$ ([B-GAL] III §3.1). Any $D \in \mathcal{D}(G)$ can be written in a unique way as a finite sum $D = \mathfrak{z}_1 \delta_{g_1} + \ldots + \mathfrak{z}_r \delta_{g_r}$ with $\mathfrak{z}_i \in U(\mathfrak{g})$ and $g_i \in G$, $\delta_g$ denoting the Dirac distribution supported at $g \in G$. Then one has $Df = \sum_i \mathfrak{z}_i(f(.g_i))$ for $f \in C^{\mathrm{an}}(G, K)$; observe that

$$\delta_g(\mathfrak{z}(f)) = (\mathrm{ad}(g)\mathfrak{z})(\delta_g(f)) .$$

This $\mathcal{D}(G)$-action commutes with the left translation action of $G$ on $C^{\mathrm{an}}(G, K)$. Moreover $\mathcal{D}(G)$ acts by continuous endomorphisms on $C^{\mathrm{an}}(G, K)$; this is again a simple application of the Banach-Steinhaus theorem (compare [Fea] 3.1.2).
The second consequence of Prop. 1' is that the map $g \mapsto g_* \eta$ from $G$ into $\Omega^d(\mathcal{X})$ is differentiable ([B-VAR] 1.1.2) for any $\eta \in \Omega^d(\mathcal{X})$. It follows that $\mathfrak{g}$ and hence $U(\mathfrak{g})$ act on $\Omega^d(\mathcal{X})$ from the left by

$$\mathfrak{x}\eta := \frac{d}{dt} \exp(t\mathfrak{x})_* \eta_{|t=0} .$$

Obviously the $G$-action and the $U(\mathfrak{g})$-action again combine into a left $\mathcal{D}(G)$-action by continuous endomorphisms on $\Omega^d(\mathcal{X})$. Note that $\Omega^d(\mathcal{X})$ as a Fréchet space is barrelled, too. We define now

$$\mathfrak{a} := \{D \in \mathcal{D}(G) : D\xi = 0\}$$

to be the annihilator ideal of $\xi$ in $\mathcal{D}(G)$; it is a left ideal. On the other hand

$$C^{\mathrm{an}}(G, K)^{\mathfrak{a}=0} := \{f \in C^{\mathrm{an}}(G, K) : \mathfrak{a}f = 0\}$$

then is a $G$-invariant closed subspace of $C^{\mathrm{an}}(G, K)$. The formula

$$[D(I(\lambda))](g) = \lambda(g_*(D\xi)) \quad \text{for} \quad D \in \mathcal{D}(G), \lambda \in \Omega^d(\mathcal{X})' \quad \text{and} \quad g \in G$$

implies that that subspace contains the image of the map $I$, i.e., that $I$ induces a $G$-equivariant continuous linear map

$$\Omega^d(\mathcal{X})' \longrightarrow C^{\mathrm{an}}(G, K)^{\mathfrak{a}=0} .$$

## 3. The kernel map

In the previous section we constructed a map $I$ from $\Omega^d(\mathcal{X})'$ to a certain space of locally analytic functions on $G$. We see this map as a "boundary value" map, but this interpretation needs clarification. In particular, the results of [SS] and [ST] suggest that a more natural "boundary" for the symmetric space $\mathcal{X}$ is the compact space $G/P$. In this section, we study a different boundary value map $I_o$, which carries $\Omega^d(\mathcal{X})'$ to (a quotient of) a space of functions on $G/P$. Our objective is to relate $I_o$ to $I$. The major complications come from the fact that the image of $I_o$ does not consist of locally analytic functions, a phenomenon essentially due to the fact that the kernel function on $G/P$ studied in [ST] is locally analytic on the big cell with continuous, not locally analytic, extension to $G/P$. We relate $I_o$ to $I$ using a "symmetrization map," due to Borel and Serre, which carries functions on $G/P$ into functions on $G$, together with a theory of "analytic vectors" in a continuous $G$-representation. One crucial consequence of our work in this section is the fact that the integral transform $I_o$ (and $I$) is injective.

Recall the definition, in [ST] Def. 27, of the integral kernel function $k(g, q)$ on $G/P \times \mathcal{X}$. This function is given by

$$k(g, .) = \begin{cases} (u_g)_* \frac{1}{\Xi_{\beta_0} \cdot \ldots \cdot \Xi_{\beta_{d-1}}} & \text{if } g = u_g w_{d+1} p \text{ is in the big cell,} \\ 0 & \text{otherwise.} \end{cases}$$

Here we rather want to consider the map

$$\begin{aligned} \mathbf{k} : G/P &\longrightarrow \Omega^d(\mathcal{X}) \\ g &\longmapsto k(g, .) d\Xi_{\beta_0} \wedge \ldots \wedge d\Xi_{\beta_{d-1}}. \end{aligned}$$

Since the numerator of the form $\xi$ is invariant under lower triangular unipotent matrices (compare the formula after Def. 28 in [ST]) we can rewrite our new map as

$$\mathbf{k}(g) = \begin{cases} (u_g)_* \xi & \text{if } g = u_g w_{d+1} p \text{ is in the big cell,} \\ 0 & \text{otherwise.} \end{cases}$$

**Proposition 1:**

*The map **k** is continuous and vanishes outside the big cell. Moreover whenever $\Omega^d(\mathcal{X})$ is equipped with the coarser topology coming from the spectral norm on $\mathcal{X}_n$ for some fixed but arbitrary $n \in \mathbb{N}$ then **k** is locally analytic on the big cell.*

Proof: The vanishing assertion holds by definition. The assertion about local analyticity of course is a consequence of Prop. 2.1'. But we will give another argument which actually produces explicitly the local series expansions. This will be needed in the subsequent considerations.

Let $U$ denote the unipotent radical of $P$. According to [ST] Lemma 12 the sets $\mathbf{B}(u,r) = uw_{d+1}t^r BP/P$, for a fixed $u \in U$, $t$ the diagonal matrix with entries $(\pi^d, \ldots, \pi, 1)$, and varying $r \in \mathbb{N}$, form a fundamental system of open neighbourhoods of the point $uw_{d+1}P/P$ in the big cell. One easily checks that

$$D(u,r) := \{v \in U : vw_{d+1}P/P \in \mathbf{B}(u,r)\}$$

is a polydisk in the affine space $U$. Hence the maps

$$\begin{array}{rcl} D(u,r) & \xrightarrow{\sim} & \mathbf{B}(u,r) \subseteq \text{ big cell} \\ v & \longmapsto & vw_{d+1}P/P \end{array}$$

constitute an atlas for the big cell as a locally analytic manifold. Fix $n \in \mathbb{N}$. We have to show that given a $u \in U$ we find an $r \in \mathbb{N}$ such that the map

$$\begin{array}{rcl} D(u,r) & \longrightarrow & \Omega^d(\mathcal{X}) \\ v & \longmapsto & \mathbf{k}(vw_{d+1}) \end{array}$$

is analytic with respect to the coarser topology on the right hand side corresponding to $n$. Recall that this amounts to the following ([B-VAR]). Let $v_{ji}$ for $0 \leq i < j \leq d$ denote the matrix entries of the matrix $v \in U$. Moreover we use the usual abbreviation

$$(v-u)^{\underline{m}} := \prod_{0 \leq i < j \leq d} (v_{ji} - u_{ji})^{m_{ji}}$$

for any multi-index $\underline{m} = (m_{10}, \ldots, m_{dd-1}) \in \mathbb{N}_0^{d(d+1)/2}$. We have to find an $r \in \mathbb{N}$ such that there is a power series expansion

$$k(vw_{d+1}, q) = \sum_{\underline{m}} (v-u)^{\underline{m}} \cdot F_{\underline{m}}(q)$$

with $F_{\underline{m}} \in \mathcal{O}(\mathcal{X})$ which is uniformly convergent on $D(u,r) \times \mathcal{X}_n$. From now on we fix $u \in U$. We choose $r \in \mathbb{N}$ such that

$$\omega(v_{ji} - u_{ji}) > 2n \quad \text{for all} \quad v \in D(u,r) \quad \text{and} \quad 0 \leq i < j \leq d.$$

We write
$$k(vw_{d+1}, q) = \prod_{i=0}^{d-1} \frac{1}{f_i(v,q)}$$

where
$$f_i(v,q) := \sum_{j=i}^{d-1} a_{ji}(v) \Xi_{\beta_j}(q) + a_{di}(v)$$

with
$$a_{ji}(v) := \begin{cases} v_{ji} & \text{for } j > i, \\ 1 & \text{for } j = i. \end{cases}$$

We also write
$$f_i(v,q) = f_i(u,q) + \sum_{j=i+1}^{d-1} b_{ji}(v) \Xi_{\beta_j}(q) + b_{di}(v)$$

with
$$b_{ji}(v) := a_{ji}(v) - u_{ji} = v_{ji} - u_{ji}.$$

Observe that
$$\omega(b_{ji}(.)) > 2n \quad \text{on} \quad D(u,r).$$

As was already discussed in the proof of [ST] Prop. 47 we have
$$\omega(f_i(u,q)) \leq n \quad \text{for} \quad q \in \mathcal{X}_n, \quad \text{and}$$
$$\omega\left( \sum_{j=i+1}^{d-1} b_{ji}(v) \Xi_{\beta_j}(q) + b_{di}(v) \right) > n \quad \text{for} \quad (v,q) \in D(u,r) \times \mathcal{X}_n.$$

Consequently
$$\frac{1}{f_i(v,q)} = \frac{1}{f_i(u,q)} \sum_{m \geq 0} (-1)^m \left( \frac{\sum_{j=i+1}^{d-1} b_{ji}(v) \Xi_{\beta_j}(q) + b_{di}(v)}{f_i(u,q)} \right)^m$$

is an expansion into a series uniformly convergent on $D(u,r) \times \mathcal{X}_n$. We rewrite this as
$$\frac{1}{f_i(v,q)}$$
$$= \sum_{m_{i+1 i},\ldots,m_{di} \geq 0} c_{\underline{m}(i)} \frac{\Xi_{m_{i+1 i}\beta_{i+1}+\ldots+m_{d-1 i}\beta_{d-1}}(q)}{f_i(u,q)^{1+m_{i+1 i}+\ldots+m_{di}}} \cdot \prod_{i<j\leq d} (v_{ji} - u_{ji})^{m_{ji}}$$

where $\underline{m}(i) := (m_{i+1\,i}, \ldots, m_{di})$ and the $c_{\underline{m}(i)}$ are certain nonzero integer coefficients. By multiplying together we obtain the expansion

$$(*) \quad k(vw_{d+1}, q) = \sum_{\underline{m}} \frac{c_{\underline{m}} \Xi_{\mu(\underline{m})}(q)}{f_0(u,q)^{s_0(\underline{m})} \cdot \ldots \cdot f_{d-1}(u,q)^{s_{d-1}(\underline{m})}} \cdot (v-u)^{\underline{m}}$$

which is uniformly convergent on $D(u,r) \times \mathcal{X}_n$; here we have set

$$\mu(\underline{m}) := m_{10}\beta_1 + (m_{20} + m_{21})\beta_2 + \ldots + (m_{d-1\,0} + \ldots + m_{d-1\,d-2})\beta_{d-1}$$

if $d > 1$, resp. $\mu(\underline{m}) := 0$ if $d = 1$, and

$$s_i(\underline{m}) := 1 + m_{i+1\,i} + \ldots + m_{di} \quad \text{for} \ \ 0 \leq i \leq d-1 \ ;$$

again the $c_{\underline{m}}$ are appropriate nonzero integer coefficients. This establishes the asserted local analyticity on the big cell. It follows immediately that $\mathbf{k}$ is continuous on the big cell (with respect to the original Fréchet topology on $\Omega^d(\mathcal{X})$). It therefore remains to prove, for all $n \in \mathbb{N}$, the continuity of $k$ viewed as a map from $G$ into $\mathcal{O}(\mathcal{X}_n)$ in all points outside the big cell. But this is the content of [ST] Lemma 45.

**Corollary 2:**

*The function $\lambda \circ \mathbf{k} : G/P \longrightarrow K$, for any continuous linear form $\lambda$ on $\Omega^d(\mathcal{X})$, is continuous, vanishes outside the big cell, and is locally analytic on the big cell.*

Proof: The continuity and the vanishing are immediately clear. The local analyticity follows by using [B-VAR] 4.2.3 and by observing that, according to Prop. 1.4, $\lambda$ comes from a continuous linear form on some $\Omega^d(\mathcal{X}_n)$.

**Proposition 3:**

*The image of $\mathbf{k}$ generates $\Omega^d(\mathcal{X})$ as a topological $K$-vector space.*

Proof: We first consider the map $k : G \longrightarrow \mathcal{O}(\mathcal{X})$. Let $\mathcal{K} \subseteq \mathcal{O}(\mathcal{X})$ be the vector subspace generated by the image of $k$ and let $\overline{\mathcal{K}}$ denote its closure. The formula $(*)$ in the proof of Proposition 1 for the matrix $u = 1$ says that, given the natural number $n \in \mathbb{N}$, we find an $r \in \mathbb{N}$ such that the expansion

$$k(vw_{d+1}, q) = \sum_{\underline{m}} \frac{c_{\underline{m}} \Xi_{\mu(\underline{m})}(q)}{\Xi_{\beta_0}(q)^{s_0(\underline{m})} \cdot \ldots \cdot \Xi_{\beta_{d-1}}(q)^{s_{d-1}(\underline{m})}} \cdot (v-1)^{\underline{m}}$$

holds uniformly for $(v, q) \in D(1, r) \times \mathcal{X}_n$. The coefficients of this expansion up to a constant are the value at $u = 1$ of iterated partial derivatives of the function $k(.w_{d+1}, .) : D(1, r) \longrightarrow \mathcal{K}$ (momentarily viewed in $\mathcal{O}(\mathcal{X}_n)$). Since increasing $n$ just means decreasing $r$ it follows that all the functions $\Xi_\mu$ with $\mu = \mu(\underline{m}) - s_0(\underline{m})\beta_0 - \ldots - s_{d-1}(\underline{m})\beta_{d-1}$ lie in $\overline{\mathcal{K}}$. This includes, for those $\underline{m}$ for which only the $m_{i+1 i}$ may be nonzero, all the functions

$$\frac{\Xi_{\beta_1}^{m_0} \cdot \ldots \cdot \Xi_{\beta_{d-1}}^{m_{d-2}}}{\Xi_{\beta_0}^{1+m_0} \cdot \ldots \cdot \Xi_{\beta_{d-1}}^{1+m_{d-1}}} = \Xi_{\alpha_0}^{m_0} \cdot \ldots \cdot \Xi_{\alpha_{d-1}}^{m_{d-1}} \cdot \frac{1}{\Xi_{\beta_0} \cdot \ldots \cdot \Xi_{\beta_{d-1}}}$$

with $m_0, \ldots, m_{d-1} \geq 0$ .

Passing now to $d$-forms we therefore know that the closed $K$-vector subspace $\Omega$ of $\Omega^d(\mathcal{X})$ generated by the image of $\mathbf{k}$ contains all forms $\Xi_\mu \xi$ where $\mu = m_0 \alpha_0 + \ldots + m_{d-1}\alpha_{d-1}$ with $m_0, \ldots, m_{d-1} \geq 0$. As a consequence of [ST] Cor. 40 the subspace $\Omega$ is $G$-invariant. By applying Weyl group elements $w$ and noting that $w_*\xi = \pm \xi$ we obtain $\Xi_\mu \xi$, for any $\mu \in X^*(\overline{T})$, in $\Omega$. Using the $G$-invariance of $\Omega$ again we then have the subset $\{u_*(\Xi_\mu \xi) : \mu \in X^*(\overline{T}), u \in U\} = \{(u_* \Xi_\mu) d\Xi_{\beta_0} \wedge \ldots \wedge d\Xi_{\beta_{d-1}} : \mu \in X^*(\overline{T}), u \in U\} \subseteq \Omega$. According to the partial fraction expansion argument in [GV] Thm. 21 the $u_* \Xi_\mu$ $K$-linearly span all rational functions of $\Xi_{\beta_0}, \ldots, \Xi_{\beta_{d-1}}$ whose denominator is a product of polynomials of degree 1. Moreover the proof of §1 Prop. 4 in [SS] shows that those la! tter functions are dense in $\mathcal{O}(\mathcal{X})$. It follows that $\Omega = \Omega^d(\mathcal{X})$. □

Put

$$C(G/P, K) := \text{ space of continuous } K\text{-valued functions on } G/P \ ;$$

it is a Banach space with respect to the supremum norm on which $G$ acts continuously by left translations. The subspace

$$C_{\text{inv}}(G/P, K) := \sum_s C(G/P_s, K) \subseteq C(G/P, K)$$

is closed; actually one has the topological direct sum decomposition

$$C(G/P, K) = C_{\text{inv}}(G/P, K) \oplus C_{\text{o}}(Pw_{d+1}P/P, K)$$

where the second summand on the right hand side is the space of $K$-valued continuous functions vanishing at infinity on the big cell ([BS] §3). We equip the quotient space $C(G/P, K)/C_{\text{inv}}(G/P, K)$ with the quotient topology. By Proposition 1 the map

$$\begin{aligned} I'_{\text{o}} : \Omega^d(\mathcal{X})' &\longrightarrow C(G/P, K) \\ \lambda &\longmapsto [g \longmapsto \lambda(\mathbf{k}(g))] \end{aligned}$$

is well defined; by [ST] Cor. 30 it is $P$-equivariant. Moreover it follows from [ST] Prop. 29.3 and the Bruhat decomposition that the induced map

$$I_{\text{o}} : \Omega^d(\mathcal{X})' \longrightarrow C(G/P, K)/C_{\text{inv}}(G/P, K)$$

is $G$-equivariant.

**Lemma 4:**

*The maps $I'_o$ and $I_o$ are continuous.*

Proof: We only need to discuss the map $I'_o$. Because of Prop. 1.4 we have to check that, for each $n \in \mathbb{N}$, the map

$$\begin{array}{rcl} \mathcal{O}(\mathcal{X}_n)' & \longrightarrow & C(G/P, K) \\ \lambda & \longmapsto & [g \longmapsto \lambda(k(g, .))] \end{array}$$

is continuous. The norm of $\lambda$ is equal to

$$c_1 := \inf\{\omega(\lambda(F)) : F \in \mathcal{O}(\mathcal{X}_n), \inf_{q \in \mathcal{X}_n} \omega(F(q)) \geq 0\} \ .$$

On the other hand the norm of the image of $\lambda$ under the above map is equal to

$$c_2 := \inf_{g \in G} \omega(\lambda(k(g, .))) \ = \ \inf_{u \in U} \omega(\lambda(k(uw_{d+1}, .)))$$

where $U$ denotes, as before, the unipotent radical of $P$. But we have

$$\inf_{\substack{u \in U \\ q \in \mathcal{X}_n}} \omega(k(uw_{d+1}, q)) \geq -dn$$

(compare the proof of [ST] Prop. 47). It follows that $c_2 \geq c_1 - dn$.

**Lemma 5:**

*The maps $I'_o$ and $I_o$ are injective.*

Proof: For $I'_o$ this is an immediate consequence of Prop. 3. According to Cor. 2 the image of $I'_o$ is contained in $C_o(Pw_{d+1}P/P, K)$ which is complementary to $C_{\mathrm{inv}}(G/P, K)$. Hence $I_o$ is injective, too. □

In order to see the relation between $I_o$ and the map $I$ in the previous section we first recall part of the content of [BS] §3:

**Fact 1:** The "symmetrization"

$$(\Sigma\phi)(g) := \sum_{w \in W} (-1)^{\ell(w)} \phi(gww_{d+1})$$

induces a $G$-equivariant injective map

$$C(G/P, K)/C_{\mathrm{inv}}(G/P, K) \xrightarrow{\Sigma} C(G, K) \ .$$

Here and in the following we let $C(Y, K)$, resp. $C_o(Y, K)$, denote, for any locally compact space $Y$, the $K$-vector space of $K$-valued continuous functions, resp. of $K$-valued continuous functions vanishing at infinity, on $Y$; the second space is a Banach space with respect to the supremum norm.

We also let
$$C(G, K) \xrightarrow{\text{res}} C(U, K)$$
$$\phi \longmapsto \phi|U$$

and

$$C_o(U, K) \longrightarrow C(G/P, K)/C_{\text{inv}}(G/P, K)$$
$$\phi \longmapsto \phi^{\#}(g) := \begin{cases} \phi(u) & \text{if } g = uw_{d+1}p \in Uw_{d+1}P, \\ 0 & \text{otherwise}. \end{cases}$$

**Fact 2:** $\#$ is an isomorphism whose inverse is $\text{res} \circ \Sigma$.
It follows in particular that $\#$ is an isometry.

Consider now the diagram

$$\begin{array}{ccc} C^{\text{an}}(G, K) & \xrightarrow{\subseteq} & C(G, K) \\ \uparrow I & & \uparrow \Sigma \\ \Omega^d(\mathcal{X})' & \xrightarrow{I_o} & C(G/P, K)/C_{\text{inv}}(G/P, K) \end{array}$$

in which all maps are $G$-equivariant and injective. We claim that the diagram is commutative; for that it suffices to prove the identity

$$(**) \qquad g_*\xi = \sum_{w \in W} (-1)^{\ell(w)} u_{gww_{d+1}*}\xi .$$

¿From Prop. 1 we know that each summand on the right hand side is a continuous function in $g \in G$ (where $u_{gw*}\xi := 0$ if $gw$ is not in the big cell). Hence it suffices to check the identity for $g$ in the dense open subset $\bigcap_{w \in W} Pw_{d+1}Pw$.
On the other hand it is an identity between logarithmic $d$-forms which can be checked after having applied the $G$-equivariant map "dis" into distributions on $G/P$; according to [ST] Remark on top of p. 423 the left hand side becomes

$$\sum_{w \in W} (-1)^{\ell(w)} \delta_{gw} = \sum_{w \in W} (-1)^{\ell(w)} \delta_{u_{gw}w_{d+1}}$$

whereas the right hand side becomes

$$\sum_{w \in W} (-1)^{\ell(w)} \cdot \sum_{v \in W} (-1)^{\ell(v)} \delta_{u_{gww_{d+1}}v} = \sum_{w \in W} (-1)^{\ell(w)} \cdot \sum_{v \in W} (-1)^{\ell(v)} \delta_{u_{gw}w_{d+1}v} .$$

The image of "dis" actually consists of linear forms on the Steinberg representation (see [ST]) and so any identity in that image can be checked by evaluation on locally constant and compactly supported functions on the big cell. But for those, all terms on the right hand side with $v \neq 1$ obviously vanish.

We view the above diagram as saying that any locally analytic function in the image of $I$ is the symmetrization of a continuous "boundary value function" on $G/P$. In order to make this more precise we first have to discuss the concept of an "analytic vector". Let $V$ be a $K$-Banach space on which $G$ acts continuously (by which we always mean that the map $G \times V \longrightarrow V$ describing the action is continuous). As in the case $V = K$ we have the Hausdorff locally convex vector space $C^{\text{an}}(G, V)$ of all $V$-valued locally $K$-analytic functions on $G$ (apart from replacing $K$ by $V$ everywhere the definition is literally the same). It is barrelled, so that the same argument as in the proof of Lemma 2.2 shows that the left translation action of $G$ on $C^{\text{an}}(G, V)$ is continuous.

**Definition:**

*A vector $v \in V$ is called analytic if the $V$-valued function $g \longmapsto gv$ on $G$ is locally analytic.*

We denote by $V_{\text{an}}$ the vector subspace of all analytic vectors in $V$. It is clearly $G$-invariant. Moreover the $G$-equivariant linear map

$$\begin{array}{rcl} V_{\text{an}} & \longrightarrow & C^{\text{an}}(G, V) \\ v & \longmapsto & [g \mapsto g^{-1}v] \end{array}$$

is injective. We always equip $V_{\text{an}}$ with the subspace topology with respect to this embedding. (Warning: That topology in general is finer than the topology which the Banach norm of $V$ would induce on $V_{\text{an}}$. Evaluating a function at $1 \in G$ defines a continuous map $C^{\text{an}}(G, V) \to V$.) Of course the $G$-action on $V_{\text{an}}$ is continuous. By functoriality any $G$-equivariant continuous linear map $L : V \to \tilde{V}$ between Banach spaces with continuous $G$-action induces a $G$-equivariant continuous linear map $L_{\text{an}} : V_{\text{an}} \to \tilde{V}_{\text{an}}$. A useful technical observation is that the locally convex vector space $V_{\text{an}}$ does not change if we pass to an open subgroup $H \subseteq G$. First of all it follows from the continuity of the $G$-action on $V$ that the function $g \mapsto g^{-1}v$ is locally analytic on $G$ if and only if its restriction to $H$ is locally analytic. Fixing a set of representatives $R$ for the cosets in $H \setminus G!$ we have the isomorphism of locally convex vector spaces

$$C^{\text{an}}(G, V) = \prod_{g \in R} C^{\text{an}}(Hg, V)$$

([Fea] 2.2.4). Hence the embedding $V_{\text{an}} \hookrightarrow C^{\text{an}}(G, V)$ coincides with the composite of the embedding $V_{\text{an}} \hookrightarrow C^{\text{an}}(H, V)$ and the "diagonal embedding"

$$C^{\text{an}}(H, V) \longrightarrow \prod_{g \in R} C^{\text{an}}(Hg, V)$$
$$f \longmapsto (g^{-1}(f(.g^{-1})))_{g \in R} \ .$$

**Remark 6:**

$V_{\text{an}}$ is closed in $C^{\text{an}}(G, V)$.

Proof: Let $(v_i)_{i \in I}$ be a Cauchy net in $V_{\text{an}}$ which in $C^{\text{an}}(G, V)$ converges to the function $f$. By evaluating at $h \in G$ we see that the net $(h^{-1}v_i)_{i \in I}$ converges to $f(h)$ in $V$. Put $v := f(1)$. Since $h$ is a continuous endomorphism of $V$ it follows on the other hand that $(h^{-1}v_i)_{i \in I}$ converges to $h^{-1}v$. Hence $f(h) = h^{-1}v$ which means that $f$ comes from $v \in V_{\text{an}}$.

**Lemma 7:**

*If each vector in $V$ is analytic then $V_{\text{an}} = V$ as topological vector spaces.*

Proof: We have to show that the map

$$V \longrightarrow C^{\text{an}}(G, V)$$
$$v \longmapsto [g \mapsto g^{-1}v]$$

is continuous. According to our earlier discussion we are allowed to replace $G$ by whatever open subgroup $H$ is convenient. By [Fea] 3.1.9 our assumption implies that the $G$-action on $V$ defines a homomorphism of Lie groups $\rho : G \longrightarrow \mathbf{GL}(V)$ (with the operator norm topology on the right hand side). On a sufficiently small compact open subgroup $H \subseteq G$ this homomorphism is given by a power series

$$\rho(g) = \sum_{\underline{n}} A_{\underline{n}} \cdot \underline{x}(g)^{\underline{n}} \text{ for } g \in H$$

which is convergent in the operator norm topology on $\text{End}_K(V)$; here $\underline{x}$ is a vector of coordinate functions from $H$ onto some polydisk of radius 1, the $\underline{n}$ are corresponding multi-indices, and the $A_{\underline{n}}$ lie in $\text{End}_K(V)$. In particular the operator norm of the $A_{\underline{n}}$ is bounded above by some constant $c > 0$. If we insert a fixed vector $v \in V$ into this power series then we obtain the expansion

$$gv = \sum_{\underline{n}} A_{\underline{n}}(v) \cdot \underline{x}(g)^{\underline{n}}$$

as a function of $g \in H$ and the spectral norm of the right hand side is bounded above by $c \cdot \|v\|$.

**Proposition 8:**

*The map $I_o$ induces a $G$-equivariant injective continuous linear map*

$$\Omega^d(\mathcal{X})' \longrightarrow [C(G/P,K)/C_{\text{inv}}(G/P,K)]_{\text{an}} \ .$$

Proof: For the purposes of this proof we use the abbreviation $V := C(G/P,K)/C_{\text{inv}}(G/P,K)$. We have to show that the image of $I_o$ is contained in $V_{\text{an}}$ and that the induced map into $V_{\text{an}}$ is continuous. As before it suffices to discuss the corresponding map $\Omega^d(\mathcal{X}_n)' \longrightarrow V_{\text{an}}$ for a fixed but arbitrary $n \in \mathbb{N}$. Both spaces, $V$ as well as $\Omega^d(\mathcal{X}_n)'$, are Banach spaces with an action of the group $GL_{d+1}(o)$; the map between them induced by $I_o$ is equivariant and continuous by Lemma 4. Since $GL_{d+1}(o)$ is open in $G$ it can be used, by the above observation, instead of $G$ to compute the locally convex vector space $V_{\text{an}}$. If we show that $GL_{d+1}(o)$ acts continuously on $\Omega^d(\mathcal{X}_n)'$ then $I_o$ certainly induces a continuous map $[\Omega^d(\mathcal{X}_n)']_{\text{an}} \longrightarrow V_{\text{an}}$. What we therefore have to show in addition is that the identity

$$[\Omega^d(\mathcal{X}_n)']_{\text{an}} \ = \ \Omega^d(\mathcal{X}_n)'$$

holds as topological vector spaces.
We know already from the proof of Prop. 2.1' that every vector in $\Omega^d(\mathcal{X}_n)$ is analytic. By [Fea] 3.1.9 this means that the $GL_{d+1}(o)$-action on $\Omega^d(\mathcal{X}_n)$ is given by a homomorphism of Lie groups

$$GL_{d+1}(o) \ \longrightarrow \ \mathbf{GL}(\Omega^d(\mathcal{X}_n))$$

(recall that the right hand side carries the operator norm topology). Since passing to the adjoint linear map is a continuous linear map between Banach spaces it follows that also the $GL_{d+1}(o)$-action on $\Omega^d(\mathcal{X}_n)'$ is given by a corresponding homomorphism of Lie groups. This means in particular that the latter action is continuous and that every vector in $\Omega^d(\mathcal{X}_n)'$ is analytic. We therefore may apply the previous lemma.

**Corollary 9:**

*The $G$-action $G \times \Omega^d(\mathcal{X})' \longrightarrow \Omega^d(\mathcal{X})'$ is continuous and, for any $\lambda \in \Omega^d(\mathcal{X})'$, the map $g \longmapsto g\lambda$ on $G$ is locally analytic.*

Proof: Since $\Omega^d(\mathcal{X})'$ is barrelled as a locally convex inductive limit of Banach spaces the first assertion follows from the second by the same argument which we have used already twice. In the proof of the previous proposition we have seen that each function $g \longmapsto g\lambda$ is locally analytic on $GL_{d+1}(o)$. But this is sufficient for the full assertion. □

Since $G/P$ is compact we may view the symmetrization map as a map

$$C(G/P,K)/C_{\text{inv}}(G/P,K) \ \xrightarrow{\Sigma} \ BC(G,K)$$

into the Banach space $BC(G, K)$ of bounded continuous functions on $G$. It is then an isometry as can be seen as follows. By its very definition $\Sigma$ is norm decreasing. On the other hand $\phi$ can be reconstructed from $\Sigma\phi$ by restriction to $U$ followed by $\#$ which again is norm decreasing. Hence $\Sigma$ must be norm preserving.

We obtain the induced continuous injective map

$$[C(G/P, K)/C_{\text{inv}}(G/P, K)]_{\text{an}} \stackrel{\Sigma_{\text{an}}}{\hookrightarrow} BC(G, K)_{\text{an}}.$$

Since any $f \in BC(G, K)_{\text{an}}$ is obtained from the locally analytic map $g \longmapsto g^{-1}f$ by composition with the evaluation map at $1 \in G$ and hence is locally analytic we see that $BC(G, K)_{\text{an}}$ in fact is contained in $C^{\text{an}}(G, K)$. We therefore can rewrite the commutative diagram which relates $I_{\text{o}}$ and $I$ in the form

$$\Omega^d(\mathcal{X})' \xrightarrow{I_{\text{o}}} [C(G/P, K)/C_{\text{inv}}(G/P, K)]_{\text{an}}$$
$$I \searrow \quad \Big\downarrow \Sigma_{\text{an}}$$
$$C^{\text{an}}(G, K) \ .$$

So far we have explained how to understand on $G/P$ the fact that the functions in the image of $I$ are locally analytic. But the latter also satisfy the differential equations from the ideal $\mathfrak{a}$ in $\mathcal{D}(G)$. How can those be viewed on $G/P$?

For any Banach space $V$ the right translation action by $G$ on $C^{\text{an}}(G, V)$ induces a corresponding action of the algebra $\mathcal{D}(G)$ by continuous endomorphisms. Any $\mathfrak{x} \in \mathfrak{g}$ acts via the usual formula

$$(\mathfrak{x}f)(g) \ = \ \frac{d}{dt}f(g\exp(t\mathfrak{x}))\big|_{t=0} \ .$$

(Compare [Fea] 3.3.4.) This clearly is functorial in $V$. If we now look at the case $BC(G, K)$ we have two embeddings

$$BC(G, K)_{\text{an}}$$
$$\swarrow \qquad \searrow$$
$$C^{\text{an}}(G, BC(G, K)) \xrightarrow{\varepsilon} C^{\text{an}}(G, K)$$

which are connected through the map $\varepsilon$ which comes by functoriality from the map $\text{ev}_1 : BC(G, K) \longrightarrow K$ evaluating a function at $1 \in G$. This latter map is $\mathcal{D}(G)$-equivariant. Hence, for any left ideal $\mathfrak{d} \subseteq \mathcal{D}(G)$, we obtain the identity

$$BC(G, K)_{\text{an}} \ \cap \ C^{\text{an}}(G, K)^{\mathfrak{d}=0} \ =$$
$$BC(G, K)_{\text{an}} \ \cap \ \{f \in C^{\text{an}}(G, BC(G, K)) : \mathfrak{d}f \subseteq \ker(\varepsilon)\} \ .$$

Using the abbreviation $V := C(G/P, K)/C_{\mathrm{inv}}(G/P, K)$ we know from Prop. 8 that $\mathrm{im}(I_{\mathrm{o}}) \subseteq V_{\mathrm{an}}$; on the other hand $\mathrm{im}(I) \subseteq C^{\mathrm{an}}(G, K)^{\mathfrak{a}=0}$. Those images correspond to each other under the map $\Sigma_{\mathrm{an}}$. It follows that $\mathrm{im}(I_{\mathrm{o}})$ is contained in the subspace

$$V_{\mathrm{an}}^{\sigma \mathfrak{a}=0} := V_{\mathrm{an}} \cap \{f \in C^{\mathrm{an}}(G, V) : \mathfrak{a} f \subseteq \ker(\sigma)\}$$

where $\sigma : C^{\mathrm{an}}(G, V) \longrightarrow C^{\mathrm{an}}(G, K)$ is the map induced by $\mathrm{ev}_1 \circ \Sigma : V \longrightarrow K$ which sends $\phi \in C(G/P, K)$ to $\sum_{w \in W}(-1)^{\ell(w)}\phi(ww_{d+1})$. We arrive at the following conclusion.

**Theorem 10:**

*We have the commutative diagram of injective continuous linear maps*

$$\Omega^d(\mathcal{X})' \xrightarrow{I_{\mathrm{o}}} [C(G/P, K)/C_{\mathrm{inv}}(G/P, K)]_{\mathrm{an}}^{\sigma \mathfrak{a}=0}$$

$$I \searrow \quad \downarrow \Sigma_{\mathrm{an}}$$

$$C^{\mathrm{an}}(G, K)^{\mathfrak{a}=0} \ .$$

We think of $I_{\mathrm{o}}$ in this form as being "the" *boundary value map*. We point out that the ideal $\mathfrak{a}$ contains the following Dirac distributions. For any $g \in G$ put $W(g) := \{w \in W : gww_{d+1} \in Pw_{d+1}P\}$ and consider the element

$$\delta(g) := \delta_g - \sum_{w \in W(g)}(-1)^{\ell(w)}\delta_{u_{gww_{d+1}}}$$

in the algebra $\mathcal{D}(G)$. By interpreting $\delta_h$ as the right translation action by $h$ those elements act on $C(G, K)$ and $BC(G, K)$. We claim that any function $\phi$ in the image of $\Sigma$ satisfies

$$\phi(g) = \sum_{w \in W(g)}(-1)^{\ell(w)}\phi(u_{gww_{d+1}}) \text{ for any } g \in G \ .$$

We may write $\phi = \Sigma \psi^{\#}$ and then compute

$$\sum_{w \in W(g)}(-1)^{\ell(w)}(\Sigma \psi^{\#})(u_{gww_{d+1}}) =$$

$$\sum_{w \in W(g)}(-1)^{\ell(w)} \sum_{v \in W}(-1)^{\ell(v)}\psi^{\#}(u_{gww_{d+1}}vw_{d+1}) =$$

$$\sum_{w \in W(g)}(-1)^{\ell(w)}\psi(u_{gww_{d+1}}) = \sum_{w \in W}(-1)^{\ell(w)}\psi^{\#}(gww_{d+1}) = (\Sigma \psi^{\#})(g) \ .$$

Since $\Sigma$ is $G$-equivariant the same identities hold for the functions $\phi(h.)$ for any $h \in G$. In other words any function $\phi$ in the image of $\Sigma$ actually satisfies

$$\delta(g)\phi \;=\; 0 \;\text{ for any } g \in G \;.$$

The relation $(**)$ implies that $\mathfrak{a}$ contains the left ideal generated by the $\delta(g)$ for $g \in G$.

## 4. The ideal $\mathfrak{b}$

As we have learned, the integral transform $I$ carries continuous linear forms on $\Omega^d(\mathcal{X})$ to locally analytic functions on $G$. Functions in the image of this map are annihilated by an ideal $\mathfrak{a}$ in the algebra of punctual distributions $\mathcal{D}(G)$. This annihilation condition means that functions in the image of $I$ satisfy a mixture of discrete relations and differential equations.

In this section, we focus our attention on the differential equations satisfied by functions in the image of $I$. By this, we mean that we will study in detail the structure of the ideal $\mathfrak{b} := \mathfrak{a} \cap U(\mathfrak{g})$. By definition, $\mathfrak{b}$ is the annihilator ideal in $U(\mathfrak{g})$ of the special differential form $\xi$. We will describe a set of generators for $\mathfrak{b}$ and use this to prove the fundamental result that the weight spaces in $U(\mathfrak{g})/\mathfrak{b}$ (under the adjoint action of the torus $\overline{T}$) are one-dimensional. We will then analyze the left $U(\mathfrak{g})$-module $U(\mathfrak{g})/\mathfrak{b}$, identifying a filtration of this module by submodules and exhibiting the subquotients of this filtration as certain explicit irreducible highest weight $U(\mathfrak{g})$-modules. At the end of the section we prove some additional technical structural results which we will need later.

The results in this section are fundamental preparation for the rest of the paper.

We begin by recalling the decomposition

$$U(\mathfrak{g}) = \oplus_{\mu \in X^*(\overline{T})} U(\mathfrak{g})_\mu$$

of $U(\mathfrak{g})$ into the weight spaces $U(\mathfrak{g})_\mu$ with respect to the adjoint action of the torus $\overline{T}$. For a root $\alpha = \varepsilon_i - \varepsilon_j$ the weight space $\mathfrak{g}_\alpha$ is the 1-dimensional space generated by the element $L_\alpha \in \mathfrak{g}$ which corresponds to the matrix with a 1 in position $(i,j)$ and zeros elsewhere; sometimes we also write $L_{ij} := L_\alpha$. Clearly a monomial $L_{\alpha_1}^{m_1} \cdot \ldots \cdot L_{\alpha_r}^{m_r} \in U(\mathfrak{g})$ has weight $m_1\alpha_1 + \ldots + m_r\alpha_r$. The Poincaré-Birkhoff-Witt theorem says that once we have fixed a total ordering of the roots $\alpha$ any element in $U(\mathfrak{g})/U(\mathfrak{g})\mathfrak{g}_o$ can be written in a unique way as a polynomial in the $L_\alpha$. We will also need the filtration $U_n(\mathfrak{g})$ of $U(\mathfrak{g})$ by degree; we write $\deg(\mathfrak{z}) := n$ if $\mathfrak{z} \in U_n(\mathfrak{g}) \setminus U_{n-1}(\mathfrak{g})$.

The form $\xi$ is invariant under $\overline{T}$. This implies that the ideal $\mathfrak{b}$ is homogeneous and contains $U(\mathfrak{g})\mathfrak{g}_o$. An elementary calculation shows that $L_\alpha$ acts on $\Omega^d(\mathcal{X})$ by

$$L_\alpha(F\xi) = (\Xi_i \frac{\partial F}{\partial \Xi_j})\xi - \Xi_\alpha F\xi \ .$$

In particular we obtain $L_\alpha \xi = -\Xi_\alpha \xi$. By iteration that formula implies that the ideal $\mathfrak{b}$ contains the following relations:

(cancellation) $\qquad\qquad L_{ij}L_{jl} \qquad\qquad$ for any indices $i \neq j \neq l$,
(sorting) $\qquad L_{ij}L_{k\ell} - L_{i\ell}L_{kj} \qquad$ for any distinct indices $(i,j,k,l)$.

Our goal is to show that the weight spaces of $U(\mathfrak{g})/\mathfrak{b}$ are 1-dimensional. For that we need to introduce one more notation. For a weight $\mu$ we put

$$d(\mu) := \sum_{m_i > 0} m_i$$

where the $m_i$ are the coefficients of $\mu$ in the linear combination

$$\mu = \sum_{i=0}^{d} m_i \varepsilon_i \ .$$

**Lemma 1:**

Let $\mathfrak{z} \in U(\mathfrak{g})$ be a monomial in the $L_\alpha$ of weight $\mu$; we then have:

i. $\deg(\mathfrak{z}) \geq d(\mu)$;

ii. write $\mathfrak{z} = \prod_{i,j} L_{ij}^{n_{ij}}$ and put $A(\mathfrak{z}) := \{i : n_{ij} > 0 \text{ for some } j\}$ and $B(\mathfrak{z}) := \{j : n_{ij} > 0 \text{ for some } i\}$; then $\deg(\mathfrak{z}) = d(\mu)$ if and only if $A(\mathfrak{z})$ and $B(\mathfrak{z})$ are disjoint.

Proof: (Recall that we have fixed a total ordering of the roots $\alpha$.) Since $\mathfrak{z}$ has weight $\mu$ we must have

$$\mu = \sum_{i,j} n_{ij}(\varepsilon_i - \varepsilon_j) \ .$$

If on the other hand we write $\mu = \sum_k m_k \varepsilon_k$ we see that

(*) $$m_k = \sum_j n_{kj} - \sum_i n_{ik} \ .$$

As a result of this expression it follows that $m_k \leq \sum_j n_{kj}$, so that

$$d(\mu) = \sum_{m_k > 0} m_k \leq \sum_k \sum_j n_{kj} = \deg(\mathfrak{z}) \ .$$

We suppose now that $A(\mathfrak{z})$ and $B(\mathfrak{z})$ are disjoint. Then $(*)$ implies that $m_k$ is positive if and only if $k \in A(\mathfrak{z})$ and for positive $m_k$ we must have

$$m_k = \sum_j n_{kj} \ .$$

Therefore

$$d(\mu) = \sum_{k \in A(\mathfrak{z})} m_k = \sum_{k \in A(\mathfrak{z})} \sum_j n_{kj} = \sum_{k,j} n_{kj} = \deg(\mathfrak{z}) \ .$$

Conversely, we suppose that $k \in A(\mathfrak{z}) \cap B(\mathfrak{z})$. Then $m_k < \sum_j n_{kj}$. Therefore, if $m_k \geq 0$, we obtain

$$d(\mu) = \sum_{m_k > 0} m_k < \sum_k \sum_j n_{kj} = \deg(\mathfrak{z}) \ .$$

If $m_k \leq 0$ a similar argument, using the fact that $d(\mu)$ may be computed from the $m_k$ with $m_k < 0$, gives the desired result.

**Lemma 2:**

*Let $\mathfrak{z} \in U(\mathfrak{g})$ be a nonzero polynomial in the $L_\alpha$ of weight $\mu$; then the coset $\mathfrak{z} + \mathfrak{b}$ contains a representative of weight $\mu$ which is a linear combination of monomials in the $L_\alpha$ of degree $d(\mu)$.*

Proof: Among all elements of $\mathfrak{z} + \mathfrak{b}_\mu$ which are polynomials in the $L_\alpha$ let $\mathfrak{x}$ be one of minimal degree. By the preceeding lemma the degree of $\mathfrak{x}$ is greater than or equal to $d(\mu)$. Let $\mathfrak{y}$ be a monomial in the $L_\alpha$ of degree $\deg(\mathfrak{x})$ which occurs with a nonzero coefficient in $\mathfrak{x}$. Assume that the sets $A(\mathfrak{y})$ and $B(\mathfrak{y})$ as defined in the preceeding lemma are not disjoint. Then there exist three indices $i \neq j \neq l$ such that $L_{ij}$ and $L_{jl}$ each occur to nonzero powers in the monomial $\mathfrak{y}$. By the commutation rules in $U(\mathfrak{g})$ we have

$$\mathfrak{y} \in U(\mathfrak{g}) L_{ij} L_{jl} + U_{n-1}(\mathfrak{g}) \ \text{ with } \ n := \deg(\mathfrak{x}) \ .$$

Hence the cancellation relations imply that $\mathfrak{y} \in \mathfrak{b} + U_{n-1}(\mathfrak{g})$. This means that modulo $\mathfrak{b}$ we may remove an appropriate scalar multiple of $\mathfrak{y}$ from $\mathfrak{x}$ and pick up only a polynomial of lower degree. But by our minimality assumption on $\deg(\mathfrak{x})$ there has to be at least one such $\mathfrak{y}$ such that $A(\mathfrak{y})$ and $B(\mathfrak{y})$ are disjoint. The previous lemma then implies that $d(\mu) = \deg(\mathfrak{y}) = \deg(\mathfrak{x})$. If we express $\mathfrak{x}$ as a linear combination of monomials in the $L_\alpha$ then each such monomial has weight $\mu$ and hence, by the previous lemma again, degree $\geq d(\mu)$. $\square$

A monomial

$$L_{i_0 j_0} L_{i_1 j_1} \ldots L_{i_m j_m}$$

will be called *sorted* if $i_0 \leq \ldots \leq i_m$ and $j_0 \leq \ldots \leq j_m$ and if those two sequences do not overlap (i.e. no $i_k$ is a $j_l$). For example, the monomial $L_{10}L_{10}L_{32}L_{32}$ is sorted with sequences 1,1,3,3, and 0,0,2,2 whereas the monomial $L_{32}L_{31}$ with sequences 3,3 and 2,1 is not sorted.

**Lemma 3:**

*Among all the monomials in the $L_\alpha$ of weight $\mu$ there is exactly one, denoted by $L_{(\mu)}$, which is sorted.*

Proof: The non-overlapping condition means that any sorted monomial must have degree $d(\mu)$. Write $\mu = \sum m_k \varepsilon_k$. Those $k$ for which $m_k$ is positive must occur as the first index in some $L_{ij}$, and therefore (by the non-overlapping condition) can only occur as first indices. Similarly, those $k$ for which $m_k$ is negative can only occur as second indices. This determines the lists of first and second indices – for example, the list of first indices consists of precisely those $k$ for which $m_k > 0$, each repeated $m_k$ times, listed in ascending order. Once these two lists are determined the corresponding monomial is determined.

**Proposition 4:**

*Let $\mathfrak{z} \in U(\mathfrak{g})$ be a polynomial in the $L_\alpha$ of weight $\mu$; we then have $\mathfrak{z}+\mathfrak{b} = aL_{(\mu)}+\mathfrak{b}$ for some $a \in K$.*

Proof: By Lemma 2 we may assume that $\mathfrak{z}$ is a monomial of degree $d(\mu)$. The sets $A(\mathfrak{z})$ and $B(\mathfrak{z})$, as defined in Lemma 1, then are disjoint, and therefore the individual $L_{ij}$ which occur in $\mathfrak{z}$ commute with one another. Consequently we may rearrange these $L_{ij}$ freely. Using this fact it is easy to see that we may use the sorting relations to transform $\mathfrak{z}$ into $L_{(\mu)}$.

**Corollary 5:**

*The weight space $(U(\mathfrak{g})/\mathfrak{b})_\mu$ in the left $U(\mathfrak{g})$-module $U(\mathfrak{g})/\mathfrak{b}$, for any $\mu \in X^*(\overline{T})$, has dimension one.*

Proof: The preceeding proposition says that the weight space in question is generated by the coset $L_{(\mu)} + \mathfrak{b}$. On the other hand an explicit computation shows that $L_{(\mu)}\xi = e2^c \Xi_\mu \xi$ with some integer $c \geq 0$ and some sign $e = \pm 1$ (both depending on $\mu$); hence $L_{(\mu)} \notin \mathfrak{b}$. □

Later on it will be more convenient to use a renormalized $L_{(\mu)}$. We let $L_\mu$ denote the unique scalar multiple of $L_{(\mu)}$ which has the property that $L_\mu \xi = -\Xi_\mu \xi$.

Although we now have a completely explicit description of $U(\mathfrak{g})/\mathfrak{b}$ its structure as a $\mathfrak{g}$-module is not yet clear. For a root $\alpha = \varepsilon_i - \varepsilon_j$ and a weight $\mu = \sum_k m_k \varepsilon_k \in X^*(\overline{T})$ our earlier formula implies

$$L_\alpha(\Xi_\mu \xi) = (m_j - 1)\Xi_{\mu+\alpha}\xi$$

and hence

(+) $$L_\alpha L_\mu \equiv (m_j - 1)L_{\mu+\alpha} \mod \mathfrak{b} .$$

If we put $J(\mu) := \{0 \leq k \leq d : m_k > 0\}$ then $J(\mu) \subseteq J(\mu + \alpha)$ provided $m_j \neq 1$. It follows that
$$\mathfrak{b}_J := \mathfrak{b} + \sum_{J \subseteq J(\mu)} KL_\mu$$
is, for any subset $J \subseteq \{0, \ldots, d\}$, a left ideal in $U(\mathfrak{g})$. We have:

- $\mathfrak{b}_{\{0,\ldots,d\}} = \mathfrak{b}$ and $\mathfrak{b}_\emptyset = U(\mathfrak{g})$;
- $\mathfrak{b}_J \subseteq \mathfrak{b}_{J'}$ if and only if $J' \subseteq J$.

For $J \neq \{0, \ldots, d\}$ we set
$$\mathfrak{b}_J^> := \sum_{J \subsetneq J'} \mathfrak{b}_{J'} .$$

Moreover we introduce the descending filtration by left ideals
$$U(\mathfrak{g}) = \mathfrak{b}_0 \supseteq \mathfrak{b}_1 \supseteq \ldots \supseteq \mathfrak{b}_{d+1} = \mathfrak{b}$$
defined by
$$\mathfrak{b}_j := \sum_{\#J \geq j} \mathfrak{b}_J .$$

The subquotients of that filtration decompose as $\mathfrak{g}$-modules into
$$\mathfrak{b}_j/\mathfrak{b}_{j+1} = \bigoplus_{\#J=j} (\mathfrak{b}_J + \mathfrak{b}_{j+1})/\mathfrak{b}_{j+1} = \bigoplus_{\#J=j} \mathfrak{b}_J/\mathfrak{b}_J^> .$$

Our aim in the following therefore is to understand the $\mathfrak{g}$-modules $\mathfrak{b}_J/\mathfrak{b}_J^>$. A trivial case is
$$\mathfrak{b}_0/\mathfrak{b}_1 = \mathfrak{b}_\emptyset/\mathfrak{b}_\emptyset^> = K .$$

We therefore assume, for the rest of this section, that $J$ is a nonempty proper subset of $\{0, \ldots, d\}$. First of all we need the maximal parabolic subalgebra of $\mathfrak{g}$ given by
$$\mathfrak{p}_J := \text{all matrices in } \mathfrak{g} \text{ with a zero entry in position } (i,j) \text{ for } i \in J \text{ and } j \notin J .$$

It follows from the above formula (+) that the subalgebra $\mathfrak{p}_J$ leaves invariant the finite dimensional subspace
$$M_J := \sum_{\mu \in B(J)} KL_\mu$$

of $\mathfrak{b}_J/\mathfrak{b}_J^>$ where

$$B(J) := \text{set of all weights } \mu = \sum_k m_k \varepsilon_k \text{ such that } J(\mu) = J \text{ and } m_k = 1 \text{ for } k \in J.$$

Using again the formula (+) the subsequent facts are straightforward. The unipotent radical

$$\mathfrak{n}_J := \text{all matrices with zero entries in position } (i,j) \text{ with } i \in J \text{ or } j \notin J$$

of $\mathfrak{p}_J$ acts trivially on $M_J$. We have the Levi decomposition $\mathfrak{p}_J = \mathfrak{l}_J + \mathfrak{n}_J$ with $\mathfrak{l}_J = \mathfrak{l}'(J) + \mathfrak{l}(J)$ where

$$\mathfrak{l}'(J) := \text{all matrices with zero entries in position } (i,j) \text{ with } i \text{ and } j \text{ not both in } J$$

and

$$\mathfrak{l}(J) := \text{all matrices with zero entries in position } (i,j) \text{ with } i \text{ or } j \in J.$$

The structure of $M_J$ as a module for the quotient $\mathfrak{p}_J/\mathfrak{n}_J = \mathfrak{l}_J$ is as follows:

- The first factor $\mathfrak{l}'(J) \cong \mathfrak{gl}_{\#J}$ acts on $M_J$ through the trace character;
- as a module for the second factor $\mathfrak{l}(J) \cong \mathfrak{gl}_{d+1-\#J}$ our $M_J$ is isomorphic to the $\#J$-th symmetric power of the contragredient of the standard representation of $\mathfrak{gl}_{d+1-\#J}$ on the $(d+1-\#J)$-dimensional $K$-vector space.

In particular, $M_J$ is an irreducible $\mathfrak{p}_J/\mathfrak{n}_J$-module. The map

$$U(\mathfrak{g}) \underset{U(\mathfrak{p}_J)}{\otimes} M_J \longrightarrow \mathfrak{b}_J/\mathfrak{b}_J^>$$

$$(\mathfrak{z}, m) \longmapsto \mathfrak{z}m$$

is surjective. In fact, $\mathfrak{b}_J/\mathfrak{b}_J^>$ is an irreducible highest weight $U(\mathfrak{g})$-module: If we put $\nu := (\sum_{k \in J} \varepsilon_k) - \#J \cdot \varepsilon_\ell$ for some fixed $\ell \notin J$, then one deduces from (+) that $U(\mathfrak{g}) \cdot L_\mu + \mathfrak{b}$, for any $\mu$ with $J(\mu) = J$, contains $L_\nu + \mathfrak{b}$. For the subset $J = \{0, \ldots, j-1\}$ the parabolic subalgebra $\mathfrak{p}_J$ is in standard form with respect to our choice of positive roots and the highest weight of $\mathfrak{b}_J/\mathfrak{b}_J^>$ is $\varepsilon_0 + \ldots + \varepsilon_{j-1} - j \cdot \varepsilon_j$.

We finish this section by establishing several facts to be used later on about the relation between the left ideals $\mathfrak{b}_J^>$ and the subalgebras $U(\mathfrak{n}_J^+)$ for

$$\mathfrak{n}_J^+ := \text{transpose of } \mathfrak{n}_J.$$

First of all, note that $\mathfrak{n}_J^+$ and hence each $U(\mathfrak{n}_J^+)$ is commutative and $\mathrm{ad}(\mathfrak{l}_J)$-invariant.

**Proposition 6:**

i. $U(\mathfrak{n}_J^+) \cap \mathfrak{b}$ is the ideal in $U(\mathfrak{n}_J^+)$ generated by the sorting relations $L_{ij}L_{k\ell} - L_{i\ell}L_{kj}$ for $i, k \in J$ and $j, l \notin J$,
ii. the cosets of the sorted monomials $L_\mu$ for $J(\mu) \subseteq J$ and $J(-\mu) \cap J = \emptyset$ form a basis of $U(\mathfrak{n}_J^+)/U(\mathfrak{n}_J^+) \cap \mathfrak{b}$ as a $K$-vector space;
iii. $U(\mathfrak{n}_J^+) \cap \mathfrak{b}_J^> = U(\mathfrak{n}_J^+) \cap \mathfrak{b}$;
iv. $U(\mathfrak{n}_J^+) \cap \mathfrak{b}) \cdot \mathfrak{l}(J) \subseteq \mathfrak{b}$;
v. $U(\mathfrak{n}_J^+) \cap \mathfrak{b}$ is $\mathrm{ad}(\mathfrak{l}(J))$-invariant.

Proof: i. Let $\mathfrak{s} \subseteq U(\mathfrak{n}_J^+)$ denote the ideal generated by those sorting relations. Using the commutativity of $U(\mathfrak{n}_J^+)$ it is easy to see that any monomial in the $L_{ij}$ in $U(\mathfrak{n}_J^+)$ can be transformed into a sorted monomial by relations in $\mathfrak{s}$. In particular any coset in $U(\mathfrak{n}_J^+) \cap \mathfrak{b}/\mathfrak{s}$ has a representative which is a linear combination of sorted monomials. But we know that the sorted monomials are linearly independent modulo $\mathfrak{b}$. We therefore must have $U(\mathfrak{n}_J^+) \cap \mathfrak{b} = \mathfrak{s}$.
ii. The argument just given also shows that the cosets of all sorted monomials contained in $U(\mathfrak{n}_J^+)$ form a basis of the quotient in question. But they are exactly those which we have listed in the assertion.
iii. The sorted monomials listed in the assertion ii. are linearly independent modulo $\mathfrak{b}_J^>$.
iv. We have to check that $\mathfrak{z} L_{rs} \in \mathfrak{b}$ for any of the sorting relations $\mathfrak{z} = L_{ij}L_{k\ell} - L_{i\ell}L_{kj}$ from i. and for any $L_{rs}$ such that $r, s \notin J$. If $r \neq j, \ell$ then $\mathfrak{z} L_{rs} = L_{rs}\mathfrak{z} \in \mathfrak{b}$. If $r = \ell$ then we may use the cancellation relations to obtain $\mathfrak{z} L_{rs} = L_{ij}L_{k\ell}L_{\ell s} - L_{kj}L_{i\ell}L_{\ell s} \in \mathfrak{b}$. Similarly if $r = j$ we have $\mathfrak{z} L_{rs} = L_{k\ell}L_{ij}L_{js} - L_{i\ell}L_{kj}L_{js} \in \mathfrak{b}$.
v. Because $\mathrm{ad}(\mathfrak{x})(\mathfrak{z}) = \mathfrak{x}\mathfrak{z} = -\mathfrak{z}\mathfrak{x}$ it follows from iv. that $\mathrm{ad}(\mathfrak{l}(J))(U(\mathfrak{n}_J^+) \cap \mathfrak{b}) \subseteq \mathfrak{b}$. But $U(\mathfrak{n}_J^+)$ is $\mathrm{ad}(\mathfrak{l}(J))$-invariant. Hence $U(\mathfrak{n}_J^+) \cap \mathfrak{b}$ is $\mathrm{ad}(\mathfrak{l}(J))$-invariant, too. □

We have $M_J \subseteq \mathfrak{b}_J/\mathfrak{b}_J^> \subseteq U(\mathfrak{n}_J^+) + \mathfrak{b}_J^>/\mathfrak{b}_J^>$. In fact $L_\mu \in U(\mathfrak{n}_J^+)$ for $\mu \in B(J)$. Let $M_J^\circ \subseteq U(\mathfrak{n}_J^+)$ denote the preimage of $M_J$ under the projection map $U(\mathfrak{n}_J^+) \longrightarrow U(\mathfrak{g})/\mathfrak{b}_J^>$.

**Lemma 7:**

i. $M_J^\circ \cdot \mathfrak{l}(J) \subseteq \mathfrak{b}_J^>$;
ii. $M_J^\circ$ and $M_J^\circ \cap \mathfrak{b}_J^> = M_J^\circ \cap \mathfrak{b}$ are $\mathrm{ad}(\mathfrak{l}(J))$-invariant;
iii. $\mathrm{ad}(\mathfrak{x})(\mathfrak{z}) = \mathfrak{x}\mathfrak{z} \bmod \mathfrak{b}_J^>$ for $\mathfrak{x} \in \mathfrak{l}(J)$ and $\mathfrak{z} \in M_J^\circ$.

Proof: i. Because of Prop. 6 iii. and iv. it suffices to show that $L_\mu L_{k\ell} \in \mathfrak{b}_J^>$ whenever $\mu \in B(J)$ and $k, \ell \notin J$. If $k \notin J(-\mu)$ then $L_\mu L_{k\ell}$ after sorting coincides up to a constant with some $L_\nu$ such that $J = J(\mu) \subsetneq J(\nu)$; hence $L_\mu L_{k\ell} \in \mathfrak{b}_J^>$ in this case. If $k \in J(-\mu)$ then $L_\mu$ has a factor $L_{ik}$ and since the factors of the monomial $L_\mu$ commute with one another we may use a cancellation relation to conclude that $L_\mu L_{k\ell} \in \mathfrak{b}$.

ii. Using Prop. 6 iii. and v. we are reduced to showing that $\mathrm{ad}(L_{k\ell})(L_\mu) \in M_J^\circ$ whenever $\mu \in B(J)$ and $k, \ell \notin J$. The monomial $L_\mu$ is of the form $L_\mu = c \cdot \prod_{i \in J} L_{is_i}$ with $s_i \notin J$ and some nonzero integer $c$. We have

$$[L_{k\ell}, L_{is_i}] = \begin{cases} -L_{i\ell} & \text{if } k = s_i, \\ 0 & \text{if } k \neq s_i. \end{cases}$$

Since $\mathrm{ad}(L_{k\ell})$ is a derivation it follows that

$$\mathrm{ad}(L_{k\ell})(L_\mu) = -c \cdot \sum_{\substack{i \in J \\ s_i = k}} L_{i\ell} \prod_{\substack{j \in J \\ j \neq i}} L_{js_j}$$

which clearly lies in $M_J^\circ$.

iii. This is an immediate consequence of the first assertion. $\square$

The last lemma shows that the structure of $M_J$ as an $\mathfrak{l}(J)$-module is induced by the adjoint action of $\mathfrak{l}(J)$ on $M_J^\circ$. Whenever convenient we will use all the notations introduced above also for the empty set $J = \emptyset$; all the above assertions become trivially true in this case.

## 5. Local duality

In this section, we study linear forms on the Banach space $\Omega_b^d(U^0)$ of bounded differential forms on the admissible open set $U^0 = r^{-1}(\overline{C}^0)$ which is the inverse image of the open standard chamber in $\overline{X}$ under the reduction map. The restriction map gives a continuous injection from $\Omega^d(\mathcal{X})$ into this Banach space, and therefore linear forms on $\Omega_b^d(U^0)$ are also elements of $\Omega^d(\mathcal{X})'$.

Our first principal result of this section identifies $\Omega_b^d(U^0)$ with the dual of the space $\mathcal{O}(B)^{\mathfrak{b}=0}$ of (globally) analytic functions on $B$ which are annihilated by the ideal $\mathfrak{b}$ studied in the preceeding section. The filtration which we introduced on $U(\mathfrak{g})/\mathfrak{b}$ then yields filtrations of $\mathcal{O}(B)^{\mathfrak{b}=0}$ and $\Omega_b^d(U^0)$. Applying our analysis of the subquotients of the filtration on $U(\mathfrak{g})/\mathfrak{b}$ from the preceeding section, we describe each subquotient of the filtration on $\mathcal{O}(B)^{\mathfrak{b}=0}$ as a space of analytic vector-valued functions on the unipotent radical of a specific maximal parabolic subgroup in $G$ satisfying certain explicit differential equations.

Of fundamental importance to this analysis are the linear forms arising from the residue map on the standard chamber.

The space $\Omega_b^d(U^0)$ of bounded $d$-forms $\eta$ on $U^0$ are those which have an expansion
$$\eta = \sum_{\nu \in X^*(\overline{T})} a(\nu)\Xi_\nu d\Xi_{\alpha_{d-1}} \wedge \ldots \wedge d\Xi_{\alpha_0}$$
such that
$$\omega_C(\eta) := \inf_\nu \{\omega(a(\nu)) - \ell(\nu)\} > -\infty \ .$$

We may and will always view $\Omega_b^d(U^0)$ as a Banach space with respect to the norm $\omega_C$ (compare [ST] Remark after Lemma 17). According to Lemma 1.2 the restriction map induces a continuous injective map $\Omega^d(\mathcal{X}) \longrightarrow \Omega_b^d(U^0)$. In [ST] Def. 19 we defined the residue of $\eta \in \Omega_b^d(U^0)$ at the pointed chamber $(\overline{C}, 0)$ by
$$\mathrm{Res}_{(\overline{C},0)} \eta := a(\alpha_d) \ .$$

It is then clear that, for any weight $\mu \in X^*(\overline{T})$,
$$\eta \longmapsto \mathrm{Res}_{(\overline{C},0)} \Xi_{-\mu} \eta = a(\alpha_d + \mu)$$
is a continuous linear form on $\Omega_b^d(U^0)$ and a fortiori on $\Omega^d(\mathcal{X})$. Applying the map $I$ we obtain the locally analytic function
$$f_\mu(g) := \mathrm{Res}_{(\overline{C},0)}(\Xi_{-\mu} \cdot g_*\xi)$$
on $G$. We collect the basic properties of these functions.

1. Under the adjoint action of $B \cap T$ the function $f_\mu$ has weight $-\mu$, i.e.,
$$f_\mu(t^{-1}gt) = \mu(t^{-1}) \cdot f_\mu(g) \text{ for } g \in G \text{ and } t \in B \cap T \ .$$

This is straightforward from [ST] Lemma 20.

2. The restriction $f_\mu|B$ of $f_\mu$ to the Iwahori subgroup $B \subseteq G$ is analytic on $B$. First of all recall that $B$ is a product of disks and annuli where the matrix entries $g_{ij}$ of $g \in B$ can be used as coordinates (the diagonal entries correspond to the annuli). By construction as well as by the formula
$$d\Xi_{\alpha_{d-1}} \wedge \ldots \wedge d\Xi_{\alpha_0} = (-1)^{d(d+1)/2} \Xi_{-\beta-\alpha_d} d\Xi_{\beta_0} \wedge \ldots \wedge d\Xi_{\beta_{d-1}}$$
we have, for a fixed $g \in G$, the expression

(a) $$(g_*\xi)|U^0 = (-1)^{d(d+1)/2} \sum_{\mu \in X^*(\overline{T})} f_\mu(g) \Xi_\mu \xi$$

in $\Omega_b^d(U^0)$. This is, of course, not a convergent expansion with respect to the norm $\omega_C$. But if we write $g_*\xi|U^0 = F(g)\xi|U^0$ then the series

$$F(g) = (-1)^{d(d+1)/2} \sum_{\mu \in X^*(\overline{T})} f_\mu(g) \Xi_\mu$$

is uniformly convergent on each affinoid subdomain of $U^0$.
On the other hand a direct calculation shows that

$$g_*\xi = \det(g) \left( \prod_{j=0}^{d} \frac{1}{f_j(g,.)} \right) d\Xi_{\beta_0} \wedge \ldots \wedge d\Xi_{\beta_{d-1}}$$

where

$$f_j(g,q) := \sum_{i=0}^{d-1} g_{ij} \Xi_{\beta_i} + g_{dj} \ .$$

Recall that $U^0$ is given by the inequalities

$$\omega(\Xi_0(q)) < \ldots < \omega(\Xi_d(q)) < 1 + \omega(\Xi_0(q)) \ .$$

It follows that for $g \in B$ the term $g_{jj}\Xi_{\beta_j}$ in the sum $f_j(g,q)$ is strictly larger in valuation than the other terms (we temporarily put $\beta_d := 0$). We therefore have, for $g \in B$ and $q \in U^0$, the geometric series expansion

$$\frac{1}{f_j(g,q)} = \frac{1}{g_{jj}\Xi_{\beta_j}} \sum_{m \geq 0} \left( -\sum_{\substack{i=0 \\ i \neq j}}^{d} \frac{g_{ij}}{g_{jj}} \Xi_{\varepsilon_i - \varepsilon_j} \right)^m \ .$$

If we multiply those expansions together and compare the result to (a) we obtain the expansion

(b) $$f_\mu(g) = \frac{\det(g)}{g_{00} \cdot \ldots \cdot g_{dd}} \cdot \sum_{\underline{m} \in I(\mu)} c_{\underline{m}} \cdot \prod_{i \neq j} \left( \frac{g_{ij}}{g_{jj}} \right)^{m_{ij}}$$

where the $c_{\underline{m}}$ are certain nonzero integers (given as a sign times a product of polynomial coefficients) and

$$I(\mu) := \text{ set of all tuples } \underline{m} = (m_{ij})_{i \neq j} \text{ consisting} \\ \text{ of integers } m_{ij} \geq 0 \text{ such that} \\ \mu = \sum_{i \neq j} m_{ij}(\varepsilon_i - \varepsilon_j) \ .$$

In order to see that this expansion actually is uniformly convergent in $g \in B$ let

$$\pi(\underline{m}) := \sum_{i < j} m_{ij} \ .$$

It is clear that if we fix $\mu$ and an $n \geq 0$ then the number of $\underline{m} \in I(\mu)$ such that $\pi(\underline{m}) = n$ is finite. But on the other hand, for $g \in B$, the matrix entries $g_{ij}$ for $i < j$ are divisible by $\pi$. Hence the valuation of the summand corresponding to the tuple $\underline{m}$ in the expansion (b) is at least $\pi(\underline{m})$.

3. The restriction of $f_\mu$ to $B$ does not vanish identically. In order to see this we make a choice of simple roots $\alpha'_0, \ldots, \alpha'_{d-1}$ with respect to which $\mu$ is positive, i.e., $\mu = n_0 \alpha'_0 + \ldots + n_{d-1} \alpha'_{d-1}$ with $n_i \geq 0$. Consider the matrix $g_0 \in B$ which has a 1 on all diagonal positions, a $\pi$ on the positions $\alpha'_0, \ldots, \alpha'_{d-1}$, and 0 elsewhere. Then

$$f_\mu(g_0) = c\pi^{n_0 + \ldots + n_{d-1}} \text{ with some nonzero } c \in \mathbb{Z} \ .$$

Let $\omega_B$ denote the spectral norm on the affinoid algebra $\mathcal{O}(B)$ of $K$-analytic functions on $B$. We have to determine the precise value of $\omega_B(f_\mu|B)$.

**Lemma 1:**

*For any $\underline{m} \in I(\mu)$ we have $\pi(\underline{m}) \geq \ell(\mu)$.*

Proof: Recall ([ST] p. 405) that

$$\ell(\mu) = -\inf\nolimits_{z \in \overline{C}} \mu(z) \ .$$

It follows that $\ell(\mu + \nu) \leq \ell(\mu) + \ell(\nu)$ holds for any $\mu, \nu \in X^*(\overline{T})$. Hence if $\mu = \sum_{i \neq j} m_{ij}(\varepsilon_i - \varepsilon_j)$ then we have

$$\ell(\mu) \leq \sum_{i \neq j} m_{ij} \ell(\varepsilon_i - \varepsilon_j) \ .$$

It therefore suffices to check that

$$\ell(\varepsilon_i - \varepsilon_j) \leq \begin{cases} 1 & \text{if } i < j, \\ 0 & \text{if } i > j. \end{cases}$$

But that is obvious from the definition of the chamber $C$.

4. We claim that
$$\omega_B(f_\mu|B) = \ell(\mu)$$
holds true. The norm $\omega_B$ on $\mathcal{O}(B)$ is multiplicative and the first factor $\det(g) \cdot (g_{00} \cdot \ldots \cdot g_{dd})^{-1}$ in the expansion (b) is a unit in $\mathcal{O}(B)$. It therefore follows from the lemma that $\omega_B(f_\mu|B) \geq \ell(\mu)$ and that it suffices to find an $\underline{m} \in I(\mu)$ such that
$$\omega(c_{\underline{m}}) + \sum_{i<j} m_{ij} = \ell(\mu) \ .$$

Let us first consider the special case where $\ell(\mu) = 0$. Then $\mu = n_0\alpha_0 + \ldots + n_{d-1}\alpha_{d-1}$ with all $n_i \geq 0$. Consider the element $h = (h_{ij}) \in B$ where $h_{ii} = 1$ for $0 \leq i \leq d$ and $h_{i+1,i} = -1$ for $0 \leq i \leq d-1$, with all other $h_{ij} = 0$. For this matrix, we compute

$$h_*\xi = (1 - \tfrac{\Xi_1}{\Xi_0})^{-1} \ldots (1 - \tfrac{\Xi_d}{\Xi_{d-1}})^{-1}\xi$$

$$= \sum_{(i_0,\ldots,i_d)} (\tfrac{\Xi_1}{\Xi_0})^{i_0} \ldots (\tfrac{\Xi_d}{\Xi_{d-1}})^{i_d}\xi$$

in $\Omega_b^d(U^0)$. Using (a) this shows that $f_\mu(h) = \pm 1$. On the other hand substituting $g = h$ in the series expansion (b) we see that $1 = \pm f_\mu(h) = \pm c_{\underline{m}}$ for the particular $\underline{m} \in I(\mu)$ corresponding to the representation $\mu = \sum_{i=0}^{d-1} n_i\alpha_i$ – that is, the $\underline{m}$ with $m_{i+1,i} = n_i$ and other $m_{ij} = 0$. This implies that

$$\ell(\mu) = 0 = \omega(c_{\underline{m}}) + \sum_{i<j} m_{ij}$$

in this case. In order to treat the general case we first make the following observations. Let $s_i \in G$ denote the permutation matrix which represents the reflection in the Weyl group corresponding to the simple root $\alpha_i$. The Coxeter element $s = s_0 \ldots s_{d-1}$ permutes the roots $\alpha_0, \ldots, \alpha_d$ cyclically. The same then is true for the element $\rho := ys$ where $y$ denotes the diagonal matrix in $G$ with entries $\pi, 1, \ldots, 1$. But $\rho$ normalizes the subgroup $B$ and in particular changes the residue of a $d$-form only by a sign ([ST] Thm. 24). Let now

$$\nu = n_0\alpha_0 + \ldots + n_d\alpha_d \text{ with all } n_i \geq 0 \text{ and } n_a = 0 \text{ for some } 0 \leq a \leq d$$

be any weight; in particular $\ell(\nu) = n_d$. The weight $\mu := \rho^{d-a}(\nu)$ then satisfies $\ell(\mu) = 0$. Defining $h \in B$ as before we have $f_\mu(h) = \pm 1$. The matrix $(h'_{ij}) = h' := \rho^{a-d}h\rho^{d-a} \in B$ is given by

$$h'_{ii} = 1,\ h'_{i+1,i} = -1 \text{ for } i \neq a,\ h'_{0d} = -\pi,\ \text{and all other } h'_{ij} = 0 \ .$$

Substituting $g = h'$ in (b) we obtain

$$f_\nu(h') = \pm c_{\underline{n}} \cdot \pi^{n_{0d}}$$

where $\underline{n} \in I(\nu)$ corresponds to the above representation of $\nu$ (in particular, $n_{0d} = n_d$). On the other hand we compute

$$\begin{aligned}
\omega(f_\nu(h')) &= \omega(\operatorname{Res}_{(\overline{C},0)} \Xi_{-\nu} \cdot (\rho^{a-d} h \rho^{d-a})_* \xi) \\
&= \omega(\operatorname{Res}_{(\overline{C},0)} \Xi_{-\nu} \cdot (\rho^{a-d} h)_* \xi) \\
&= \omega(\operatorname{Res}_{(\overline{C},0)} (\rho^{d-a})_* \Xi_{-\nu} \cdot h_* \xi) \\
&= \omega(\prod_{i=1}^{d-a} s^i(\nu)(y^{-1}) \cdot \operatorname{Res}_{(\overline{C},0)} \Xi_{-\mu} \cdot h_* \xi) \\
&= -\sum_{i=1}^{d-a} \omega(s^i(\nu)(y)) + \omega(f_\mu(h)) \\
&= -(n_{d-1} - n_d) - (n_{d-2} - n_{d-1}) - \ldots - (n_a - n_{a+1}) \\
&= n_d = \ell(\nu) \ .
\end{aligned}$$

5. Since $B$ is open in $G$ the enveloping algebra $U(\mathfrak{g})$ also acts by left invariant differential operators on $\mathcal{O}(B)$. It is an immediate consequence of the definition that

$$f_\mu | B \in \mathcal{O}(B)^{\mathfrak{b}=0} \ .$$

**Proposition 2:**

*For any $\mu \in X^*(\overline{T})$ the weight space of weight $-\mu$ in $\mathcal{O}(B)^{\mathfrak{b}=0}$ with respect to the adjoint action of $B \cap T$ is the 1-dimensional subspace generated by $f_\mu|B$.*

Proof: We consider the pairing

$$\begin{array}{ccc}
U(\mathfrak{g})/\mathfrak{b} \times \mathcal{O}(B)^{\mathfrak{b}=0} & \longrightarrow & K \\
(\mathfrak{z}, f) & \longmapsto & (\mathfrak{z} f)(1) \ .
\end{array}$$

It is nondegenerate on the right by Taylor's formula. It also is invariant with respect to the adjoint action of $B \cap T$ on both sides. Hence the induced map

$$\mathcal{O}(B)^{\mathfrak{b}=0} \hookrightarrow \operatorname{Hom}_K(U(\mathfrak{g})/\mathfrak{b}, K)$$

is injective and respects weight spaces. It then follows from Corollary 4.5 that the weight spaces on the left hand side are at most 1-dimensional. But we know that $f_\mu|B$ is nonvanishing. □

The meaning of that proposition is that any function $f \in \mathcal{O}(B)^{\mathfrak{b}=0}$ has an expansion of the form

$$f = \sum_{\mu \in X^*(\overline{T})} b(\mu)(f_\mu|B)$$

with $b(\mu) \in K$ such that $\omega(b(\mu)) + \ell(\mu) \longrightarrow \infty$ with respect to the Fréchet filter of complements of finite subsets in $X^*(\overline{T})$. First expand $f$ into a series in the matrix entries and then collect all terms of a specific weight $-\mu$. We obtain in this way an expansion

$$f = \sum_{\mu} \tilde{f}_\mu \text{ with } \omega_B(\tilde{f}_\mu) \longrightarrow \infty .$$

Since the $U(\mathfrak{g})$-action on $\mathcal{O}(B)$ is by continuous endomorphisms and since the ideal $\mathfrak{b}$ is homogeneous in the weight space decomposition the equation $\mathfrak{b}f = 0$ implies $\mathfrak{b}\tilde{f}_\mu = 0$ for any $\mu$. It therefore follows from the proposition that $\tilde{f}_\mu = b(\mu)(f_\mu|B)$ for some $b(\mu) \in K$.
If we now consider a $d$-form

$$\eta = \sum_{\mu \in X^*(\overline{T})} a(\mu) \Xi_\mu \xi \in \Omega_\mathfrak{b}^d(U^0)$$

then we see that

$$\langle \eta, f \rangle := \sum_{\mu} a(\mu) b(\mu)$$

converges in $K$. In this way we obtain a bilinear pairing

$$\langle \, , \, \rangle : \Omega_\mathfrak{b}^d(U^0) \times \mathcal{O}(B)^{\mathfrak{b}=0} \longrightarrow K .$$

Actually the following stronger statement is immediately clear.

**Proposition 3:**

*The pairing $\langle \, , \, \rangle$ induces a topological isomorphism*

$$[\mathcal{O}(B)^{\mathfrak{b}=0}]' = \Omega_\mathfrak{b}^d(U^0) .$$

The connection between this local duality and the map $I$ from the second section is provided by the diagram

$$\begin{array}{ccc}
\Omega^d(\mathcal{X})' & \xrightarrow{I} & C^{\mathrm{an}}(G,K)^{\mathfrak{a}=0} \\
& & \downarrow \text{restriction} \\
(\text{restriction})' \uparrow & & C^{\mathrm{an}}(B,K)^{\mathfrak{b}=0} \\
& & \uparrow \subseteq \\
\Omega_\mathfrak{b}^d(U^0)' & \xleftarrow{\langle , \rangle} & \mathcal{O}(B)^{\mathfrak{b}=0}
\end{array}$$

which, by the very construction of the above pairing, is commutative up to sign. The ideal filtration $\mathfrak{b} \subseteq \ldots \subseteq \mathfrak{b}_j \subseteq \ldots \subseteq U(\mathfrak{g})$ gives rise to a filtration

$$\mathcal{O}(B)^{\mathfrak{b}=0} \supseteq \ldots \supseteq \mathcal{O}(B)^{\mathfrak{b}_{j+1}=0} \supseteq \mathcal{O}(B)^{\mathfrak{b}_j=0} \supseteq \ldots \supseteq \mathcal{O}(B)^{\mathfrak{b}_0=0} = \{0\}$$

as well as, by duality, to a "local" filtration

$$\Omega_b^d(U^0) = \Omega_b^d(U^0)^0 \supseteq \ldots \supseteq \Omega_b^d(U^0)^j \supseteq \ldots \supseteq \Omega_b^d(U^0)^{d+1} = \{0\}$$

with

$$\Omega_b^d(U^0)^j := [\mathcal{O}(B)^{\mathfrak{b}=0}/\mathcal{O}(B)^{\mathfrak{b}_j=0}]'\ .$$

We need to understand how the properties of the ideal filtration which we have established in the previous section translate into properties of the other filtrations. Recall that the "bases" $\{f_\mu\}$ of $\mathcal{O}(B)^{\mathfrak{b}=0}$ and $\{L_\mu\}$ of $U(\mathfrak{g})/\mathfrak{b}$, respectively, are "dual" to each other in the sense that

$$L_\nu f_\mu(1) = \begin{cases} \pm 1 & \text{for } \nu = \mu\ , \\ 0 & \text{for } \nu \neq \mu\ . \end{cases}$$

If $f \in \mathcal{O}(B)^{\mathfrak{b}=0}$ has the expansion $f = \sum_\mu b(\mu)(f_\mu|B)$ we therefore have

$$(L_\nu f)(1) = \sum_\mu b(\mu)(L_\nu f_\mu)(1) = \pm b(\nu)\ .$$

¿From this one easily deduces that

- $\mathcal{O}(B)^{\mathfrak{b}_J=0} = \{f \in \mathcal{O}(B)^{\mathfrak{b}=0} : f = \sum_{J \not\subseteq J(\mu)} b(\mu)(f_\mu|B)\}$;

- any coset in $\mathcal{O}(B)^{\mathfrak{b}_J^>=0}/\mathcal{O}(B)^{\mathfrak{b}_J=0}$ has a unique representative of the form

$$f = \sum_{J(\mu)=J} b(\mu)(f_\mu|B)\ .$$

This leads to the fact that the map

$$\mathcal{O}(B)^{\mathfrak{b}_{j+1}=0}/\mathcal{O}(B)^{\mathfrak{b}_j=0} \xrightarrow{\cong} \bigoplus_{\#J=j} \mathcal{O}(B)^{\mathfrak{b}_J^>=0}/\mathcal{O}(B)^{\mathfrak{b}_J=0}$$

$$f = \sum_{\#J(\mu)=j} b(\mu)(f_\mu|B) \longmapsto (\sum_{J(\mu)=J} b(\mu)(f_\mu|B))_J$$

is a continuous linear isomorphism. We will give a reinterpretation of the right hand side which reflects the fact that $\mathfrak{b}_J/\mathfrak{b}_J^>$ is a quotient of the generalized Verma module $U(\mathfrak{g}) \otimes_{U(\mathfrak{p}_J)} M_J$ via the map which sends $\mathfrak{z} \otimes m$ to $\mathfrak{z}m$. Set

$$\mathfrak{d}_J := \ker(U(\mathfrak{g}) \otimes_{U(\mathfrak{p}_J)} M_J \longrightarrow \mathfrak{b}_J/\mathfrak{b}_J^>)\ .$$

By the Poincaré-Birkhoff-Witt theorem the inclusion $U(\mathfrak{n}_J^+) \subseteq U(\mathfrak{g})$ induces an isomorphism $U(\mathfrak{n}_J^+) \underset{K}{\otimes} M_J \xrightarrow{\cong} U(\mathfrak{g}) \underset{U(\mathfrak{p}_J)}{\otimes} M_J$. In this section we always will view $\mathfrak{d}_J$ as a subspace of $U(\mathfrak{n}_J^+) \underset{K}{\otimes} M_J$.

Let $U_J^+$ be the unipotent subgroup in $G$ whose Lie algebra is $\mathfrak{n}_J^+$, and let $\mathcal{O}(U_J^+ \cap B)$ denote the $K$-affinoid algebra of $K$-analytic functions on the polydisk $U_J^+ \cap B$. Consider the pairing

$$\langle\,,\,\rangle : (U(\mathfrak{n}_J^+) \underset{K}{\otimes} M_J) \times (\mathcal{O}(U_J^+ \cap B) \underset{K}{\otimes} M_J') \longrightarrow \mathcal{O}(U_J^+ \cap B)$$

$$(\mathfrak{z} \otimes m, e \otimes E) \longmapsto E(m) \cdot \mathfrak{z}e$$

and define the Banach space

$$\mathcal{O}(U_J^+ \cap B, M_J')^{\mathfrak{d}_J=0} := \{\varepsilon \in \mathcal{O}(U_J^+ \cap B) \underset{K}{\otimes} M_J' : \langle \mathfrak{d}_J, \varepsilon \rangle = 0\}\,.$$

Let also $\{L_\mu^*\}_{\mu \in B(J)}$ denote the basis of $M_J'$ dual to the basis $\{L_\mu\}_\mu$ of $M_J$.

**Proposition 4:**

*The map*

$$\nabla_J : \mathcal{O}(B)^{\mathfrak{b}_J^>=0} / \mathcal{O}(B)^{\mathfrak{b}_J=0} \xrightarrow{\cong} \mathcal{O}(U_J^+ \cap B, M_J')^{\mathfrak{d}_J=0}$$

$$f \longmapsto \sum_{\mu \in B(J)} [(L_\mu f) | U_J^+ \cap B] \otimes L_\mu^*$$

*is an isomorphism of Banach spaces.*

Proof: For $\mathfrak{z} = \sum_\nu \mathfrak{z}_{(\nu)} \otimes L_\nu \in \mathfrak{d}_J \subseteq U(\mathfrak{n}_J^+) \underset{K}{\otimes} M_J$ we have

$$\langle \mathfrak{z}, \sum_\mu [(L_\mu f)|U_J^+ \cap B] \otimes L_\mu^* \rangle = \sum_{\mu,\nu} L_\mu^*(L_\nu) \cdot (\mathfrak{z}_{(\nu)} L_\mu f)|U_J^+ \cap B$$

$$= (\sum_\nu \mathfrak{z}_{(\nu)} L_\nu) f | U_J^+ \cap B = 0$$

since $\sum_\nu \mathfrak{z}_{(\nu)} L_\nu \in U(\mathfrak{n}_J^+) \cap \mathfrak{b}_J^>$. Morover for $\mu \in B(J)$ we have $L_\mu \in \mathfrak{b}_J$. Hence the map $\nabla_J$ is well defined. It clearly is continuous. The Banach space on the left hand side of the assertion has the orthonormal basis $\pi^{-\ell(\nu)} f_\nu | B$ for $J(\nu) = J$. Concerning the right hand side we observe that the above pairing composed with the evaluation in 1 induces an injection

$$\mathcal{O}(U_J^+ \cap B) \underset{K}{\otimes} M_J' \hookrightarrow \mathrm{Hom}_K(U(\mathfrak{n}_J^+) \underset{K}{\otimes} M_J, K)$$

which restricts to an injection

$$\mathcal{O}(U_J^+ \cap B, M_J')^{\mathfrak{d}_J=0} \hookrightarrow \mathrm{Hom}_K(\mathfrak{b}_J/\mathfrak{b}_J^>, K) \ .$$

Hence the only weights which can occur in the right hand side are those $\nu$ with $J(-\nu) = J$ and the corresponding weight spaces are at most 1-dimensional. Moreover the same argument as after Prop. 2 shows that the occurring weight vectors (scaled appropriately) form an orthonormal basis. Since $\nabla_J$ visibly preserves weights the assertion follows once we show that

$$\nabla_J(f_\nu|B) \neq 0 \ \text{ for any } \nu \text{ with } J(\nu) = J \ .$$

All that remains to be checked therefore is the existence, for a given $\nu$ with $J(\nu) = J$, of a $\mu \in B(J)$ such that $L_\mu f_\nu$ dies not vanish identically on $U_J^+ \cap B$.

The weight $\nu$ is of the form $\nu = \sum_{j=0}^{d} n_j \varepsilon_j$ with $n_j > 0$ for $j \in J$ and $n_j \leq 0$ for $j \notin J$. We have

$$\#J \leq \sum_{j \in J} n_j = -\sum_{j \notin J} n_j \ .$$

Choose integers $n_j \leq m_j \leq 0$ for $j \notin J$ such that $\#J = -\sum_{j \notin J} m_j$ and define

$$\mu := \sum_{j \in J} \varepsilon_j + \sum_{j \notin J} m_j \varepsilon_j \in B(J) \ .$$

Observe that $J(\nu-\mu) \subseteq J$ and $J(\mu-\nu) \cap J = \emptyset$. This means that $L_{\nu-\mu} \in U(\mathfrak{n}_J^+)$. It suffices to check that $L_{\nu-\mu} L_\mu f_\nu(1) \neq 0$. We compute

$$L_{\nu-\mu} L_\mu f_\nu(1) = \mathrm{Res}_{(\overline{C},0)} \Xi_{-\nu} \cdot L_{\nu-\mu} L_\mu \xi$$

$$= -\mathrm{Res}_{(\overline{C},0)} \Xi_{-\nu} \cdot L_{\nu-\mu}(\Xi_\mu \xi) \ .$$

As a consequence of the formula (+) in section 4 we have $L_{\nu-\mu}(\Xi_\mu \xi) = m \cdot \Xi_\nu \xi$ for some nonzero integer $m$. Hence we obtain

$$L_{\nu-\mu} L_\mu f_\nu(1) = -m \cdot \mathrm{Res}_{(\overline{C},0)} \xi = \pm m \neq 0 \ . \qquad \square$$

As a consequence of this discussion we in particular have a unique continuous linear map

$$D_J : \mathcal{O}(U_J^+ \cap B, M_J')^{\mathfrak{d}_J=0} \longrightarrow [\Omega_b^d(U^0)^j / \Omega_b^d(U^0)^{j+1}]' \ ,$$

where $j := \#J$, which sends the weight vector $\sum_{\mu \in B(J)} [(L_\mu f_\nu)|U_J^+ \cap B] \otimes L_\mu^*$, for $\nu$ with $J(\nu) = J$, to the linear form $\lambda_\nu(\eta) := \mathrm{Res}_{(\overline{C},0)} \Xi_{-\nu}\eta$.

## 6. The global filtration

In this section, we find a $G$-invariant filtration on the full space $\Omega^d(\mathcal{X})$ that is compatible with the local filtration discussed in the previous section. This "global" filtration is defined first on the subspace of $\Omega^d(\mathcal{X})$ consisting of algebraic $d$-forms having poles along a finite set of $K$-rational hyperplanes; the filtration on the full space is obtained by passing to the closure. We obtain at the same time a filtration on the dual space $\Omega^d(\mathcal{X})'$. A key tool in our description of this filtration is a "partial fractions decomposition" due to Gelfand and Varchenko.

At the end of the section, we apply general results from the theory of topological vector spaces (in the non-archimedean situation) to show that the subquotients of the global filtration on $\Omega^d(\mathcal{X})$ are reflexive Fréchet spaces whose duals can be computed by the subquotients of the dual filtration.

Let us first recall some general notions from algebraic geometry. Let $\mathcal{L}$ be an invertible sheaf on $\mathbb{P}^d_{/K}$. With any regular meromorphic section $s$ of $\mathcal{L}$ over $\mathbb{P}^d_{/K}$ we may associate a divisor $\mathrm{div}(s)$ (compare EGA IV.21.1.4). One has $\mathrm{div}(s') = \mathrm{div}(s)$ if and only if $s' = ts$ for some invertible regular function $t$ on $\mathbb{P}^d$. Let $\{Y_i\}_{i \in I}$ be the collection of prime divisors on $\mathbb{P}^d_{/K}$ and write

$$\mathrm{div}(s) = \sum_{i \in I} n_i Y_i$$

where almost all of the integers $n_i$ are zero. One has

$$\sum_i n_i = n \text{ if } \mathcal{L} \cong \mathcal{O}(n)$$

([Har] II.6.4). We put

$$\mathrm{div}(s)_\infty := -\sum_{\substack{i \\ n_i < 0}} n_i Y_i$$

and

$$\imath_o(s) := \#\{i \in I : n_i < 0\} \ .$$

By convention let $\mathrm{div}(0)_\infty := 0$ and $\imath_o(0) = 0$. We want to apply these notions in the case of the canonical invertible sheaf $\mathcal{L} = \Omega^d \cong \mathcal{O}(-d-1)$ on $\mathbb{P}^d_{/K}$. A regular meromorphic global section $\eta$ in this case is a $d$-form $\eta = F\xi$ such that $F$ is a nonzero rational function on $\mathbb{P}^d_{/K}$. We will study the subspace

$$\Omega^d_{\text{alg}}(\mathcal{X}) := \quad \text{all regular meromorphic global sections}$$
$$\eta \text{ of } \Omega^d \text{ such that } \text{div}(\eta)_\infty \text{ is supported on}$$
$$\text{a union of } K\text{-rational hyperplanes in } \mathbb{P}^d$$
$$\text{together with the zero section}$$

of "algebraic forms" in $\Omega^d(\mathcal{X})$. For any $\eta \in \Omega^d_{\text{alg}}(\mathcal{X})$ we introduce its index as being the nonnegative integer

$$\imath(\eta) := \min \max_k \imath_o(\eta_k)$$

where the minimum is taken over all representations $\eta = \sum_k \eta_k$ of $\eta$ as a finite sum of other $\eta_k = \Omega^d_{\text{alg}}(\mathcal{X})$. By definition we have

$$\imath(\eta + \eta') \leq \max(\imath(\eta), \imath(\eta')) .$$

Hence $\Omega^d_{\text{alg}}(\mathcal{X})$ is equipped with the filtration

$$\ldots \supseteq \Omega^d_{\text{alg}}(\mathcal{X})^0 \supseteq \ldots \supseteq \Omega^d_{\text{alg}}(\mathcal{X})^d \supseteq \Omega^d_{\text{alg}}(\mathcal{X})^{d+1} = \{0\}$$

by the subspaces

$$\Omega^d_{\text{alg}}(\mathcal{X})^j := \{\eta \in \Omega^d_{\text{alg}}(\mathcal{X}) \,:\, \imath(\eta) \leq d + 1 - j\} .$$

**Lemma 1:**

*The index $\imath(\eta)$ is $G$-invariant and takes values between $1$ and $d+1$ for all nonzero $\eta \in \Omega^d_{\text{alg}}(\mathcal{X})$.*

Proof: The $G$-invariance is clear since $G$ preserves $K$-rational hyperplanes. The upper bound for the index follows from the existence of a partial fraction decomposition ([GV] Thm. 21) which says that $\Omega^d_{\text{alg}}(\mathcal{X})$ as a vector space is spanned by the forms $u_*(\Xi_\mu \xi) = (u_* \Xi_{\mu-\beta}) d\Xi_{\beta_0} \wedge \ldots \wedge d\Xi_{\beta_{d-1}}$ with $\mu \in X^*(\overline{T})$ and $u \in P$ unipotent. Each $\Xi_\mu \xi$ has poles along at most the $d+1$ coordinate hyperplanes defined by the equations $\Xi_i = 0$ for $i = 0, \ldots, d$. □

It follows that the subspace $\Omega^d_{\text{alg}}(\mathcal{X})$ together with its filtration is $G$-invariant. Moreover the filtration is finite with $\Omega^d_{\text{alg}}(\mathcal{X}) = \Omega^d_{\text{alg}}(\mathcal{X})^0$. In order to obtain finer information we need to take a closer look at that partial fraction decomposition. First we introduce, for any subset $J \subseteq \{0, \ldots, d\}$, the subgroup

$$U(J) := \quad \text{all lower triangular unipotent matrices } u = (u_{ij})$$
$$\text{such that } u_{ij} = 0 \text{ whenever } i > j \text{ and } j \in J$$

of $U$. In particular $U(\{0,\ldots,d\}) = U(\{0,\ldots,d-1\}) = \{1\}$ and $U(\emptyset) = U(\{d\}) = U$.

**Proposition 2:**

*Every differential form $\eta \in \Omega^d_{\mathrm{alg}}(\mathcal{X})$ may be written as a sum*

$$\eta = \sum_{J \subseteq \{0,\ldots,d\}} \sum_{\substack{\mu \in X^*(\overline{T}) \\ J(\mu) = J}} \sum_{u \in U(J)} A(\mu, u) u_*(\Xi_\mu \xi)$$

*where the coefficients $A(\mu, u) \in K$ are zero for all but finitely many pairs $(\mu, u)$; furthermore, such an expression is unique.*

Proof: We write $\eta = F d\Xi_{\beta_0} \wedge \ldots \wedge d\Xi_{\beta_{d-1}}$ and apply that partial fraction decomposition to $F$ obtaining an expression

$$F = \sum_{\mu \in X^*(\overline{T})} \sum_{u \in U} B(\mu, u) u_* \Xi_\mu \, .$$

The uniqueness part of that Thm. 21 in [GV] says that such an expression even exists and is unique under the following additional requirement: If $\mu = \sum_k m_k \varepsilon_k$ then we sum only over those $u \in U$ whose $k$-th column consists of zeroes except for the diagonal entry $u_{kk} = 1$ for every $k$ such that $m_k \geq 0$. But $\{k : m_k \geq 0\} = J(\mu + \beta)$ so that the condition on $u$ becomes exactly that $u \in U(J(\mu+\beta))$. Because of $(u_* \Xi_\mu) d\Xi_{\beta_0} \wedge \ldots \wedge d\Xi_{\beta_{d-1}} = u_*(\Xi_{\mu+\beta} \xi)$ we obtain the desired unique expression if we put $A(\mu, u) := B(\mu - \beta, u)$. □

Let us temporarily introduce as another invariant of a form $\eta \in \Omega^d_{\mathrm{alg}}(\mathcal{X})$ the linear subvariety

$$Z(\eta) := \text{ the intersection of all hyperplanes contained in the support of } \mathrm{div}(\eta)_\infty$$

in $\mathbb{P}^d_{/K}$. One obviously has:

- codim $Z(\eta) \leq \imath_{\mathrm{o}}(\eta)$;
- codim $Z(g_*(\Xi_\mu \xi)) = \imath_{\mathrm{o}}(g_*(\Xi_\mu \xi))$ for any $g \in G$ and $\mu \in X^*(\overline{T})$.

Write

$$\eta = F_{\mathrm{hom}}(\Xi_0, \ldots, \Xi_d) \cdot \sum_{i=0}^{d} (-1)^i \Xi_i d\Xi_0 \wedge \ldots \wedge \widehat{d\Xi_i} \wedge \ldots \wedge d\Xi_d$$

as a homogeneous form on affine space $\mathbb{A}^{d+1}$ and apply the partial fraction decomposition in [GV] to $F_{\text{hom}}$. Then, at each stage of the construction of the partial fraction decomposition of $F_{\text{hom}}$, the linear forms occurring in the denominator of any term are linear combinations of those in the denominator of $F_{\text{hom}}$. This means that

- $Z(\eta) \subseteq Z(A(\mu, u) u_*(\Xi_\mu \xi))$.

Together these three observations imply that

$$\imath_{\text{o}}(\eta) \geq \imath_{\text{o}}(A(\mu, u) u_*(\Xi_\mu \xi)) .$$

It then follows from the unicity of the partial fraction decomposition that we actually have

$$\imath(\eta) = \max_{\mu, u} \imath_{\text{o}}(A(\mu, u) u_*(\Xi_\mu \xi)) .$$

But $\text{div}(\xi)_\infty = \sum_{i=0}^{d} \{\Xi_i = 0\}$ and therefore $\imath_{\text{o}}(\Xi_\mu \xi) = d + 1 - \#J(\mu)$. We obtain the following explicit formula

$$\imath(\eta) = \max\{d + 1 - \#J(\mu) : \mu \text{ such that } A(\mu, u) \neq 0 \text{ for some } u \in U(J(\mu))\}$$

for the index of any $d$-form $\eta \neq 0$. Another consequence of this discussion that we will need later is the inequality

$$\imath(\eta) \leq \text{codim}\, Z(\eta) .$$

**Corollary 3:**

$\Omega_{\text{alg}}^d(\mathcal{X})^j$ as a $K$-vector space is spanned by the forms $u_*(\Xi_\mu \xi)$ where $(\mu, u) \in X^*(\overline{T}) \times U$ runs over those pairs for which $u \in U(J(\mu))$ and $\#J(\mu) \geq j$; in particular

$$\Omega_{\text{alg}}^d(\mathcal{X})^j = \sum_{g \in G} g_*(\mathfrak{b}_j \xi) .$$

Our "global" $G$-equivariant filtration

$$\Omega^d(\mathcal{X}) = \Omega^d(\mathcal{X})^0 \supseteq \ldots \supseteq \Omega^d(\mathcal{X})^d \supseteq \Omega^d(\mathcal{X})^{d+1} = \{0\}$$

of $\Omega^d(\mathcal{X})$ now is defined by taking closures

$$\Omega^d(\mathcal{X})^j := \text{ closure of } \Omega_{\text{alg}}^d(\mathcal{X})^j \text{ in } \Omega^d(\mathcal{X}) .$$

The dual filtration

$$\{0\} = \Omega^d(\mathcal{X})'_0 \subseteq \Omega^d(\mathcal{X})'_1 \subseteq \ldots \subseteq \Omega^d(\mathcal{X})'_{d+1} = \Omega^d(\mathcal{X})'$$

is given by
$$\Omega^d(\mathcal{X})'_j := [\Omega^d(\mathcal{X})/\Omega^d(\mathcal{X})^j]' \ .$$

The second statement in Corollary 3 immediately implies that the latter filtration corresponds under our map $I$ to the filtration of $C^{\mathrm{an}}(G, K)$ defined through annihilation conditions with respect to the left invariant differential operators in the ideal sequence $\mathfrak{b}_0 \supseteq \ldots \supseteq \mathfrak{b}_{d+1} = \mathfrak{b}$, i.e.,

$$I(\Omega^d(\mathcal{X})'_j) \subseteq C^{\mathrm{an}}(G, K)^{\mathfrak{b}_j=0} \text{ for } 0 \leq j \leq d+1 \ .$$

The compatibility between the local and the global filtration is established in the subsequent lemma.

**Lemma 4:**

$\Omega^d(\mathcal{X})^j \subseteq \Omega^d(\mathcal{X}) \cap \Omega^d_b(U^0)^j$.

Proof: Consider any $d$-form $u_*(\Xi_\mu \xi)$ with $u \in U(J(\mu))$ and write

$$u_*(\Xi_\mu \xi)|U^0 = \sum_{\nu \in X^*(\overline{T})} a(\nu) \Xi_\nu \xi \ .$$

We claim that $a(\nu) \neq 0$ implies that $\#J(\nu) \geq \#J(\mu)$. In order to see this let $\mu = \sum_k m_k \varepsilon_k$. Because of the condition on $u$ we have

$$u_*(\Xi_\mu \xi) = \Big( \prod_{k \notin J(\mu)} \Xi_k \Big) \Big( \prod_{k \in J(\mu)} \Xi_k^{m_k} \Big) \Big( \prod_{k \notin J(\mu)} u_* \Xi_k^{m_k - 1} \Big) \xi \ .$$

The first two products together contain each $\Xi_k$ with a positive exponent. In the third product the exponents are negative. On $U^0$ the summands of the linear form $u_* \Xi_k = \Xi_k + u_{k+1 k} \Xi_{k+1} + \ldots + u_{dk} \Xi_d$ differ pairwise in valuation. Hence after factoring out the largest summand we can develop $(u_* \Xi_k)^{-1}$, on $U^0$, into a geometric series. The terms of the resulting series have powers of a single $\Xi_{k'}$ in the denominator. It follows that each of the $d + 1 - \#J(\mu)$ factors in the third product can cancel out at most one of the $\Xi_k$'s in the first two products so that at least $d + 1 - (d + 1 - \#J(\mu)) = \#J(\mu)$ others remain. This establishes our claim which was that

$$u_*(\Xi_\mu \xi)|U^0 \in \Omega^d_b(U^0)^{\#J(\mu)} \text{ for } u \in U(J(\mu)) \ .$$

(For this slight reformulation one only has to observe that $\Omega^d_b(U^0)^j$ has a non-vanishing weight space exactly for those $\nu$ with $\#J(\nu) \geq j$.) It is then a consequence of Cor. 3 that
$$\Omega^d_{\mathrm{alg}}(\mathcal{X})^j \subseteq \Omega^d_b(U^0)^j \ .$$

As a simultaneous kernel of certain among the continuous linear forms $\eta \longmapsto \text{Res}_{(\overline{C},0)} \Xi_{-\mu} \eta$ on $\Omega_b^d(U^0)$ the right hand side is closed in $\Omega_b^d(U^0)$. It therefore follows that

$$\Omega^d(\mathcal{X})^j \subseteq \Omega_b^d(U^0)^j .$$
□

As a consequence of this fact we may view the map $D_J$ from the end of section 5 as a continuous linear map

$$D_J : \mathcal{O}(U_J^+ \cap B, M_J')^{\partial_J = 0} \longrightarrow [\Omega^d(\mathcal{X})^j / \Omega^d(\mathcal{X})^{j+1}]' ,$$

where $j := \#J$, which sends the weight vector $\sum_{\mu \in B(J)} [(L_\mu f_\nu)|U_J^+ \cap B] \otimes L_\mu^*$, for $\nu$ with $J(\nu) = J$, to the linear form $\lambda_\nu(\eta) := \text{Res}_{(\overline{C},0)} \Xi_{-\nu} \eta$.

We finish this section by collecting the basic properties which the subquotients of our global filtration have as locally convex vector spaces.

**Proposition 5:**

*Each subquotient $\Omega^d(\mathcal{X})^i / \Omega^d(\mathcal{X})^j$ for $0 \leq i \leq j \leq d+1$ is a reflexive Fréchet space; in particular its strong dual is barrelled and complete.*

Proof: In section 1 we deduced the reflexivity of $\Omega^d(\mathcal{X})$ from the fact that it is the projective limit of a sequence of Banach spaces with compact transition maps. It is a general fact (the proofs of Theorems 2 and 3 in [Kom] carry over literally to the nonarchimedean situation) that in such a Fréchet space every closed subspace along with its corresponding quotient space are projective limits of this type, too.

**Lemma 6:**

*Let $A : V \longrightarrow \tilde{V}$ be a strict continuous linear map between the $K$-Fréchet spaces $V$ and $\tilde{V}$; if $\tilde{V}$ is reflexive then the dual map $A' : \tilde{V}' \longrightarrow V'$ between the strong duals is strict as well.*

Proof: (Recall that $A$ is strict if on $\text{im}(A)$ the quotient topology from $V$ coincides with the subspace topology from $\tilde{V}$.) The subspace $\text{im}(A)$ of $\tilde{V}$ being a quotient of the Fréchet space $V$ is complete by the open mapping theorem and hence is closed. Let now $\tilde{\Sigma} \subseteq \tilde{V}'$ be any open $o$-submodule. We have to find an open $o$-submodule $\Sigma \subseteq V'$ such that $A'(\tilde{\Sigma}) \supseteq \text{im}(A') \cap \Sigma$. We may assume that $\ker(A') \subseteq \tilde{\Sigma}$. By the definition of the strong dual we also may assume that $\tilde{\Sigma} = \tilde{\Gamma}^\circ := \{\lambda \in \tilde{V}' : |\lambda(\tilde{v})| \leq 1 \text{ for any } \tilde{v} \in \tilde{\Gamma}\}$ for some closed and bounded $o$-submodule $\tilde{\Gamma} \subseteq \tilde{V}$. Since $\tilde{V}$ is reflexive $\tilde{\Gamma}$ is weakly compact ([Tie] Thms 4.20.b, 4.21, and 4.25.2) and hence compact ([DeG] Prop. 3.b). Since $\text{im}(A)$ is closed in

!V the Hahn-Banach theorem ([Tie] Thm. 3.5) implies that $\ker(A')^o = \text{im}(A)$. Using [Tie] Thm. 4.14 we deduce form the inclusion $\ker(A') \subseteq \tilde{\Sigma}$ that

$$\tilde{\Gamma} = \tilde{\Gamma}^{oo} = \tilde{\Sigma}^o \subseteq \ker(A')^o = \text{im}(A) \ .$$

In fact, $\tilde{\Gamma}$ is a compact subset of $\text{im}(A)$. According to [B-GT] IX 2.10, Prop. 18 we find a compact subset $\Gamma \subseteq V$ such that $A(\Gamma) = \tilde{\Gamma}$. Then $\Sigma := \Gamma^o$ is an open $o$-submodule in $V'$ such that $\text{im}(A') \cap \Sigma = A'(\tilde{\Sigma})$.

**Proposition 7:**

i. For $0 \leq j \leq d+1$ the natural map $\Omega^d(\mathcal{X})'_j \hookrightarrow \Omega^d(\mathcal{X})'$ is a topological embedding as a closed subspace;

ii. for $0 \leq i \leq j \leq d+1$ the natural map $\Omega^d(\mathcal{X})'_j/\Omega^d(\mathcal{X})'_i \xrightarrow{\cong} [\Omega^d(\mathcal{X})^i/\Omega^d(\mathcal{X})^j]'$ is a topological isomorphism.

Proof: i. This follows immediately from Prop. 5 and Lemma 6. ii. The natural exact sequence

$$0 \longrightarrow \Omega^d(\mathcal{X})^i/\Omega^d(\mathcal{X})^j \longrightarrow \Omega^d(\mathcal{X})/\Omega^d(\mathcal{X})^j \longrightarrow \Omega^d(\mathcal{X})/\Omega^d(\mathcal{X})^i \longrightarrow 0$$

consists of strict linear maps between Fréchet spaces which are reflexive by Prop. 5. The dual sequence is exact by Hahn-Banach and consists of strict linear maps by Lemma 6.

**Corollary 8:**

*If $V$ denotes one of the locally convex vector spaces appearing in the previous Proposition then the $G$-action $G \times V \longrightarrow V$ is continuous and the map $g \longmapsto g\lambda$ on $G$, for any $\lambda \in V$, is locally analytic.*

Proof: Because of Prop. 7 this is a consequence of Cor. 3.9.

### 7. The top filtration step

The purpose of this section is to describe the first stage of the global filtration in various different ways. This information (for all the $p$-adic symmetric spaces of dimension $\leq d$) will be used in an essential way in our computation of all the stages of the global filtration in the last section.

**Theorem 1:**

*The following three subspaces of $\Omega^d(\mathcal{X})$ are the same:*

1. *The subspace $d(\Omega^{d-1}(\mathcal{X}))$ of exact forms in $\Omega^d(\mathcal{X})$;*

2. *The first stage $\Omega^d(\mathcal{X})^1$ in the global filtration;*

3. *The subspace of forms $\eta$ such that $\operatorname{Res}_{(\overline{C},0)} g_*\eta = 0$ for any $g \in G$.*

*In particular all three are closed subspaces.*

The proof requires a series of preparatory statements which partly are of interest in their own right. We recall right away that any exact form of course has vanishing residues. The subspace $\Omega^d(\mathcal{X})^1$ is closed by construction. The subspace in 3. is closed as the simultaneous kernel of a family of continuous linear forms.

**Lemma 2:**

*An algebraic differential form $\eta \in \Omega^d_{\operatorname{alg}}(\mathcal{X})$ is exact if and only if $\eta$ belongs to $\Omega^d_{\operatorname{alg}}(\mathcal{X})^1$.*

Proof: Suppose first that $\eta$ is exact. Expand $\eta$ in its partial fractions decomposition (Prop. 6.2). From Cor. 6.3 we see that $\eta$ is congruent to a finite sum of logarithmic forms $u_*\xi$ modulo $\Omega^d_{\operatorname{alg}}(\mathcal{X})^1$, where $u$ is in the subgroup $U$ of lower triangular unipotent matrices. However, by [ST] Thm. 24, Cor. 40, and Cor. 50 the forms $u_*\xi$ are linearly independent modulo exact forms. Since $\eta$ is exact, therefore, no logarithmic terms can appear in its partial fractions expansion and $\eta$ belongs to $\Omega^d_{\operatorname{alg}}(\mathcal{X})^1$. Conversely it suffices, by Cor. 6.3 and G-invariance, to consider a form $\Xi_\mu \xi$ with $\mu \neq 0$. Since the Weyl group acts through the sign character on $\xi$ we may use G-equivariance again and assume that $\varepsilon_0$ occurs in $\mu$ with a positive coefficient $m_0 > 0$. Then $\Xi_\mu \xi = d\theta$ with $\theta := \frac{1}{m_0} \Xi_{\beta_0} \Xi_{\mu-\beta} d\Xi_{\beta_1} \wedge \ldots \wedge d\Xi_{\beta_{d-1}}$. $\square$

In the following we let $\Omega^d(\mathcal{X}_n)^j$, for $n \in \mathbb{N}$, denote the closure of $\Omega^d_{\operatorname{alg}}(\mathcal{X}_n)^j$ in the Banach space $\Omega^d(\mathcal{X}_n)$.

**Lemma 3:**

*For a form $\eta \in \Omega^d(\mathcal{X})$ we have:*

i. *$\eta$ is exact if and only if $\eta|\mathcal{X}_n$ is exact for any $n \in \mathbb{N}$;*

ii. *$\eta \in \Omega^d(\mathcal{X})^1$ if and only if $\eta|\mathcal{X}_n \in \Omega^d(\mathcal{X}_n)^1$ for any $n \in \mathbb{N}$.*

Proof: i. By the formula on the bottom of p. 64 in [SS] we have

$$H^*_{DR}(\mathcal{X}) = \varprojlim_n H^*_{DR}(\mathcal{X}_n^{\mathrm{o}})$$

where the $\mathcal{X}_n^{\mathrm{o}} \subseteq \mathcal{X}$ are certain admissible open subvarieties such that

- $\mathcal{X} = \bigcup_n \mathcal{X}_n^{\mathrm{o}}$ is an admissible covering, and

- $\mathcal{X}_{n-1} \subseteq \mathcal{X}_n^\circ \subseteq \mathcal{X}_n$.

The second property of course implies that

$$\varprojlim_n H^*_{DR}(\mathcal{X}_n^\circ) = \varprojlim_n H^*_{DR}(\mathcal{X}_n) .$$

ii. This follows by a standard argument about closed subspaces of projective limits of Banach spaces (compare the proof of Thm. 2 in [Kom]). □

The main technique for the proof of Theorem 1 will be a "convergent partial fractions" decomposition for rigid $d$-forms on $\mathcal{X}$. We begin by recalling the explicit description of rigid forms on $\mathcal{X}_n$ given in [SS] p. 53. Fix a set $\mathcal{H} = \{\ell_0, \ldots \ell_s\}$ of unimodular representatives for the hyperplanes modulo $\pi^{n+1}$ in such a way that it contains the coordinate hyperplanes $\{\Xi_i = 0\}$ for $0 \leq i \leq d$. A rigid $d$-form $\eta$ on the affinoid $\mathcal{X}_n$ is represented by a convergent expansion

$$(*) \qquad \eta = \sum_{I,J} a_{I,J} \frac{\Xi_0^{j_0} \cdot \ldots \cdot \Xi_d^{j_d}}{\ell_0^{i_0} \cdot \ldots \cdot \ell_s^{i_s}} \Theta$$

in homogeneous coordinates where $I$ and $J$ run over all $(s+1)$-tuples $(i_0, \ldots, i_s)$ and $(d+1)$-tuples $(j_0, \ldots, j_d)$ of non-negative integers respectively with $\sum i_k - \sum j_k = d+1$ and where

$$\Theta := \sum_{i=0}^{d} (-1)^i \Xi_i d\Xi_0 \wedge \ldots \wedge \widehat{d\Xi_i} \wedge \ldots \wedge d\Xi_d .$$

The convergence means that the coefficients $a_{I,J}$ satisfy $\omega(a_{I,J}) - n(\sum_{k=0}^{d} j_k) \longrightarrow \infty$ as $\sum_{k=0}^{d} j_k \longrightarrow \infty$.

**Lemma 4:**

*In the expansion $(*)$ we may assume $a_{I,J} = 0$ unless the corresponding set of "denominator forms" $\{\ell_k : i_k \geq 1\}$ is linearly independent.*

Proof: Suppose that $\ell_0, \ldots, \ell_r$ are linearly dependent, and that $i_k \geq 1$ for $0 \leq k \leq r$. Write

$$\sum_{k=0}^{r} b_k \ell_k = 0$$

with the $b_k \in o$ and at least one $b_k = 1$. Suppose for example that $b_0 = 1$. Then

$$\ell_0 = -\sum_{k=1}^{r} b_k \ell_k$$

and
$$\frac{\Xi_0^{j_0} \cdot \ldots \cdot \Xi_d^{j_d}}{\ell_0^{i_0} \cdot \ldots \cdot \ell_r^{i_r}} \Theta = -\sum_{k=1}^{r} \frac{b_k \Xi_0^{j_0} \cdot \ldots \cdot \Xi_d^{j_d}}{\ell_0^{i_0+1} \ldots \ell_k^{i_k-1} \ldots \ell_r^{i_r}} \Theta \ .$$

The individual terms on the right side of this sum have the same degree as the term on the left. This, together with the fact that the $b_k$ belong to $o$, implies that the expression on the right may be substituted into the series expansion for $\eta$ and the sum re-arranged. Further, this process may be iterated until the denominators occurring on the right side are linearly independent. $\square$

Using this Lemma, we see that any $\eta \in \Omega^d(\mathcal{X}_n)$ can be written as a finite sum of forms

$$(**) \qquad \eta_L = \sum_{I,J} a_{I,J} \frac{\Xi_0^{j_0} \cdot \ldots \cdot \Xi_d^{j_d}}{\ell_0^{i_0} \cdot \ldots \cdot \ell_r^{i_r}} \Theta$$

where $L = \{\ell_0, \ldots, \ell_r\}$ is a fixed linearly independent set chosen from $\mathcal{H}$, $I$ runs through the $(r+1)$-tuples of positive integers, and $\omega(a_{I,J}) - n(\sum_{k=0}^{d} j_k) \longrightarrow \infty$ as $\sum_{k=0}^{d} j_k \longrightarrow \infty$.

**Lemma 5:**

*A form $\eta_L$ as in $(**)$ belongs to $\Omega^d(\mathcal{X}_n)^{d+1-\#L}$.*

Proof: This is clear from the inequality $\imath(\eta) \leq \operatorname{codim} Z(\eta)$ in section 6.

**Definition:**

*A form $\eta \in \Omega^d(\mathcal{X}_n)$ is called decomposable if it has a convergent expansion of the form*

$$\eta = \sum_{g \in G} \sum_{\mu \in X^*(\overline{T})} c(g,\mu)(g_*(\Xi_\mu \xi)|\mathcal{X}_n)$$

*where*

1. *$c(g,\mu) \in K$ and $= 0$ for all but finitely many $g \in G$ which are independent of $\mu$,*

2. *if $c(g,\mu) \neq 0$ for some $\mu$ then the columns of the matrix $g$ are unimodular,*

3. *$\omega(c(g,\mu)) - nd(\mu) \longrightarrow \infty$ as $d(\mu) \longrightarrow \infty$ ($d(\mu)$ was defined in section 4, just before the statement of Lemma 4.1).*

**Lemma 6:**

*Suppose that $\eta_L$ is given by a series as in $(**)$ on $\mathcal{X}_{2n}$. Then the restriction of $\eta_L$ to $\mathcal{X}_{n-1}$ is either decomposable or may be written as a (finite) sum of series $\eta_{L'}$ converging on $\mathcal{X}_{n-1}$ and with $\#L' < \#L$.*

Proof: The dichotomy in the statement of the Lemma arises out of the following two possibilities:

*Case I.* There is a unimodular relation

$$\sum_{k=0}^{r} b_k \ell_k \equiv 0 (\mathrm{mod}\, \pi^n)\ .$$

*Case II.* Whenever there is a relation

$$\sum_{k=0}^{r} b_k \ell_k \equiv 0 (\mathrm{mod}\, \pi^n)\ , \text{ with } b_k \in o\ ,$$

we must have all $b_k$ divisible by $\pi$.

Let us treat Case I first. Suppose that $b_0$ is a unit in the unimodular relation, and write

$$\ell_0 \ =\ -\sum_{k=1}^{r} (b_k/b_0)\ell_k + \pi^n h\ .$$

To simplify the notation, set

$$\ell := -\sum_{k=1}^{r} (b_k/b_0)\ell_k\ .$$

The fact that $\ell_0$ is unimodular means that $\ell$ is unimodular as well. We have

$$\frac{1}{\ell_0}\ =\ \frac{1}{\ell}\,(1 + \pi^n h/\ell)^{-1}\ ,$$

and, since $\pi^n h/\ell$ has sup-norm $\leq |\pi| < 1$ on $\mathcal{X}_{n-1}$, using the geometric series we see that we may rewrite the series expansion for $\eta_L$ so that it converges on $\mathcal{X}_{n-1}$:

$$\eta_L|\mathcal{X}_{n-1}\ =\ \sum a'_{I,J}\,\frac{\Xi_0^{j_0}\cdot\ldots\cdot\Xi_d^{j_d}}{\ell^i \ell_1^{i_1}\cdot\ldots\cdot\ell_r^{i_r}}\,\Theta$$

But since $\ell$ is a linear combination of $\ell_k$ for $k \neq 0$, the proof of Lemma 4 shows that $\eta_L|\mathcal{X}_{n-1}$ is a sum of series $\eta_{L'}$ where $L'$ is a linearly independent subset of the dependent set $\{\ell, \ell_1, \ldots, \ell_r\}$; such a set has fewer than $r+1$ elements.

For Case II we take a different approach. Apply elementary divisors to find linear forms $f_0, \ldots, f_d$ which form a basis for the $o$-lattice spanned by $\Xi_0, \ldots, \Xi_d$ and such that $\pi^{e_0} f_0, \ldots, \pi^{e_r} f_r$ form a basis for the span of $\ell_0, \ldots, \ell_r$. Since any monomial in the $\Xi_i$ is an integral linear combination of monomials in the $f_i$, we may rewrite $\eta_L$ using the $f_i$ for coordinates:

$$\eta_L = \sum a''_{I,J} \frac{f_0^{j_0} \cdot \ldots \cdot f_d^{j_d}}{\ell_0^{i_0} \cdot \ldots \cdot \ell_r^{i_r}} \Theta .$$

Using our Case II hypothesis, we know that $e_k < n$ for $0 \leq k \leq r$. Therefore $\pi^{n-1} f_k$, for each $0 \leq k \leq r$, is an integral linear combination of $\ell_0, \ldots, \ell_r$. Let $g \in G$ be the matrix such that $g_* \Xi_i = \ell_i$ for $0 \leq i \leq r$ and $g_* \Xi_i = f_i$ for $r+1 \leq i \leq d$. By construction the columns of $g$ are unimodular. Rewriting the series for $\eta_L$ in terms of the $\pi^{1-n} \ell_0, \ldots, \pi^{1-n} \ell_r, f_{r+1}, \ldots, f_d$ we see that

$$(***) \qquad \eta_L = \sum a''_{I,J} \det(g)^{-1} \pi^{(1-n)(\sum_{k=0}^{r} j_k)} \sum_\mu c_{\mu,I,J} g_*(\Xi_\mu \xi)$$

where each of the inner sums is finite and the coefficients $c_{\mu,I,J}$ are integral. Since the original sum for $\eta_L$ converges on $\mathcal{X}_{2n}$, we have

$$\omega(a''_{I,J}) = H(\sum_{k=0}^{d} j_k) + 2n(\sum_{k=0}^{d} j_k)$$

where $H(m)$ is a function which goes to infinity as $m$ goes to infinity. But then

$$\omega(a''_{I,J} \det(g)^{-1} \pi^{(1-n)(\sum_{k=0}^{r} j_k)}) \geq H(\sum_{k=0}^{d} j_k) + 2n(\sum_{k=0}^{d} j_k) + (1-n)(\sum_{k=0}^{d} j_k) + C$$

which shows that, after rearrangement according to $\mu$, the series $(***)$ converges on $\mathcal{X}_{n+1}$ (if $c_{\mu,I,J} \neq 0$ then $\sum_{k=0}^{d} j_k \geq d(\mu)$). Thus in Case II $\eta_L$ is decomposable on $\mathcal{X}_{n+1}$.

**Proposition 7:**

*Let $\eta$ be a rigid $d$-form on $\mathcal{X}$; then $\eta | \mathcal{X}_n$, for any $n > 0$, is decomposable.*

Proof: This follows by induction from Lemma 6. Indeed, any rigid form $\eta$ on $\mathcal{X}_m$ with $m := 2^{d+1}(n+2)$ is decomposable on $\mathcal{X}_n$.

**Definition:**

A form $\eta \in \Omega^d_{\text{alg}}(\mathcal{X})$ is called logarithmic if it lies in the smallest $G$-invariant vector subspace containing $\xi$.

**Corollary 8:**

Let $\eta$ be a rigid $d$-form on $\mathcal{X}$. Then, for any $n > 0$, the restriction of $\eta$ to $\mathcal{X}_n$ has a decomposition

$$\eta = \eta_0 + \eta_1$$

where $\eta_0$ is the restriction of a logarithmic form and $\eta_1$ is an exact form in $\Omega^d(\mathcal{X}_n)^1$.

Proof: Applying the convergent partial fractions decomposition of Prop. 7, write $\eta$ on $\mathcal{X}_{n+1}$ as

$$\eta|\mathcal{X}_{n+1} = \sum_g c(g,0)(g_*\xi|\mathcal{X}_{n+1}) + \sum_g \sum_{\mu \neq 0} c(g,\mu)(g_*(\Xi_\mu \xi)|\mathcal{X}_{n+1}) \ .$$

Let $\eta_0$ be the first of these sums, and $\eta_1$ the second. Clearly $\eta_0$ is logarithmic and $\eta_1$, by Cor. 6.3, belongs to $\Omega^d(\mathcal{X}_{n+1})^1$. Thus we need only show that $\eta_1$ is exact on $\mathcal{X}_n$. However, one sees easily that the series for $\eta_1$ may be integrated term-by-term to obtain a rigid $(d-1)$-form $\theta$ on $\mathcal{X}_n$ with $d\theta = \eta_1$ (compare the proof of Lemma 2). □

As a last preparation we need the following result on logarithmic forms.

**Proposition 9:**

For any $n > 0$ we have:

i. There is a compact open set $V_n \subset U$ such that

$$u_*\xi \in \Omega^d(\mathcal{X}_n)^1 \cap d(\Omega^{d-1}(\mathcal{X}_n))$$

for all $u \in U \backslash V_n$;

ii. there is a finite set $u^{(1)}, \ldots, u^{(k)}$ of elements of $U$ and a disjoint covering of $V_n$ by sets $\{D(u^{(\ell)}, r)\}_{\ell=1}^k$ such that

$$v_*\xi \equiv u_*^{(\ell)}\xi \pmod{\Omega^d(\mathcal{X}_n)^1 \cap d(\Omega^{d-1}(\mathcal{X}_n))}$$

if $v \in D(u^{(\ell)}, r)$;

iii. the image of

$$\Omega^d(\mathcal{X}) \to \Omega^d(\mathcal{X}_n)/(\Omega^d(\mathcal{X}_n)^1 \cap d(\Omega^{d-1}(\mathcal{X}_n)))$$

and the space $\Omega^d(\mathcal{X}_n)/\Omega^d(\mathcal{X}_n)^1$ both are finite dimensional; more precisely, the classes of the forms $u_*^{(1)}\xi, \ldots, u_*^{(k)}\xi$ span both spaces.

Proof: i. In homogeneous coordinates, we write

$$u_*\xi = \frac{\Theta}{\ell_0 \cdots \ell_d}$$

where $\ell_j = \sum_{i=j}^d u_{ij}\Xi_i$ and the $u_{ij}$ are the matrix entries of the lower triangular unipotent matrix $u$. Let

$$V_n := \{u \in U : \omega(u_{lk}) \geq -(n+1)d \text{ for all } d \geq l \geq k \geq 0\}.$$

We claim $V_n$ has the desired property. Suppose that $u \notin V_n$, so that, for some pair $d \geq l > k \geq 0$ we have $\omega(u_{lk}) < -(n+1)d$. Focus attention for the moment on the linear form $\ell_k$. Since $u_{kk} = 1$, we may split the set of row indices $k, \ldots, d$ into two nonempty sets $A$ and $B$ such that

$$\inf_{l \in A} \omega(u_{lk}) > \sup_{l \in B} \omega(u_{lk}) + n + 1.$$

We point out two facts for later use. First, the index $k$ automatically belongs to the set $A$, and so $\ell_k^B$ is a linear combination of the $\Xi_i$ with $i > k$. Second, and for the same reason, the set of linear forms $\{\ell_j\}_{j \neq k} \cup \{\ell_k^A\}$ is a triangular basis for the full space of $K$-linear forms in the $\Xi_i$. Continuing with the main line of argument, write

$$\ell_k = \ell_k^A + \ell_k^B = (\sum_{l \in A} u_{lk}\Xi_l) + (\sum_{l \in B} u_{lk}\Xi_l).$$

Then

$$\frac{1}{\ell_k} = \frac{1}{\ell_k^B}\left(\frac{1}{1 + (\ell_k^A/\ell_k^B)}\right).$$

The linear forms $\pi^{-\inf_{l \in A} \omega(u_{lk})}\ell_k^A$ and $\pi^{-\inf_{l \in B} \omega(u_{lk})}\ell_k^B$ are unimodular. From this, we obtain the following estimate on $\mathcal{X}_n$:

(1)
$$\begin{aligned}
\omega(\ell_k^A/\ell_k^B) &\geq \inf_{l \in A} \omega(u_{lk}) - \inf_{l \in B} \omega(u_{lk}) - n \\
&\geq \inf_{l \in A} \omega(u_{lk}) - \sup_{l \in B} \omega(u_{lk}) - n \\
&> 1 \, .
\end{aligned}$$

At this point, it will be convenient to change from homogeneous to inhomogeneous coordinates. Let

$$\overline{\ell}_j := \ell_j/\Xi_d = \sum_{i=0}^{d-1} u_{ij}\Xi_{\beta_i} + u_{dj},$$

and similarly let $\bar{\ell}_k^A := \ell_k^A/\Xi_d$ and $\bar{\ell}_k^B := \ell_k^B/\Xi_d$. Then we may expand the form $u_*\xi$ as a convergent series on $\mathcal{X}_n$:

$$u_*\xi = \sum_{m=0}^{\infty} c_m F_m d\Xi_{\beta_0} \wedge \cdots \wedge d\Xi_{\beta_{d-1}} \tag{2}$$

where the coefficients $c_m \in \mathbb{Z}$,

$$F_m = \frac{(\bar{\ell}_k^A)^m}{\bar{\ell}_0 \cdots (\bar{\ell}_k^B)^{m+1} \cdots \bar{\ell}_{d-1}}$$

and $\bar{\ell}_k^B$ has taken the place of $\bar{\ell}_k$ in the denominators of these forms (observe that $\Theta = (-1)^d \Xi_d^{d+1} d\Xi_{\beta_0} \wedge \cdots \wedge d\Xi_{\beta_{d-1}}$). Our estimate (1) tells us that there is a constant $C$ so that the functions $F_m$ satisfy $\inf_{q \in \mathcal{X}_n} \omega(F_m(q)) \geq m - C$ in the sup norm on $\mathcal{X}_n$.

To finish the proof, we will show that the expansion (2) may be integrated term by term on $\mathcal{X}_n$. This shows that $u_*\xi$ is exact on $\mathcal{X}_n$. In addition, since it proves that each algebraic form in the expansion (2) is exact, we see from Lemma 2 that these forms belong to $\Omega_{\text{alg}}^d(\mathcal{X}_n)^1$ and so $u_*\xi$ belongs to $\Omega^d(\mathcal{X}_n)^1$ as well.

As we remarked earlier, the forms $\ell_0, \ldots, \ell_k^A, \ldots, \ell_d$ are a triangular basis for the space of all linear forms. Therefore, we may choose $v \in U$ so that $v_*\Xi_j = \ell_j$ for all $0 \leq j \leq d$ except for $j = k$, and $v_*\Xi_k = \ell_k^A$. Let $f := v_*^{-1}(\ell_k^B/\Xi_d)$. The form $f$ does not involve $\Xi_{\beta_k}$. Then we compute

$$F_m = v_*\left(\frac{\Xi_{\beta_k}^m}{\Xi_{\beta_0} \cdots \Xi_{\beta_{k-1}} f^{m+1} \Xi_{\beta_{k+1}} \cdots \Xi_{\beta_{d-1}}}\right).$$

Using this and the estimate for the $F_m$, we see that

$$\theta = \left(\sum_{m=0}^{\infty} \frac{c_m}{m+1} F_m\right) v_*((-1)^k \Xi_{\beta_k} d\Xi_{\beta_0} \wedge \cdots \wedge \widehat{d\Xi_{\beta_k}} \wedge \cdots \wedge d\Xi_{\beta_{d-1}})$$

is a convergent expansion for a rigid $(d-1)$-form $\theta$ on $\mathcal{X}_n$, and that $d\theta = u_*\xi$.

ii. In the notation of Prop. 3.1, let $u^{(1)}, \ldots, u^{(k)}$ be finitely many elements of $U$ so that the open sets $\{D(u^{(\ell)}, r)\}_{\ell=1}^k$ form a disjoint covering of $V_n$ and so that, for each $\ell = 1, \ldots, k$,

$$\omega(v_{ji} - u_{ji}^{(\ell)}) > 2(n+1) \text{ for all } v \in D(u^{(\ell)}, r) \text{ and all } 0 \leq i < j \leq d. \tag{3}$$

Then, for $v \in D(u^{(\ell)}, r)$, we have the uniformly convergent expansion (*) on $\mathcal{X}_n$ from the proof of Prop. 3.1, where, to simplify the notation, we write $u = u^{(\ell)}$:

$$v_*\xi = k(vw_{d+1}, \cdot) d\Xi_{\beta_0} \wedge \cdots \wedge d\Xi_{\beta_{d-1}}$$

with
$$k(vw_{d+1}, q) = \sum_{\underline{m}} c_{\underline{m}} h_{\underline{m}} \cdot (v-u)^{\underline{m}},$$

$$h_{\underline{m}} = \frac{\Xi_{\mu(\underline{m})}(q)}{f_0(u,q)^{s_0(\underline{m})} \cdot \ldots \cdot f_{d-1}(u,q)^{s_{d-1}(\underline{m})}},$$

and $c_{\underline{m}} \in \mathbb{Z}$. In this expansion, the term with $\underline{m} = (0, \ldots, 0)$ is $u_*\xi = u_*^{(\ell)}\xi$. Also, comparing the estimate in (3) with those used in Prop. 3.1, we see that we have
$$\inf_{q \in \mathcal{X}_n} \omega(h_{\underline{m}}(q) \cdot (v-u)^{\underline{m}}) \geq \big(\sum_{0 \leq i < j \leq d} m_{ji}\big) - nd.$$

We claim that, except for the term with $\underline{m} = (0, \ldots, 0)$, this series may be integrated term by term to yield a convergent $(d-1)$-form on $\mathcal{X}_n$. This means that $(v_*\xi - u_*\xi)|\mathcal{X}_n$ is an exact form, and further that (just as in the proof of the first assertion) each term in the expansion of $v_*\xi - u_*\xi$ is an exact algebraic form, so that $v_*\xi - u_*\xi$ belongs to $\Omega^d(\mathcal{X}_n)^1$. In other words,
$$v_*\xi \equiv u_*^{(\ell)}\xi \pmod{\Omega^d(\mathcal{X}_n)^1 \cap d(\Omega^{d-1}(\mathcal{X}_n))}.$$

To prove our claim, let $S_j$ be the set of $\underline{m}$ such that $s_i(\underline{m}) = 1$ for $i = 0, \ldots, j-1$ but $s_j(\underline{m}) > 1$. Let
$$F_j := \sum_{\underline{m} \in S_j} h_{\underline{m}} \cdot (v-u)^{\underline{m}}$$

and
$$\eta_j := F_j d\Xi_{\beta_0} \wedge \cdots \wedge d\Xi_{\beta_{d-1}}.$$

Because
$$v_*\xi = \eta_0 + \ldots + \eta_{d-1} + u_*^{(\ell)}\xi,$$

it suffices to integrate each $\eta_j$ term by term. Notice that if $\underline{m} \in S_j$, then $\Xi_{\mu(\underline{m})}$ does not involve any of $\Xi_{\beta_i}$ for $i = 0, \ldots, j$. We may choose a matrix $g \in U$ so that $g_*\Xi_{\beta_i} = f_i(u, \cdot)$ for $i = 0, \ldots, j$ and $g_*\Xi_{\beta_i} = \Xi_{\beta_i}$ for $i = j+1, \ldots, d-1$. Now set
$$G_j := \sum_{\underline{m} \in S_j} \frac{c_{\underline{m}}}{1 - s_j(\underline{m})} h_{\underline{m}}(v-u)^{\underline{m}}.$$

The estimate on the sup norm for $h_{\underline{m}}$ implies that this is the convergent expansion of a rigid function on $\mathcal{X}_n$. Therefore
$$\theta_j := (-1)^j G_j g_*(\Xi_{\beta_j} d\Xi_{\beta_0} \wedge \cdots \wedge \widehat{d\Xi_{\beta_j}} \wedge \cdots \wedge d\Xi_{\beta_{d-1}})$$

is a rigid $(d-1)$-form on $\mathcal{X}_n$. Furthermore, a simple computation shows that $d\theta_j = \eta_j$. Indeed, a typical term in the series for $\theta_j$ is

$$(4) \qquad \frac{c_{\underline{m}}}{1-s_j(\underline{m})} h_{\underline{m}}(v-u)^{\underline{m}} g_*((-1)^j \Xi_{\beta_j} d\Xi_{\beta_0} \wedge \cdots \wedge \widehat{d\Xi_{\beta_j}} \wedge \cdots \wedge d\Xi_{\beta_{d-1}}) \ .$$

Let

$$H_{\underline{m}} := \left( \frac{\Xi_{\mu(\underline{m})}}{\Xi_{\beta_0} \cdots \Xi_{\beta_{j-1}} \Xi_{\beta_j}^{s_j(\underline{m})} f_{j+1}(u,\cdot)^{s_{j+1}(\underline{m})} \cdots f_{d-1}(u,\cdot)^{s_{d-1}(\underline{m})}} \right),$$

so that $h_{\underline{m}} = g_* H_{\underline{m}}$. Then the term in (4) is

$$(-1)^j \frac{c_{\underline{m}}}{1-s_j(\underline{m})} (v-u)^{\underline{m}} g_*(H_{\underline{m}} \Xi_{\beta_j} d\Xi_{\beta_0} \wedge \cdots \wedge \widehat{d\Xi_{\beta_j}} \wedge \cdots \wedge d\Xi_{\beta_{d-1}}) \ .$$

We leave it as an exercise to verify that applying $d$ to this expression one obtains the term

$$c_{\underline{m}} h_{\underline{m}} (v-u)^{\underline{m}} d\Xi_{\beta_0} \wedge \cdots \wedge d\Xi_{\beta_{d-1}} \ .$$

iii. By Cor. 8 and [ST] Cor. 40, a form $\eta \in \Omega^d(\mathcal{X})$, restricted to $\mathcal{X}_n$, may be written

$$(\eta|\mathcal{X}_n) = \eta_0 + \eta_1$$

where

$$\eta_1 \in \Omega^d(\mathcal{X}_n)^1 \cap d(\Omega^{d-1}(\mathcal{X}_n))$$

and $\eta_0$ is (the restriction of) a finite sum of logarithmic forms $u_*\xi$. Thus the image of

$$\Omega^d(\mathcal{X}) \to \Omega^d(\mathcal{X}_n)/(\Omega^d(\mathcal{X}_n)^1 \cap d(\Omega^{d-1}(\mathcal{X}_n)))$$

is spanned by logarithmic forms $u_*\xi$. Similarly, from Prop. 3.3 we know that the logarithmic forms $u_*\xi$ generate $\Omega^d(\mathcal{X})$ as a topological vector space. Since the image of $\Omega^d(\mathcal{X})$ in $\Omega^d(\mathcal{X}_n)$ under restriction is dense the same forms $u_*\xi$ generate the quotient $\Omega^d(\mathcal{X}_n)/\Omega^d(\mathcal{X}_n)^1$ as a Banach space. In both cases we hence may conclude that, using the first assertion, the $u_*\xi$ for $u \in V_n$ and then, using the second assertion, even the $u_*^{(1)}\xi, \ldots, u_*^{(k)}\xi$ span the two vector spaces in question. $\square$

Proof of Theorem 1:

We show that each $\eta$ in the third space also lies in the intersection of the first two spaces. By Lemma 3, it suffices to show that the restriction of $\eta$ to $\mathcal{X}_n$ belongs to $E_n := \Omega^d(\mathcal{X}_n)^1 \cap d(\Omega^{d-1}(\mathcal{X}_n)))$ for all $n > 0$. We fix an $n$ and choose a finite set $u^{(1)}, \ldots, u^{(k)}$ of elements of $V_n$ as in Prop. 9. We also choose $m \geq n$ so that the image of $\mathcal{X}_m$ in the building contains the chambers $u^{(\ell)}(\overline{C}, 0)$.

Apply Cor. 8 to write $\eta|\mathcal{X}_m = \eta_0 + \eta_1$ on $\mathcal{X}_m$, with $\eta_0$ logarithmic and $\eta_1 \in E_m \subset E_n$. Our hypothesis on $m$ implies that the linear form $\operatorname{Res}_{u^{(\ell)}(\overline{C},0)}$ is continuous on $\Omega^d(\mathcal{X}_m)$, and since $\eta_1$ is exact on $\mathcal{X}_m$ we must have

$$\operatorname{Res}_{u^{(\ell)}(\overline{C},0)}(\eta_1) = 0 \text{ for } \ell = 0, \ldots, k.$$

Since all residues of $\eta$ are zero, we conclude that

$$\operatorname{Res}_{u^{(\ell)}(\overline{C},0)}(\eta_0) = 0 \text{ for } \ell = 0, \ldots, k.$$

We now need to show that, under our residue hypothesis, the restriction to $\mathcal{X}_n$ of the logarithmic form $\eta_0$ belongs to $E_n$. Since $\eta_0$ is a logarithmic form, we may write it as a sum of forms $u_*\xi$ with $u \in U$ ([ST] Cor. 40), and for our purposes we may (by Prop. 9.i) assume that all $u \in V_n$. Thus, for each $\ell$, we have finitely many distinct $v_{\ell j} \in V_n$ and constants $c_{\ell j}$ so that

$$\eta_0 = \sum_{\ell=1}^{k} \sum_{j=0}^{s_\ell} c_{\ell j}((v_{\ell j})_*\xi)$$

where, for $j = 0, \ldots, s_\ell$, we have $v_{\ell j} \in D(u^{(\ell)}, r)$. By the proofs of Facts A and B of [ST], page 430-431, we see that

$$\operatorname{Res}_{u^{(\ell)}(\overline{C},0)}(\eta_0) = \sum_{j=0}^{s_\ell} \operatorname{Res}_{u^{(\ell)}(\overline{C},0)} c_{\ell j}((v_{\ell j})_*\xi)$$
$$= \sum_{j=0}^{s_\ell} c_{\ell j}$$
$$= 0.$$

It then follows from Prop. 9.ii that

$$\eta_0|\mathcal{X}_n \equiv \sum_{\ell=1}^{k} \sum_{j=0}^{s_\ell} c_{\ell j} u_*^{(\ell)}\xi \pmod{E_n}$$
$$\equiv 0 \pmod{E_n}$$

as claimed. □

From section 3, in particular Lemma 3.5, we have the injective $G$-equivariant map
$$\begin{array}{rcl} I_o : \Omega^d(\mathcal{X})' & \hookrightarrow & C_o(U, K) \cong C(G/P, K)/C_{\mathrm{inv}}(G/P, K) \\ \lambda & \longmapsto & [u \mapsto \lambda(u_*\xi)]. \end{array}$$

Let $C^\infty(G/P, K) \subseteq C(G/P, K)$ denote the subspace of all locally constant functions and put $C^\infty_{\mathrm{inv}}(G/P, K) := C^\infty(G/P, K) \cap C_{\mathrm{inv}}(G/P, K)$ and $C^\infty_{\mathrm{o}}(U, K) := C^\infty(U, K) \cap C_{\mathrm{o}}(U, K)$. The quotient

$$\mathrm{St} := C^\infty(G/P, K)/C^\infty_{\mathrm{inv}}(G/P, K)$$

is an irreducible smooth $G$-representation known as the Steinberg representation of the group $G$. The above isomorphism for the target of $I_{\mathrm{o}}$ restricts to an isomorphism

$$C^\infty_{\mathrm{o}}(U, K) \cong \mathrm{St} .$$

**Proposition 10:**

If $\lambda \in \Omega^d(\mathcal{X})'$ vanishes on exact forms then the function $I_{\mathrm{o}}(\lambda)$ on $U$ is locally constant with compact support.

Proof: Such a linear form $\lambda$ extends continuously to $\Omega^d(\mathcal{X}_n)$ for some $n$. Since, by Thm. 1, it vanishes on $\Omega^d(\mathcal{X})^1$, it vanishes on $\Omega^d(\mathcal{X}_n)^1$. Then from Prop. 9.i it vanishes on $u_*\xi$ outside of $V_n$, and therefore the function in question is compactly supported. Prop. 9.ii shows that there is a finite disjoint covering of $V_n$ by sets $D(u^{(\ell)}, r)$ such that $\lambda(v_*\xi) = \lambda(u^{(\ell)}_*\xi)$ for $v \in D(u^{(\ell)}, r)$. Therefore the function in question is locally constant. $\square$

It follows that $I_{\mathrm{o}}$ induces an injective $G$-equivariant map

$$[\Omega^d(\mathcal{X}) / \genfrac{}{}{0pt}{}{\text{exact}}{\text{forms}}]' \hookrightarrow C^\infty_{\mathrm{o}}(U, K) \cong \mathrm{St} .$$

Since $\mathrm{Res}_{(\overline{C},0)}\xi$ is nonzero the left hand side contains a nonzero vector. But the right hand side is algebraically irreducible as a $G$-representation. Hence we see that this map must be bijective.
¿From our Theorem 1 and from the nonarchimedean version of [Kom] Thm. 3 we have the identifications of locally convex vector spaces

$$\Omega^d(\mathcal{X}) / \genfrac{}{}{0pt}{}{\text{exact}}{\text{forms}} = \Omega^d(\mathcal{X})/\Omega^d(\mathcal{X})^1 = \varprojlim_n \Omega^d(\mathcal{X}_n)/\Omega^d(\mathcal{X}_n)^1 .$$

On the other hand, Prop. 9 says that, for any $n > 0$, the space $\Omega^d(\mathcal{X}_n)/\Omega^d(\mathcal{X}_n)^1$ is finite dimensional. We conclude that $\Omega^d(\mathcal{X}) / \genfrac{}{}{0pt}{}{\text{exact}}{\text{forms}}$, resp. its dual space, is a projective, resp. injective, limit of finite dimensional Hausdorff spaces. In particular the topology on $[\Omega^d(\mathcal{X}) / \genfrac{}{}{0pt}{}{\text{exact}}{\text{forms}}]'$ is the finest locally convex topology. In this way we have computed the top step of our filtration as a topological vector space.

**Theorem 11:**

*The $G$-equivariant map*

$$[\Omega^d(\mathcal{X}) / \substack{\text{exact} \\ \text{forms}}]' \xrightarrow{\cong} \text{St}$$

$$\lambda \longmapsto [u \mapsto \lambda(u_*\xi)]$$

*is an isomorphism; morover, the topology of the strong dual on the left hand side is the finest locally convex topology.*

### 8. The partial boundary value maps

In this section we will introduce and study, for any $0 \leq j \leq d$, a "partial boundary value map" $I^{[j]}$ from $[\Omega^d(\mathcal{X})^j/\Omega^d(\mathcal{X})^{j+1}]'$ into a space of functions on $G$. Recall that we denoted by $\mathfrak{p}_J$, for any subset $J \subseteq \{0, \ldots, d\}$, the parabolic subalgebra in $\mathfrak{g}$ of all matrices which have a zero entry in position $(i,j)$ for $i \in J$ and $j \notin J$; moreover $\mathfrak{n}_J \subseteq \mathfrak{p}_J$ denoted the unipotent radical. Let $P_J \subseteq G$ be the parabolic subgroup whose Lie algebra is $\mathfrak{p}_J$ and let $U_J \subseteq P_J$ be its unipotent radical. We have the Levi decomposition $P_J = U_J L_J$ with $L_J := L'(J) \times L(J)$ and

$$L'(J) := \text{ all matrices in } G \text{ with}$$
$$\phantom{L'(J) :=} - \text{ a zero entry in position } (i,j)$$
$$\phantom{L'(J) := -} \text{ for } i \neq j \text{ and not both in } J, \text{ and}$$
$$\phantom{L'(J) :=} - \text{ an entry 1 in position } (i,i) \text{ for } i \notin J$$

and

$$L(J) := \text{ all matrices in } G \text{ with}$$
$$\phantom{L(J) :=} - \text{ a zero entry in position } (i,j)$$
$$\phantom{L(J) := -} \text{ for } i \neq j \text{ and } i \text{ or } j \in J, \text{ and}$$
$$\phantom{L(J) :=} - \text{ an entry 1 in position } (i,i) \text{ for } i \in J \ .$$

Clearly, $L'(J) \cong GL_{\#J}(K)$ and $L(J) \cong GL_{d+1-\#J}(K)$. With these new notations, the subgroup $U(J)$ from section 6 is the subgroup $U(J) = U \cap L(J)$ of lower triangular unipotent matrices in $L(J)$, and $\mathfrak{l}_J$, $\mathfrak{l}'(J)$, and $\mathfrak{l}(J)$ are the Lie algebras of $L_J$, $L'(J)$, and $L(J)$ respectively. In the following we are mostly interested in the subsets $\underline{j} := \{0, \ldots, j-1\}$ for $0 \leq j \leq d$. Let

$$V_j := \text{ closed subspace of } \Omega^d(\mathcal{X})^j/\Omega^d(\mathcal{X})^{j+1} \text{ spanned by}$$
$$\phantom{V_j :=} \text{ the forms } g_*(\Xi_\mu \xi) \text{ for } \mu \in B(\underline{j}) \text{ and } g \in L(\underline{j})$$

viewed as a locally convex vector space with respect to the subspace topology.

**Lemma 1:**

*i. The subgroup $P_{\underline{j}}$ preserves $V_j$;*

*ii.* $U_{\underline{j}}L'(\underline{j})$ *acts through the determinant character on* $V_j$.

Proof: Only the second assertion requires a proof. We have $\Xi_\mu \xi = \Xi_{\mu-\beta} d\Xi_{\beta_0} \wedge \ldots \wedge d\Xi_{\beta_{d-1}}$. For $\mu \in B(\underline{j})$ the product $\Xi_{\mu-\beta}$ does not contain any $\Xi_i$ for $i \in \underline{j}$. On the other hand the elements $h \in U_{\underline{j}}L'(\underline{j})$ have columns $i$ for $i \notin \underline{j}$ consisting of zeroes except the entry 1 in position $(i,i)$. It follows that $h_*\Xi_{\mu-\beta} = \Xi_{\mu-\beta}$ for those $h$ and $\mu$. And on $d\Xi_{\beta_0} \wedge \ldots \wedge d\Xi_{\beta_{d-1}}$ such an $h$ acts through multiplication by $\det(h)$ (see the last formula on p. 416 in [ST]). Since $U(\underline{j})$ normalizes $U_{\underline{j}}L'(\underline{j})$ we more generally obtain

$$h_*(u_*(\Xi_\mu \xi)) \ = \ \det(h) \cdot u_*(\Xi_\mu \xi) \text{ for } h \in U_{\underline{j}}L'(\underline{j}),\ u \in U(\underline{j}),\text{ and } \mu \in B(\underline{j})\ . \quad \square$$

In order to compute the space $V_j$ we use the rigid analytic morphism

$$\mathrm{pr}_j : \quad \mathcal{X} \quad \longrightarrow \quad \mathcal{X}^{d+1-j}$$

$$q = [q_0 : \ldots : q_d] \quad \longmapsto \quad [q_j : \ldots : q_d]\ ;$$

here $\mathcal{X}^{d+1-j}$ denotes the $p$-adic symmetric space of the group $GL_{d+1-j}(K)$. This morphism is $P_{\underline{j}}$-equivariant if $P_{\underline{j}}$ acts on $\mathcal{X}$, resp. on $\mathcal{X}^{d+1-j}$, through the inclusion $P_{\underline{j}} \subseteq G$, resp. the projection $P_{\underline{j}} \twoheadrightarrow L(\underline{j}) \cong GL_{d+1-j}(K)$. In section 4 we introduced the irreducible $\mathfrak{p}_{\underline{j}}$-submodule $M_{\underline{j}}$ of $\mathfrak{b}_{\underline{j}}/\mathfrak{b}_{\underline{j}}^>$. For general reasons, it integrates to a rational representation of $P_{\underline{j}}$. We will work with the following explicit model for this representation. Consider an element $g = (g_{rs}) \in L(\underline{j})$. The adjoint action of $g^{-1}$ on any $L_{i\ell} \in \mathfrak{n}_{\underline{j}}^+$, i.e., with $0 \leq i < j \leq \ell \leq d$, is given by

$$\mathrm{ad}(g^{-1})L_{i\ell} \ = \ g_{\ell j}L_{ij} + \ldots + g_{\ell d}L_{id}\ .$$

We may deduce from this that the adjoint action of $L(\underline{j})$ on $U(\mathfrak{g})$ preserves $M_{\underline{j}}^\circ$ as well as $M_{\underline{j}}^\circ \cap \mathfrak{b}_{\underline{j}}^> = U(\mathfrak{n}_{\underline{j}}^+) \cap \mathfrak{b}$. Indeed, the sorting relations $L_{ik}L_{i'\ell} - L_{i\ell}L_{i'k}$ generate $U(\mathfrak{n}_{\underline{j}}^+) \cap \mathfrak{b}_{\underline{j}}^> = U(\mathfrak{n}_{\underline{j}}^+) \cap \mathfrak{b}$ according to Prop. 4.6, and the image of such a relation $\mathrm{ad}(g^{-1})(L_{ik}L_{i'\ell} - L_{i\ell}L_{i'k})$ is a linear combination of sorting relations of the same type involving only $i$ and $i'$ as first indices. It follows that

$$g(\mathfrak{z} + \mathfrak{b}_{\underline{j}}^>) := \mathrm{ad}(g)(\mathfrak{z}) + \mathfrak{b}_{\underline{j}}^> \quad \text{for } \mathfrak{z} \in M_{\underline{j}}^\circ$$

is a well defined action of the group $L(\underline{j})$ on the space $M_{\underline{j}}$. We extend this to a rational representation of $P_{\underline{j}}$ by letting $U_{\underline{j}}L'(\underline{j})$ act through the determinant character. The corresponding derived action of the Lie algebra $\mathfrak{p}_{\underline{j}}$ on $M_{\underline{j}}$ is trivial on $\mathfrak{n}_{\underline{j}}$, is through the trace character on $\mathfrak{l}'(\underline{j})$, and on $\mathfrak{l}(\underline{j})$ is induced by the adjoint action. But in Lemma 4.7 we have seen that this latter action coincides with the left multiplication action.

We now consider the continuous linear map

$$A_j : \Omega^{d-j}(\mathcal{X}^{d+1-j}) \otimes_K M_{\underline{j}} \longrightarrow \Omega^d(\mathcal{X})/\Omega^d(\mathcal{X})^{j+1}$$

$$\eta \otimes (L_\mu + \mathfrak{b}_{\underline{j}}^>) \longmapsto L_\mu \left( \frac{d\Xi_{\beta_0}}{\Xi_{\beta_0}} \wedge \ldots \wedge \frac{d\Xi_{\beta_{j-1}}}{\Xi_{\beta_{j-1}}} \wedge \mathrm{pr}_j^*(\eta) \right) .$$

According to Lemma 1.3 the $P_{\underline{j}}$-action on both sides (diagonally on the left side) is continuous. In the following we will use the abbreviations

$$\xi_{d-j} := \frac{d\Xi_{\beta_j}}{\Xi_{\beta_j}} \wedge \ldots \wedge \frac{d\Xi_{\beta_{d-1}}}{\Xi_{\beta_{d-1}}} \text{ as a } (d-j)\text{-form on } \mathcal{X}^{d+1-j}$$

and

$$\xi^{(j)} := \frac{d\Xi_{\beta_0}}{\Xi_{\beta_0}} \wedge \ldots \wedge \frac{d\Xi_{\beta_{j-1}}}{\Xi_{\beta_{j-1}}} \text{ as a } j\text{-form on } \mathcal{X} .$$

For $g \in G$ we have

$$g_* \xi = \det(g) \cdot \left( \prod_{i=0}^d \frac{\Xi_i}{g_* \Xi_i} \right) \cdot \xi .$$

For $g \in L(\underline{j})$ there is a corresponding formula for $g_* \xi_{d-j}$ and the two together imply

$$g_* \xi = \xi^{(j)} \wedge \mathrm{pr}_j^*(g_* \xi_{d-j}) .$$

By Prop. 3.3 the $u_* \xi_{d-j}$ for $u \in U(\underline{j})$ generate a dense subspace in $\Omega^{d-j}(\mathcal{X}^{d+1-j})$. After we establish various properties of the map $A_j$, this fact will allow us to assume that $\eta = u_* \xi_{d-j}$ for some $u \in U(\underline{j})$. First of all we note that the definition of $A_j$ is independent of the particular representative $L_\mu$ for the coset $L_\mu + \mathfrak{b}_{\underline{j}}^>$ as long as this representative is chosen in $M_{\underline{j}}^o$: For $\mathfrak{z} \in M_{\underline{j}}^o \cap \mathfrak{b}_{\underline{j}}^> = U(\mathfrak{n}_{\underline{j}}^+) \cap \mathfrak{b}$ and $u \in U(\underline{j})$ we have $\mathrm{ad}(u^{-1})(\mathfrak{z}) \in \mathfrak{b}$ and consequently

$$\mathfrak{z}(\xi^{(j)} \wedge \mathrm{pr}_j^*(u_* \xi_{d-j})) = \mathfrak{z}(u_* \xi) = u_*([\mathrm{ad}(u^{-1})(\mathfrak{z})]\xi) = 0 .$$

Next we compute

$$\begin{aligned}
A_j(g_* h_* \xi_{d-j} \otimes \mathrm{ad}(g)(L_\mu)) &= [\mathrm{ad}(g)(L_\mu)](\xi^{(j)} \wedge \mathrm{pr}_j^*(g_* h_* \xi_{d-j})) \\
&= [\mathrm{ad}(g)(L_\mu)](g_* h_* \xi) \\
&= g_*(L_\mu(h_* \xi)) \\
&= g_*(L_\mu(\xi^{(j)} \wedge \mathrm{pr}_j^*(h_* \xi_{d-j}))) \\
&= g_*(A_j(h_* \xi_{d-j} \otimes L_\mu))
\end{aligned}$$

for $g, h \in L(\underline{j})$. This shows that the map $A_j$ is $L(\underline{j})$-equivariant. As special cases of the above identity we have

$$A_j(g_* \xi_{d-j} \otimes L_\mu) = g_*([\mathrm{ad}(g^{-1})(L_\mu)]\xi)$$

and
$$A_j(g_*\xi_{d-j} \otimes \mathrm{ad}(g)(L_\mu)) \;=\; g_*(L_\mu\xi) \;=\; -g_*(\Xi_\mu\xi)$$

for $g \in L(\underline{j})$ and $\mu \in B(\underline{j})$. The former, together with the fact that $M_{\underline{j}}^\circ \cdot \xi \subseteq \sum_{\mu \in B(\underline{j})} K \cdot \Xi_\mu \xi$, shows that the image of $A_j$ is contained in $V_j$. The latter shows that this image is dense in $V_j$. By Lemma 1.ii, the group $U_{\underline{j}} L'(\underline{j})$ acts on the domain of $A_j$, as well as on $V_j$, through the determinant character. Hence $A_j$ in fact is $P_{\underline{j}}$-equivariant.

By Thm. 7.1 the exact $(d-j)$-forms on $\mathcal{X}^{d+1-j}$ coincide with the subspace $\Omega^{d-j}(\mathcal{X}^{d+1-j})^1$. According to Cor. 6.3, this latter space is topologically generated, as an $L(\underline{j})$-representation, by the forms $\Xi_\nu \xi_{d-j}$ for the weights $0 \neq \nu = \sum_{k=j}^{d} n_k \varepsilon_k \in X^*(\overline{T})$. We have $A_j(\Xi_\nu \xi_{d-j} \otimes L_\mu) = L_\mu(\Xi_\nu \xi)$. Let $\mu = \varepsilon_0 + \ldots + \varepsilon_{j-1} - \sum_{k=j}^{d} m_k \varepsilon_k$ with $m_k \geq 0$. By an iteration of the formula (+) in section 4 one has

$$L_\mu(\Xi_\nu \xi) \;=\; c(\mu) \cdot (\prod_{k=j}^{d} \prod_{m=1}^{m_k} (n_k - m)) \cdot \Xi_{\mu+\nu}\xi$$

for some constant $c(\mu) \in K^\times$. There are two cases to distinguish. If $n_k \leq m_k$ for all $j \leq k \leq d$ then we choose a $j \leq \ell \leq d$ such that $n_\ell \geq 1$ and see that the product on the right hand side of the above identity contains the factor 0. Hence $L_\mu(\Xi_\nu \xi) = 0$ in this case. Otherwise there is some $j \leq \ell \leq d$ such that $m_\ell < n_\ell$. Then $J(\mu+\nu) \supseteq \{0, \ldots, j-1, \ell\}$ so that, by Cor. 6.3, $\Xi_{\mu+\nu}\xi$ and hence $L_\mu(\Xi_\nu \xi)$ lies in $\Omega^d(\mathcal{X})^{j+1}$. This shows that $\begin{pmatrix} \text{exact} \\ \text{forms} \end{pmatrix} \otimes M_{\underline{j}}$ lies in the kernel of $A_j$.

**Remark 2:**

$V_j$ is topologically generated by the forms $u_*(\Xi_\mu \xi)$ for $u \in U(\underline{j})$ and $\mu \in B(\underline{j})$.

Proof: We in fact will show that in $\Omega^d(\mathcal{X})$ any form $g_*(\Xi_\mu \xi)$ with $g \in L(\underline{j})$ and $\mu \in B(\underline{j})$ is a (finite) linear combination of forms $u_*(\Xi_\nu \xi)$ with $u \in U(\underline{j})$ and $\nu \in B(\underline{j})$. First of all we have

$$g_*(\Xi_\mu \xi) \;=\; -g_*(L_\mu \xi) \;=\; -[\mathrm{ad}(g)(L_\mu)](g_*\xi)\;.$$

¿From the discussion after Lemma 3.5 we know that $g_*\xi_{d-j}$ is an alternating sum of forms $u_*\xi_{d-j}$ with $u \in U(\underline{j})$. Using the identity $g_*\xi = \xi^{(j)} \wedge \mathrm{pr}_j^*(g_*\xi_{d-j})$ again we see that $g_*\xi$ is an alternating sum of forms $u_*\xi$ with $u \in U(\underline{j})$. Inserting this into the above equation we are reduced to treating a form $[\mathrm{ad}(g)(L_\mu)](u_*\xi) = u_*([\mathrm{ad}(u^{-1}g)(L_\mu)]\xi)$. But $\mathrm{ad}(u^{-1}g)(L_\mu)$ lies in $\sum_{\nu \in B(\underline{j})} K \cdot L_\nu + \mathfrak{b}$.

**Proposition 3:**

*The linear map $A_j$ induces a $P_{\underline{j}}$-equivariant topological isomorphism*

$$[\Omega^{d-j}(\mathcal{X}^{d+1-j})/\begin{array}{c}\text{exact}\\\text{forms}\end{array}]\otimes_K M_{\underline{j}} \xrightarrow{\cong} V_j \;.$$

Proof: So far we know that $A_j$ induces a continuous $P_{\underline{j}}$-equivariant map with dense image between the two sides in the assertion. For simplicity we denote this latter map again by $A_j$. Both sides are Fréchet spaces (the left hand side as a consequence of Thm. 7.1). We claim that it suffices to show that the dual map $A'_j$ is surjective. We only sketch the argument since it is a straightforward nonarchimedean analog of [B-TVS] IV.28, Prop. 3. Let us assume $A'_j$ to be surjective for the moment being. The Hahn-Banach theorem ([Tie] Thm. 3.6) then immediately implies that $A_j$ is injective. Actually $A'_j : V'_j = \mathrm{im}(A_j)' \xrightarrow{\cong} V'$ then is a linear bijection where we abbreviate by $V$ the space on the left hand side of the assertion. This means that $A'_j$ induces a topological isomorphism $\mathrm{im}(A_j)'_s \longrightarrow V'_s$ between the weak dual spaces. Since the Mackey topology ([Tie] p. 282) is defined in terms of the weak dual it follows that $A!_j : V \longrightarrow \mathrm{im}(A_j)$ is a homeomorphism for the Mackey topologies. But on metrizable spaces the Mackey topology coincides with the initial topology ([Tie] Thm. 4.22). Therefore $A_j : V \xrightarrow{\cong} \mathrm{im}(A_j)$ is a topological isomorphism for the initial topologies. With $V$ also $\mathrm{im}(A_j)$ then is complete. Because of the density we have to have $\mathrm{im}(A_j) = V_j$.

Before we establish the surjectivity of $A'_j$ we interrupt the present proof in order to discuss the strong dual of the left hand side in our assertion.

Let

$$\mathrm{St}_{d+1-j} := C^\infty(L(\underline{j})/L(\underline{j})\cap P, K)/C^\infty_{\mathrm{inv}}(L(\underline{j})/L(\underline{j})\cap P, K)$$

denote the Steinberg representation of the group $L(\underline{j})$ equipped with the finest locally convex topology (in particular, $\mathrm{St}_1$ is the trivial character of the group $K^\times$). Recall that identifying $U(\underline{j})$ with the big cell in $L(\underline{j})/L(\underline{j})\cap P$ induces an isomorphism $\mathrm{St}_{d+1-j} \cong C^\infty_o(U(\underline{j}), K)$. We know from Thm. 7.11 that

$$[\Omega^{d-j}(\mathcal{X}^{d+1-j})/\begin{array}{c}\text{exact}\\\text{forms}\end{array}]' \xrightarrow{\cong} C^\infty_o(U(\underline{j}), K) \cong \mathrm{St}_{d+1-j}$$

$$\lambda \longmapsto [u \mapsto \lambda(u_*\xi_{d-j})]$$

is a $L(\underline{j})$-equivariant topological isomorphism. In particular, the strong dual of the left hand side in Prop. 3 carries the finest locally convex topology and may be identified with the space $\mathrm{Hom}_K(M_{\underline{j}}, \mathrm{St}_{d+1-j})$ of all $K$-linear maps from $M_{\underline{j}}$ into $\mathrm{St}_{d+1-j}$. With this identification, the map $A'_j$ becomes the map

$$I^{[j]}_o : V'_j \longrightarrow \mathrm{Hom}_K(M_{\underline{j}}, \mathrm{St}_{d+1-j})$$

$$\lambda \longmapsto \{L_\mu \mapsto [u \mapsto \lambda(L_\mu(u_*\xi))]\}$$

and ist surjectivity will be proved in the course of the proof of Prop. 4 below.

Recall that $M_{\underline{j}}$ is isomorphic to the contragredient of the $j$-th symmetric power $\mathrm{Sym}^j(K^{d+1-j})$ of the standard representation of $L(\underline{j}) \cong GL_{d+1-j}(K)$ on $K^{d+1-j}$.

**Proposition 4:**

i. $V_j$ is a reflexive Fréchet space;

ii. the linear map

$$I_{\mathrm{o}}^{[j]} : V_j' \xrightarrow{\cong} \mathrm{Hom}_K(M_{\underline{j}}, \mathrm{St}_{d+1-j})$$

$$\lambda \longmapsto \{L_\mu \longmapsto [u \mapsto \lambda(L_\mu(u_*\xi))]\}$$

is a $P_{\underline{j}}$-equivariant isomorphism;

iii. the topology of $V_j'$ is the finest locally convex one;

iv. $V_j' \cong \mathrm{St}_{d+1-j} \underset{K}{\otimes} \mathrm{Sym}^j(K^{d+1-j})$ (with $U_{\underline{j}}L'(\underline{j})$ acting on the right hand side through the inverse of the determinant character);

v. $V_j \cong \mathrm{Hom}_K(\mathrm{St}_{d+1-j}, M_{\underline{j}})$ (with the weak topology on the right hand side).

Proof: The first assertion follows by the same argument as for Prop. 6.5. The only other point to establish is the surjectivity of $I_{\mathrm{o}}^{[j]}$. This then settles Prop. 3 which in turn implies the rest of the present assertions by dualizing.

Let $\varphi \in C_{\mathrm{o}}^\infty(U(\underline{j}), K) \cong \mathrm{St}_{d+1-j}$ denote the characteristic function of the compact open subgroup $U(\underline{j}) \cap B$ in $U(\underline{j})$. Since $\mathrm{St}_{d+1-j}$ is an irreducible (in the algebraic sense) $L(\underline{j})$-representation it is generated by $\varphi$ as a $L(\underline{j})$-representation. Hence the finitely many linear maps

$$E_\mu : M_{\underline{j}} \longrightarrow \mathrm{St}_{d+1-j}$$

$$L_\nu \longmapsto \begin{cases} \varphi & \text{if } \nu = \mu, \\ 0 & \text{otherwise} \end{cases}$$

for $\mu \in B(\underline{j})$ generate $\mathrm{Hom}_K(M_{\underline{j}}, \mathrm{St}_{d+1-j})$ as a $L(\underline{j})$-representation. For the surjectivity of $I_{\mathrm{o}}^{[j]}$ it therefore suffices, by $L(\underline{j})$-equivariance, to find a preimage for each $E_\mu$. At the beginning of section 5, we introduced the continuous linear forms

$$\eta \longmapsto \mathrm{Res}_{(\overline{C},0)} \Xi_{-\mu} \eta$$

on $\Omega^d(\mathcal{X})$ for any $\mu \in X^*(\overline{T})$. In terms of the pairing $\langle\,,\,\rangle$ defined before Prop. 5.3 this linear form is given as

$$\eta \longmapsto \langle \eta | U^\circ, f_\mu | B \rangle .$$

We now fix a $\mu \in B(\underline{j})$. Since $f_\mu|B$ has weight $-\mu$ we have $(L_\nu(f_\mu|B))(1) = 0$ for all $\nu \neq \mu$ (compare the proof of Prop. 5.2); in particular $(\mathfrak{z}(f_\mu|B))(1) = 0$ for any $\mathfrak{z} \in \mathfrak{b}_{j+1}$. Taylor's formula then implies that

$$f_\mu|B \in \mathcal{O}(B)^{\mathfrak{b}_{j+1}=0} \ .$$

By Lemma 6.4, the above linear form vanishes on $\Omega^d(\mathcal{X})^{j+1}$ and consequently induces a continuous linear form $\lambda_\mu$ on $V_j$. We compute

$$I_o^{[j]}(\lambda_\mu)(L_\nu)(u) \ = \ \operatorname{Res}_{(\overline{C},0)} \Xi_{-\mu} \cdot L_\nu(u_*\xi) \ = \ \operatorname{Res}_{u^{-1}(\overline{C},0)} \theta$$

with

$$\theta := (u^{-1}\Xi_{-\mu}) \cdot (\operatorname{ad}(u^{-1})(L_\nu))(\xi) \ .$$

Since, by Thm. 7.1, forms in $\Omega^d(\mathcal{X})^1$ have no residues it suffices to determine $\theta$ modulo $\Omega^d(\mathcal{X})^1$. The subspace $M_{\underline{j}}^* := \sum\limits_{\mu \in B(\underline{j})} K \cdot \Xi_{-\mu}$ of $\mathcal{O}(\mathcal{X})$ is $L(\underline{j})$-invariant. In fact, one easily computes that, for $g = (g_{rs}) \in L(\underline{j})$ and $0 \leq i < j \leq \ell \leq d$, one has

$$g_*\Xi_{-(\varepsilon_i-\varepsilon_\ell)} \ = \ g_{j\ell}\Xi_{-(\varepsilon_i-\varepsilon_j)} + \ldots + g_{d\ell}\Xi_{-(\varepsilon_i-\varepsilon_d)} \ .$$

This formula and our previous formula for $\operatorname{ad}(g^{-1})L_{i\ell}$ together show that the pairing

$$M_{\underline{j}} \times M_{\underline{j}}^* \ \longrightarrow \ K$$

$$(L_\mu + \mathfrak{b}_{\underline{j}}^>, \Xi_{-\nu}) \ \longmapsto \ \begin{cases} 1 & \text{if } \mu = \nu, \\ 0 & \text{otherwise} \end{cases}$$

is $L(\underline{j})$-equivariant. It therefore exhibits $M_{\underline{j}}^*$ as the $L(\underline{j})$-representation dual to $M_{\underline{j}}$. The point of this pairing is that, by Cor. 6.3, we have $\Xi_{-\mu} \cdot (L_\nu\xi) = -\Xi_{\nu-\mu}\xi \in \Omega^d(\mathcal{X})^1$ for $\mu \neq \nu$. Applying this together with the equivariance to the above form $\theta$ we obtain that

$$\theta \in \begin{cases} -\xi + \Omega^d(\mathcal{X})^1 & \text{if } \mu = \nu, \\ \Omega^d(\mathcal{X})^1 & \text{if } \mu \neq \nu \end{cases}$$

and consequently that

$$I_o^{[j]}(\lambda_\mu)(L_\nu)(u) \ = \ \begin{cases} -\operatorname{Res}_{u^{-1}(\overline{C},0)} \xi & \text{if } \mu = \nu, \\ 0 & \text{if } \mu \neq \nu \ . \end{cases}$$

By [St] Lemma 23 the form $\xi$ has residues only on the standard apartment and those are equal to $\pm 1$. The chamber $u^{-1}\overline{C}$ lies in the standard apartment if and only if $u$ fixes $\overline{C}$. It follows that $I_o^{[j]}(\lambda_\mu)(L_\mu)$ is supported on $U(\underline{j}) \cap B$ where it is a constant function with value $\pm 1$. All in all we see that

$$I_o^{[j]}(\lambda_\mu) \ = \ \pm E_\mu$$

(the sign depending on the parity of $d$). $\square$

The natural $P_{\underline{j}}$-equivariant linear map

$$[\Omega^d(\mathcal{X})^j/\Omega^d(\mathcal{X})^{j+1}]' \longrightarrow V_j'$$

is surjective (by Hahn-Banach) and is strict (by the same argument as for Prop. 6.7). Moreover both sides are inductive limits of sequences of Banach spaces (see the proof of Prop. 6.5) and are locally analytic $P_{\underline{j}}$-representations in the sense of Cor. 6.8. Therefore the assumptions of the Frobenius reciprocity theorem 4.2.6 in [Fea] are satisfied and we obtain the $G$-equivariant continuous linear map

$$I^{[j]} : [\Omega^d(\mathcal{X})^j/\Omega^d(\mathcal{X})^{j+1}]' \longrightarrow C^{\mathrm{an}}(G, P_{\underline{j}}; V_j')$$
$$\lambda \longmapsto [g \mapsto (g^{-1}\lambda)|V_j] \, .$$

Here $C^{\mathrm{an}}(G, P_{\underline{j}}; V_j')$ – the "induced representation in the locally analytic sense" – denotes the vector space of all locally analytic maps $f : G \longrightarrow V_j'$ such that $f(gh) = h^{-1}(f(g))$ for any $g \in G$ and $h \in P_{\underline{j}}$ on which $G$ acts by left translations. Its natural locally convex topology is constructed in [Fea] 4.1.3 (to avoid confusion we should point out that [Fea] uses a more restrictive notion of a $V$-valued locally analytic map but which coincides with the notion from Bourbaki provided $V$ is quasi-complete – see loc.cit. 2.1.4 and 2.1.7).

**Lemma 5:**

$I^{[j]}$ *is injective.*

Proof: It is an immediate consequence of Cor. 6.3 that $\sum_{g \in G} g(V_j)$ is dense in $\Omega^d(\mathcal{X})^j/\Omega^d(\mathcal{X})^{j+1}$. $\square$

In order to describe the image of $I^{[j]}$ we first need to understand in which sense we can impose left invariant differential equations on vectors in an induced representation. For any Hausdorff locally convex $K$-vector space $V$ the right translation action of $G$ on $C^{\mathrm{an}}(G, V) := C^{\mathrm{an}}(G, \{1\}; V)$ is differentiable and induces an action of $U(\mathfrak{g})$ by left invariant and continuous operators ([Fea] 3.3.4). For $V := V_j' \cong \mathrm{Hom}_K(M_{\underline{j}}, \mathrm{St}_{d+1-j})$ we therefore may consider the $K$-bilinear map

$$\langle \, , \, \rangle : (U(\mathfrak{g}) \underset{K}{\otimes} M_{\underline{j}}) \times C^{\mathrm{an}}(G, \mathrm{Hom}_K(M_{\underline{j}}, \mathrm{St}_{d+1-j})) \longrightarrow C^{\mathrm{an}}(G, \mathrm{St}_{d+1-j})$$
$$(\mathfrak{z} \otimes m, f) \longmapsto [g \mapsto (\mathfrak{z}f)(g)(m)] \, .$$

Note that, for a fixed $\mathfrak{z} \in U(\mathfrak{g}) \underset{K}{\otimes} M_{\underline{j}}$, the "differential operator"

$$\langle \mathfrak{z}, \ \rangle : C^{\mathrm{an}}(G, \mathrm{Hom}_K(M_{\underline{j}}, \mathrm{St}_{d+1-j})) \longrightarrow C^{\mathrm{an}}(G, \mathrm{St}_{d+1-j})$$

is continuous and $G$-equivariant (for the left translation actions). The action of $P_{\underline{j}}$ on $\mathrm{Hom}_K(M_{\underline{j}}, \mathrm{St}_{d+1-j}) = M'_{\underline{j}} \underset{K}{\otimes} \mathrm{St}_{d+1-j}$ is differentiable and the derived action of $\mathfrak{p}_{\underline{j}}$ is given by

(1) $$(\mathfrak{x} E)(m) = -E(\mathfrak{x} m)$$

for $\mathfrak{x} \in \mathfrak{p}_{\underline{j}}, E \in \mathrm{Hom}_K(M_{\underline{j}}, \mathrm{St}_{d+1-j})$, and $m \in M_{\underline{j}}$. This is immediate from the fact that any vector in $\mathrm{St}_{d+1-j}$ is fixed by an open subgroup of $P_{\underline{j}}$ so that the derived action of $\mathfrak{p}_{\underline{j}}$ on $\mathrm{St}_{d+1-j}$ is trivial.

Now recall that the induced representation $C^{\mathrm{an}}(G, P_{\underline{j}}; \mathrm{Hom}_K(M_{\underline{j}}, \mathrm{St}_{d+1-j}))$ is the closed subspace of $C^{\mathrm{an}}(G, \mathrm{Hom}_K(M_{\underline{j}}, \mathrm{St}_{d+1-j}))$ of all those maps $f$ which satisfy $f(gh) = h^{-1}(f(g))$ for $g \in G$ and $h \in P_{\underline{j}}$. For such an $f$ we therefore have

$$\begin{aligned}(\mathfrak{x} f)(g) &= \tfrac{d}{dt} f(g \exp(t\mathfrak{x}))\big|_{t=0} \\ &= \tfrac{d}{dt} \exp(t\mathfrak{x})^{-1}(f(g))\big|_{t=0} \\ &= -\mathfrak{x}(f(g))\end{aligned}$$

for $\mathfrak{x} \in \mathfrak{p}_{\underline{j}}$ and slightly more generally

(2) $$\begin{aligned}(\mathfrak{z}(\mathfrak{x} f))(g) &= \tfrac{d}{dt}(\mathfrak{x} f)(g \exp(t\mathfrak{z}))\big|_{t=0} \\ &= -\tfrac{d}{dt}\mathfrak{x}(f(g \exp(t\mathfrak{z})))\big|_{t=0} \\ &= -\mathfrak{x}(\tfrac{d}{dt} f(g \exp(t\mathfrak{z}))\big|_{t=0}) \\ &= -\mathfrak{x}((\mathfrak{z} f)(g))\end{aligned}$$

for $\mathfrak{x} \in \mathfrak{p}_{\underline{j}}$ and $\mathfrak{z} \in \mathfrak{g}$; the third equality is a consequence of the continuity of the operator $\mathfrak{x}$. Combining (1) and (2) we obtain

$$(\mathfrak{z}(\mathfrak{x} f))(g)(m) = (-\mathfrak{x}((\mathfrak{z} f)(g)))(m) = ((\mathfrak{z} f)(g))(\mathfrak{x} m)$$

or equivalently

$$\langle \mathfrak{z}\mathfrak{x} \otimes m, f \rangle = \langle \mathfrak{z} \otimes \mathfrak{x} m, f \rangle$$

for $f \in C^{\mathrm{an}}(G, P_{\underline{j}}; \mathrm{Hom}_K(M_{\underline{j}}, \mathrm{St}_{d+1-j}))$, $m \in M_{\underline{j}}$, $\mathfrak{x} \in \mathfrak{p}_{\underline{j}}$, and $\mathfrak{z} \in \mathfrak{g}$. This means that the above pairing restricts to a pairing

$$\langle \, , \, \rangle : (U(\mathfrak{g}) \underset{U(\mathfrak{p}_{\underline{j}})}{\otimes} M_{\underline{j}}) \times C^{\mathrm{an}}(G, P_{\underline{j}}; \mathrm{Hom}_K(M_{\underline{j}}, \mathrm{St}_{d+1-j})) \longrightarrow C^{\mathrm{an}}(G, \mathrm{St}_{d+1-j})$$

and enables us to consider, for any subset $\mathfrak{d} \subseteq U(\mathfrak{g}) \underset{U(\mathfrak{p}_{\underline{j}})}{\otimes} M_{\underline{j}}$, the $G$-invariant closed subspace

$$C^{\mathrm{an}}(G, P_{\underline{j}}; \mathrm{Hom}_K(M_{\underline{j}}, \mathrm{St}_{d+1-j}))^{\mathfrak{d}=0} := \\ \{f \in C^{\mathrm{an}}(G, P_{\underline{j}}; \mathrm{Hom}_K(M_{\underline{j}}, \mathrm{St}_{d+1-j})) \ : \ \langle \mathfrak{z}, f \rangle = 0 \text{ for any } \mathfrak{z} \in \mathfrak{d}\} \ .$$

The relevant subset for our purposes is the kernel

$$\mathfrak{d}_{\underline{j}} = \ker(U(\mathfrak{g}) \underset{U(\mathfrak{p}_{\underline{j}})}{\otimes} M_{\underline{j}} \longrightarrow \mathfrak{b}_{\underline{j}}/\mathfrak{b}_{\underline{j}}^{>})$$

of the natural surjection sending $\mathfrak{z} \otimes m$ to $\mathfrak{z}m$. By the Poincaré-Birkhoff-Witt theorem the inclusion $U(\mathfrak{n}_{\underline{j}}^+) \subseteq U(\mathfrak{g})$ induces an isomorphism $U(\mathfrak{n}_{\underline{j}})^+ \underset{K}{\otimes} M_{\underline{j}} \overset{\cong}{\longrightarrow} U(\mathfrak{g}) \underset{U(\mathfrak{p}_{\underline{j}})}{\otimes} M_{\underline{j}}$. We mostly will view $\mathfrak{d}_{\underline{j}}$ as a subspace of $U(\mathfrak{n}_{\underline{j}}^+) \underset{K}{\otimes} M_{\underline{j}}$.

**Theorem 6:**

*The map $I^{[j]}$ (together with $I_{\mathrm{o}}^{[j]}$) induces a $G$-equivariant topological isomorphism*

$$I^{[j]} : [\Omega^d(\mathcal{X})^j/\Omega^d(\mathcal{X})^{j+1}]' \overset{\cong}{\longrightarrow} C^{\mathrm{an}}(G, P_{\underline{j}}; \mathrm{Hom}_K(M_{\underline{j}}, \mathrm{St}_{d+1-j}))^{\mathfrak{d}_{\underline{j}}=0}$$

$$\lambda \longmapsto [g \mapsto I_{\mathrm{o}}^{[j]}((g^{-1}\lambda)|V_j)] \ .$$

Proof: We start by showing that the image of $I^{[j]}$ satisfies the relations $\mathfrak{d}_{\underline{j}} = 0$. Let $\mathfrak{z} = \underset{\mu \in B(\underline{j})}{\sum} \mathfrak{z}_{(\mu)} \otimes L_\mu \in \mathfrak{d}_{\underline{j}} \subseteq U(\mathfrak{n}_{\underline{j}}^+) \underset{K}{\otimes} M_{\underline{j}}$; then $\mathfrak{z} = \underset{\mu}{\sum} \mathfrak{z}_{(\mu)} L_\mu \in U(\mathfrak{n}_{\underline{j}}^+) \cap \mathfrak{b}_{\underline{j}}^{>} = U(\mathfrak{n}_{\underline{j}}^+) \cap \mathfrak{b}$ (Prop. 4.6.iii). Note that

$$[I^{[j]}(\lambda)(g)](L_\mu)(u) \ = \ (g^{-1}\lambda)(L_\mu(u_*\xi)) \ = \ \lambda(g_*(L_\mu(u_*\xi)))$$

for $g \in G$, $\mu \in B(\underline{j})$, and $u \in U(\underline{j})$. We compute

$$\langle \mathfrak{z}, I^{[j]}(\lambda) \rangle(g)(u) \ = \ \sum_\mu [(\mathfrak{z}_{(\mu)}(I^{[j]}(\lambda)))(g)](L_\mu)(u)$$

$$= \ \sum_\mu \lambda(g_*(\mathfrak{z}_{(\mu)} L_\mu(u_*\xi)))$$

$$= \ \lambda(g_*(\mathfrak{z}(u_*\xi)))$$

$$= \ \lambda(g_* u_*((\mathrm{ad}(u^{-1})(\mathfrak{z}))\xi))$$

which is zero because $U(\mathfrak{n}_{\underline{j}}^+) \cap \mathfrak{b}$ is $\mathrm{ad}(U(\underline{j}))$-invariant as we have seen earlier in this section.

We know already that $I^{[j]}$ is continuous, $G$-equivariant, and injective. Next we establish surjectivity. Let $f$ be a map in the right hand side of the assertion. By a series of simplifications we will show that it suffices to consider an $f$ of a very particular form for which we then will exhibit an explicit preimage under $I^{[j]}$. We show first that we may assume that

- $f$ is supported on $BP_{\underline{j}}$ and
- $f|U_{\underline{j}}^+ \cap B$ is analytic (not merely locally analytic).

By the Iwasawa decomposition we have the finite disjoint open covering

$$G/P_{\underline{j}} = \dot{\bigcup_g} gBP_{\underline{j}}/P_{\underline{j}}$$

where $g$ runs through a set of representatives for the cosets in $GL_{d+1}(o)/B$. As before let $U_{\underline{j}}^+$ denote the transpose of $U_{\underline{j}}$. Then $U_{\underline{j}}^{(0)} := U_{\underline{j}}^+ \cap B$ is the congruence subgroup of all matrices in $U_{\underline{j}}^+$ whose non-diagonal entries are integral multiples of $\pi$. Consider the higher congruence subgroups $U_{\underline{j}}^{(n)}$, for $n \geq 0$, of all matrices in $U_{\underline{j}}^+$ whose non-diagonal entries are integral multiples of $\pi^n$. These $U_{\underline{j}}^{(n)}$ are polydisks in an obvious way, and we have $U_{\underline{j}}^{(n)} = y^n(U_{\underline{j}}^+ \cap B)y^{-n}$ where $y \in G$ is the diagonal matrix with entries $(\pi, \ldots, \pi, 1, \ldots, 1)$. The Iwahori decomposition for $B$ implies that the map

$$\begin{array}{rcl} gU_{\underline{j}}^{(0)} & \xrightarrow{\sim} & gBP_{\underline{j}}/P_{\underline{j}} \\ gu & \longmapsto & guP_{\underline{j}} \end{array}$$

is a homeomorphism. Our map $f$ restricted to $gU_{\underline{j}}^{(0)}$ still only is locally analytic. But we find a sufficiently big $n \in \mathbb{N}$ such that $f|ghU_{\underline{j}}^{(n)}$ is analytic for all $g$ as above and all $h$ in a system of representatives for the cosets in $U_{\underline{j}}^{(0)}/U_{\underline{j}}^{(n)}$. If we put

$$f_{g,h} := ((gh)^{-1}f)|U_{\underline{j}}^{(n)}P_{\underline{j}} \text{ extended by zero}$$

then these maps lie in the right hand side of our assertion and we have

$$f = \sum_{g,h}(gh)f_{g,h} \;.$$

The reason for this of course is that

$$G = \dot{\bigcup_{g,h}} ghU_{\underline{j}}^{(n)}P_{\underline{j}}$$

is a disjoint finite open covering. By linearity and $G$-equivariance of $I^{[j]}$ it therefore suffices to find a preimage for each $f_{g,h}$. This means we may assume that our map $f$ is supported on $U_{\underline{j}}^{(n)} P_{\underline{j}}$ and is analytic on $U_{\underline{j}}^{(n)}$. Using $G$-equivariance again, we may translate $f$ by $y^{-n}$ so that it has the desired properties.
For our next reduction, we will show that we may further assume that

- $f$ is supported on $BP_{\underline{j}}$ with $f|U_{\underline{j}}^+ \cap B = \varepsilon \otimes \varphi$ for some $\varepsilon \in \mathcal{O}(U_{\underline{j}}^+ \cap B, M'_{\underline{j}})^{\mathfrak{d}_{\underline{j}}=0}$ and $\varphi \in \mathrm{St}_{d+1-j}$.

If we consider an analytic map on $U_{\underline{j}}^+ \cap B$ with values in the locally convex vector space $\mathrm{Hom}_K(M_{\underline{j}}, \mathrm{St}_{d+1-j})$ then the coefficients in its power series expansion multiplied by appropriate powers of $\pi$ form a bounded subset of $\mathrm{Hom}_K(M_{\underline{j}}, \mathrm{St}_{d+1-j})$. The topology of that vector space is the finest locally convex one. Hence any bounded subset and therefore the set of coefficients lies in a finite dimensional subspace. This means that our $f|U_{\underline{j}}^+ \cap B$ is an element of $\mathcal{O}(U_{\underline{j}}^+ \cap B) \underset{K}{\otimes} \mathrm{Hom}_K(M_{\underline{j}}, \mathrm{St}_{d+1-j})$. Moreover, viewing $\mathfrak{d}_{\underline{j}}$ as a subspace of $U(\mathfrak{n}_{\underline{j}}^+) \underset{K}{\otimes} M_{\underline{j}}$ it is clear that with respect to the obvious pairing

$$\langle\,,\,\rangle : (U(\mathfrak{n}_{\underline{j}}^+) \underset{K}{\otimes} M_{\underline{j}}) \times (\mathcal{O}(U_{\underline{j}}^+ \cap B) \underset{K}{\otimes} \mathrm{Hom}_K(M_{\underline{j}}, \mathrm{St}_{d+1-j}))$$

$$\longrightarrow \mathcal{O}(U_{\underline{j}}^+ \cap B) \underset{K}{\otimes} \mathrm{St}_{d+1-j}$$

$$(\mathfrak{z} \otimes m, e \otimes E) \longmapsto \mathfrak{z} e \otimes E(m)$$

we have $\langle \mathfrak{d}_{\underline{j}}, f|U_{\underline{j}}^+ \cap B\rangle = 0$. We now decompose

$$f|U_{\underline{j}}^+ \cap B = \sum_i e_i \otimes E_i$$

into a finite sum with $e_i \in \mathcal{O}(U_{\underline{j}}^+ \cap B)$ and $E_i \in \mathrm{Hom}_K(M_{\underline{j}}, \mathrm{St}_{d+1-j})$ such that the images $E_i(M_{\underline{j}})$ are linearly independent 1-dimensional subspaces of $\mathrm{St}_{d+1-j}$. Then each $e_i \otimes E_i$ satisfies the relations $\langle \mathfrak{d}_{\underline{j}}, e_i \otimes E_i\rangle = 0$. We define maps $f_i$ on $G$ with values in $\mathrm{Hom}_K(M_{\underline{j}}, \mathrm{St}_{d+1-j})$ by setting

$$f_i(uh) := e_i(u) \cdot h^{-1}(E_i) \text{ for } u \in U_{\underline{j}}^+ \cap B \text{ and } h \in P_{\underline{j}}$$

and extending this by zero to $G$. Since the map $h \mapsto h^{-1}(E_i)$ is locally analytic on $P_{\underline{j}}$ it easily follows that $f_i \in C^{\mathrm{an}}(G, P_{\underline{j}}; \mathrm{Hom}_K(M_{\underline{j}}, \mathrm{St}_{d+1-j}))$. By construction $f_i$ is supported on $BP_{\underline{j}}$ with $f_i|U_{\underline{j}}^+ \cap B = e_i \otimes E_i$. Clearly

$$f = \sum_i f_i \ .$$

We claim that each $f_i$ satisfies the relations $\mathfrak{d}_{\underline{j}} = 0$. This will be a consequence of the following observation. The group $P_{\underline{j}}$ acts diagonally on $U(\mathfrak{g}) \underset{U(\mathfrak{p}_{\underline{j}})}{\otimes} M_{\underline{j}}$ via $h(\mathfrak{z} \otimes m) := \mathrm{ad}(h)\mathfrak{z} \otimes hm$. The point to observe is that the subspace $\mathfrak{d}_{\underline{j}}$ is $P_{\underline{j}}$-invariant. Note first that because $U(\mathfrak{n}_{\underline{j}}^+) \cap \mathfrak{b}_{\underline{j}}^> \subseteq \mathfrak{b}$ (Prop. 4.6. iii) an element $\sum_\mu \mathfrak{z}_{(\mu)} \otimes L_\mu \in U(\mathfrak{n}_{\underline{j}}^+) \underset{K}{\otimes} M_{\underline{j}}$ lies in $\mathfrak{d}_{\underline{j}}$ if and only if $\sum_\mu \mathfrak{z}_{(\mu)} L_\mu \xi = 0$. Let now $\sum_\mu \mathfrak{z}_{(\mu)} \otimes L_\mu \in \mathfrak{d}_{\underline{j}} \subseteq U(\mathfrak{n}_{\underline{j}}^+) \underset{K}{\otimes} M_{\underline{j}}$ and $h \in P_{\underline{j}}$. We distinguish two cases. If $h \in L(\underline{j})$ then using the $\mathrm{ad}(L(\underline{j}))$-invariance of $U(\mathfrak{n}_{\underline{j}}^+) \cap \mathfrak{b}$ we obtain

$$(\sum_\mu \mathrm{ad}(h)(\mathfrak{z}_{(\mu)}) \cdot hL_\mu)\xi = (\mathrm{ad}(h)(\sum_\mu \mathfrak{z}_{(\mu)} L_\mu))\xi = 0 .$$

If $h \in L'(\underline{j})U_{\underline{j}}$ then using Lemma 1.ii we obtain

$$(\sum_\mu \mathrm{ad}(h)(\mathfrak{z}_{(\mu)}) \cdot hL_\mu)\xi = \det(h) \cdot h_*(\sum_\mu \mathfrak{z}_{(\mu)} h_*^{-1} L_\mu \xi)$$

$$= h_*(\sum_\mu \mathfrak{z}_{(\mu)} L_\mu \xi) = 0 .$$

Going back to our maps $f_i$ and letting again $\sum_\mu \mathfrak{z}_{(\mu)} \otimes L_\mu \in \mathfrak{d}_{\underline{j}} \subseteq U(\mathfrak{n}_{\underline{j}}^+) \underset{K}{\otimes} M_{\underline{j}}$ we now compute

$$((\sum_\mu \mathfrak{z}_{(\mu)} \otimes L_\mu) f_i)(uh) = \sum_\mu (\mathfrak{z}_{(\mu)} f_i)(uh)(L_\mu)$$

$$= \sum_\mu ((\mathrm{ad}(h)\mathfrak{z}_{(\mu)})e_i)(u) \cdot h^{-1}(E_i)(L_\mu)$$

$$= h^{-1}(\sum_\mu ((\mathrm{ad}(h)\mathfrak{z}_{(\mu)})e_i)(u) \cdot E_i(hL_\mu))$$

$$= h^{-1}(\langle \sum_\mu \mathrm{ad}(h)\mathfrak{z}_{(\mu)} \otimes hL_\mu, e_i \otimes E_i \rangle(u))$$

$$= h^{-1}(\langle h(\sum_\mu \mathfrak{z}_{(\mu)} \otimes L_\mu), e_i \otimes E_i \rangle(u))$$

$$= 0 .$$

This establishes our claim.

We want to further normalize the component $\varphi$ in this last expression. Let $\varphi_\mathrm{o} \in C_\mathrm{o}^\infty(U(\underline{j}), K) \cong \mathrm{St}_{d+1-j}$ denote the characteristic function of $U(\underline{j}) \cap B$. Then $\varphi$ can be written as a linear combination of vectors of the form $g^{-1}\varphi_\mathrm{o}$ with $g \in L(\underline{j})$. A straightforward argument shows that $f$ can be decomposed accordingly so that we may assume $\varphi = g^{-1}\varphi_\mathrm{o}$ for some $g \in L(\underline{j})$. We now find a finite disjoint open covering

$$g(U_{\underline{j}}^+ \cap B)P_{\underline{j}} = \dot{\bigcup_i} u_i y^n (U_{\underline{j}}^+ \cap B)P_{\underline{j}}$$

with appropriate $n \in \mathbb{N}$ and $u_i \in U_{\underline{j}}^+$. The map $gf$ is supported on $gBP_{\underline{j}}$ and its restriction $gf|g(U_{\underline{j}}^+ \cap B)g^{-1}$ is analytic with values in $M'_{\underline{j}} \otimes K\varphi_o$. If we put

$$f_i := ((u_i y^n)^{-1} gf)|BP_{\underline{j}} \text{ extended by zero}$$

then these maps lie in the induced representation on the right hand side of our assertion and we have

$$f = \sum_i g^{-1} u_i y^n f_i \ .$$

The restriction of $f_i$ to $U_{\underline{j}}^+ \cap B$ satisfies

$$f_i(u) = (gf)(u_i y^n u) = \pi^{-jn} \cdot (gf)(u_i y^n u y^{-n}) \ .$$

But $u_i y^n u y^{-n} \in gBP_{\underline{j}} \cap U_{\underline{j}}^+ \subseteq g(U_{\underline{j}}^+ \cap B)g^{-1}$. It follows that $f_i|U_{\underline{j}}^+ \cap B$ is analytic with values in $M'_{\underline{j}} \otimes K\varphi_o$. At this point we have arrived at the conclusion that we may assume that

- $f$ is supported on $BP_{\underline{j}}$ with $f|U_{\underline{j}}^+ \cap B = \varepsilon \otimes \varphi_o$ for some $\varepsilon \in \mathcal{O}(U_{\underline{j}}^+ \cap B, M'_{\underline{j}})^{\partial_{\underline{j}}=0}$.

We rephrase the above discussion in the following way. We have the linear map

$$\text{Ext}_{\underline{j}} : \mathcal{O}(U_{\underline{j}}^+ \cap B, M'_{\underline{j}})^{\partial_{\underline{j}}=0} \longrightarrow C^{\text{an}}(G, P_{\underline{j}}; \text{Hom}_K(M_{\underline{j}}, \text{St}_{d+1-j}))^{\partial_{\underline{j}}=0}$$

defined by

$$\text{Ext}_{\underline{j}}(\varepsilon)(g) := \begin{cases} h^{-1}(\varepsilon(u) \otimes \varphi_o) & \text{for } g = uh \text{ with } u \in U_{\underline{j}}^+ \cap B, h \in P_{\underline{j}} \ , \\ 0 & \text{otherwise} \ . \end{cases}$$

Its image generates the right hand side (algebraically) as a $G$-representation. An argument analogous to the proof of [Fea] 4.3.1 shows that $\text{Ext}_{\underline{j}}$ is continuous. On the other hand, in section 6 after Lemma 4 we had constructed a continuous linear map

$$D_{\underline{j}} : \mathcal{O}(U_{\underline{j}}^+ \cap B, M'_{\underline{j}})^{\partial_{\underline{j}}=0} \longrightarrow [\Omega^d(\mathcal{X})^j/\Omega^d(\mathcal{X})^{j+1}]' \ .$$

The surjectivity of $I^{[j]}$ therefore will follow from the identity

$$\text{Ext}_{\underline{j}} = I^{[j]} \circ D_{\underline{j}} \ .$$

By the continuity of all three maps involved it suffices to check this identity on weight vectors. Fix a weight $\nu$ with $J(\nu) = \{0, \ldots, j-1\}$. By construction the map $D_{\underline{j}}$ sends the weight vector $\sum_{\mu \in B(\underline{j})} [(L_\mu f_\nu)|U_{\underline{j}}^+ \cap B] \otimes L_\mu^*$ to the linear form $\lambda_\nu(\eta) = \mathrm{Res}_{(\overline{C},0)} \Xi_{-\nu}\eta$. What we therefore have to check is that $I^{[j]}(\lambda_\nu)$ is supported on $BP_{\underline{j}}$ with

$$I^{[j]}(\lambda_\nu)|U_{\underline{j}}^+ \cap B = \sum_{\mu \in B(\underline{j})} [(L_\mu f_\nu)|U_{\underline{j}}^+ \cap B] \otimes L_\mu^* \otimes \varphi_{\mathrm{o}}\,.$$

By definition we have

$$\begin{aligned}
[I^{[j]}(\lambda_\nu)(g)](L_\mu)(u) &= [(g^{-1}\lambda_\nu)|V_j](L_\mu(u_*\xi)) \\
&= \lambda_\nu(g(L_\mu(u_*\xi))) \\
&= \mathrm{Res}_{(\overline{C},0)} \Xi_{-\nu} \cdot g_*(L_\mu(u_*\xi)) \\
&= \mathrm{Res}_{(\overline{C},0)} \Xi_{-\nu} \cdot g_* u_*((\mathrm{ad}(u^{-1})(L_\mu))\xi)
\end{aligned}$$

for $\mu \in B(\underline{j})$ and $u \in U(\underline{j}) \subseteq P_{\underline{j}}$. First we deal with the vanishing of this expression for $g \notin BP_{\underline{j}}$. Observe that

- $g \notin BP_{\underline{j}}$ if and only if $gu \notin BP_{\underline{j}}$, and
- $\mathrm{ad}(u^{-1})(L_\mu)\xi \in \sum_{\mu' \in B(\underline{j})} K \cdot \Xi_{\mu'}\xi$.

Hence it suffices to show that

$$\mathrm{Res}_{(\overline{C},0)} \Xi_{-\nu} \cdot g_*(\Xi_\mu \xi) = 0 \text{ for } g \notin BP_{\underline{j}}\,.$$

We distinguish two cases. First assume that $g \notin U_{\underline{j}}^+ P_{\underline{j}}$. The divisor $\mathrm{div}(\Xi_{-\nu} \cdot g_*(\Xi_\mu \xi))_\infty$ is supported among the hyperplanes $\Xi_0 = 0, \ldots, \Xi_{j-1} = 0, g_*\Xi_j = 0, \ldots, g_*\Xi_d = 0$. Those are linearly dependent if $g \notin U_{\underline{j}}^+ P_{\underline{j}}$ and hence have a nonempty intersection, i.e., $Z(\Xi_{-\nu} \cdot g_*(\Xi_\mu \xi)) \neq \emptyset$. According to the discussion after Prop. 6.2 the form $\Xi_{-\nu} \cdot g_*(\Xi_\mu \xi)$ therefore lies in $\Omega_{\mathrm{alg}}^d(\mathcal{X})^1$, hence is exact by Lemma 7.2, and consequently has zero residue. Second we consider the case $g \in U_{\underline{j}}^+ \setminus (U_{\underline{j}}^+ \cap B)$. Then $g$ fixes $\Xi_0, \ldots, \Xi_{j-1}$ so that $g^{-1}\Xi_{-\nu}$ is a linear combination of $\Xi_{-\nu'}$ with $J(\nu') \subseteq \{0, \ldots, j-1\}$. It follows that $\Xi_{-\nu} \cdot g_*(\Xi_\mu \xi)$ is a linear combination of forms $\Xi_{\nu''}\xi$ among which the only possible non-exact one is $\xi$ (compare the proof of Lemma 7.2). We obtain

$$\mathrm{Res}_{(\overline{C},0)} \Xi_{-\nu} \cdot g_*(\Xi_\mu \xi) = \mathrm{Res}_{g^{-1}(\overline{C},0)}(g^{-1}\Xi_{-\nu})\Xi_\mu \xi = c \cdot \mathrm{Res}_{g^{-1}(\overline{C},0)} \xi$$

with some constant $c \in K$. But $\xi$ has residues only on the standard apartment and $g^{-1}(\overline{C}, 0)$ lies in the standard apartment only if $g \in U_{\underline{j}}^+ \cap B$. This establishes the assertion about the support of $I^{[j]}(\lambda_\nu)$.

Fix now a $g \in U_{\underline{j}}^+ \cap B$ and let $u \in U(\underline{j})$. Repeating the last argument for $gu$ instead of $g$ we obtain that $\operatorname{Res}_{(\overline{C},0)} \Xi_{-\nu} \cdot g_* u_*(\Xi_\mu \xi) = 0$ unless $gu$ and hence $u$ fixes $(\overline{C}, 0)$. This means that, for $g \in U_{\underline{j}}^+ \cap B$, the function $[I^{[j]}(\lambda_\nu)(g)](L_\mu) \in C_\circ^\infty(U(\underline{j}), K)$ vanishes outside $U(\underline{j}) \cap B$. For $u \in U(\underline{j}) \cap B$ we have

$$[I^{[j]}(\lambda_\nu)(g)](L_\mu)(u) = \operatorname{Res}_{(\overline{C},0)}(u^{-1}g^{-1}\Xi_{-\nu})((\operatorname{ad}(u^{-1})(L_\mu))\xi)$$

$$= \sum_{J(\nu') \subseteq \underline{j}} c(\nu') \operatorname{Res}_{(\overline{C},0)}(u^{-1}\Xi_{-\nu'})((\operatorname{ad}(u^{-1})(L_\mu))\xi)$$

where

$$g^{-1}\Xi_{-\nu} = \sum_{J(\nu') \subseteq \underline{j}} c(\nu')\Xi_{-\nu'} .$$

If $\nu' \in B(\underline{j})$ then we computed the corresponding summand already in the proof of Prop. 8.4 and, in particular, showed that it is independent of $u \in U(\underline{j}) \cap B$. On the other hand the subspace $\sum_{\substack{J(\nu') \subseteq \underline{j} \\ \nu' \notin B(\underline{j})}} K \cdot \Xi_{-\nu'}$ of $\mathcal{O}(\mathcal{X})$ is preserved by the action of $U(\underline{j})$. This means that, for $\nu' \notin B(\underline{j})$, the form $(u^{-1}\Xi_{-\nu'})((\operatorname{ad}(u^{-1})(L_\mu))\xi)$ cannot contain $\xi$ and therefore must have zero residue. This computation says that, for fixed $g \in U_{\underline{j}}^+ \cap B$ and fixed $\mu \in B(\underline{j})$, the function $[I^{[j]}(\lambda_\nu)(g)](L_\mu)(u)$ is constant in $u \in U(\underline{j}) \cap B$. In other words we have

$$I^{[j]}(\lambda_\nu)(g) = \sum_{\mu \in B(\underline{j})} [I^{[j]}(\lambda_\nu)(g)](L_\mu)(1) \otimes L_\mu^* \otimes \varphi_\circ$$

for $g \in U_{\underline{j}}^+ \cap B$. But using the various definitions we compute

$$[I^{[j]}(\lambda_\nu)(g)](L_\mu)(1) = \operatorname{Res}_{(\overline{C},0)} \Xi_{-\nu} \cdot g_*(L_\mu \xi) = (L_\mu f_\nu)(g) .$$

This establishes the surjectivity and hence bijectivity of the map $I^{[j]}$. Finally, that $I^{[j]}$ is open and hence a topological isomorphism is a consequence of the open mapping theorem in the form given in [GK] Thm. 3.1($A_3$) provided we show that both sides of $I^{[j]}$ are (LB)-spaces, i.e., a locally convex inductive limit of a sequence of Banach spaces. For the left hand side this fact is implicitly contained in our earlier arguments: In the proof of Prop. 6.5 we had noted that $\Omega^d(\mathcal{X})^j/\Omega^d(\mathcal{X})^{j+1}$ is the projective limit of a sequence of Banach spaces with compact transition maps. We certainly may assume in addition that these transition maps have dense images. By the same argument as in the proof of Prop. 2.4 it then follows that the strong dual $[\Omega^d(\mathcal{X})^j/\Omega^d(\mathcal{X})^{j+1}]'$ is an (LB)-space. We now turn to the right hand side. Using [GKPS] Thm. 3.1.16 (compare also [Kom] Thm. 7') it suffices to show that $C^{\operatorname{an}}(G, P_{\underline{j}}; \operatorname{Hom}_K(M_{\underline{j}}, \operatorname{St}_{d+1-j}))$ is the locally convex inductive limi! t of a sequence of Banach spaces with compact transition maps. To see this it is convenient to identify this space, as a locally convex vector space (without the $G$-action), with the space $C^{\mathrm{an}}(G/P_{\underline{j}}, \mathrm{Hom}_K(M_{\underline{j}}, \mathrm{St}_{d+1-j}))$ of all locally analytic functions on $G/P_{\underline{j}}$ with values in $\mathrm{Hom}_K(M_{\underline{j}}, \mathrm{St}_{d+1-j})$. The recipe how to do this is given in [Fea] 4.3.1. One fixes a section $\imath$ of the projection map $G \twoheadrightarrow G/P_{\underline{j}}$ such that

$$\begin{array}{rcl} G/P_{\underline{j}} \times P_{\underline{j}} & \xrightarrow{\sim} & G \\ (gP_{\underline{j}}, h) & \longmapsto & \imath(gP_{\underline{j}})h \end{array}$$

is an isomorphism of locally analytic manifolds ([Fea] 4.1.1). We then have the continuous injection

$$\begin{array}{rcl} C^{\mathrm{an}}(G, P_{\underline{j}}; V) & \longrightarrow & C^{\mathrm{an}}(G/P_{\underline{j}}, V) \\ f & \longmapsto & [gP_{\underline{j}} \mapsto f(\imath(gP_{\underline{j}}))] \end{array}$$

writing $V := \mathrm{Hom}_K(M_{\underline{j}}, \mathrm{St}_{d+1-j})$ for short. In fact we will need that $V$ is of the form $V = V_{\mathrm{fin}} \underset{K}{\otimes} V_{\mathrm{sm}}$ for two $P_{\underline{j}}$-representations $V_{\mathrm{fin}}$ and $V_{\mathrm{sm}}$ which are finite dimensional algebraic and smooth, respectively. If $V_f$ runs over the finite dimensional subspaces of $V_{\mathrm{sm}}$ then

$$V = \varinjlim_{V_f} V_{\mathrm{fin}} \underset{K}{\otimes} V_f$$

and each $V_{\mathrm{fin}} \underset{K}{\otimes} V_f$ is invariant under some open subgroup of $P_{\underline{j}}$. A possible inverse of the above map has to be given by

$$\phi \longmapsto f_\phi(g) := (g^{-1}\imath(gP_{\underline{j}}))(\phi(gP_{\underline{j}})) \,.$$

Since

$$C^{\mathrm{an}}(G/P_{\underline{j}}, V) = \varinjlim_{V_f} C^{\mathrm{an}}(G/P_{\underline{j}}; V_{\mathrm{fin}} \underset{K}{\otimes} V_f)$$

it suffices to check that

$$\begin{array}{rcl} C^{\mathrm{an}}(G/P_{\underline{j}}; V_{\mathrm{fin}} \underset{K}{\otimes} V_f) & \longrightarrow & C^{\mathrm{an}}(G, P_{\underline{j}}; V) \\ \phi & \longmapsto & f_\phi \end{array}$$

is well defined and continuous. Consider the obvious bilinear map

$$\beta : [V_{\mathrm{fin}} \underset{K}{\otimes} V_f] \times [\mathrm{End}_K(V_{\mathrm{fin}}) \underset{K}{\otimes} \mathrm{Hom}_K(V_f, V_{\mathrm{sm}})] \longrightarrow V$$

between vector spaces equipped with the finest locally convex topology. By [Fea] 2.4.3 (the condition BIL is trivially satisfied) it induces a continuous bilinear map

$$C^{\mathrm{an}}(G/P_{\underline{j}}; V_{\mathrm{fin}} \otimes V_f) \times C^{\mathrm{an}}(P_{\underline{j}}, \mathrm{End}(V_{\mathrm{fin}}) \otimes \mathrm{Hom}(V_f, V_{\mathrm{sm}}))$$
$$\longrightarrow \quad C^{\mathrm{an}}(G/P_{\underline{j}} \times P_{\underline{j}}, V) \, .$$
$$(\phi, \Psi) \quad \longmapsto \quad \beta \circ (\phi \times \bar{\Psi})$$

Using the section $\imath$ we obtain the continuous bilinear map

$$\hat{\beta} : C^{\mathrm{an}}(G/P_{\underline{j}}; V_{\mathrm{fin}} \otimes V_f) \times C^{\mathrm{an}}(P_{\underline{j}}, \mathrm{End}(V_{\mathrm{fin}}) \otimes \mathrm{Hom}(V_f, V_{\mathrm{sm}})) \longrightarrow C^{\mathrm{an}}(G, V)$$

defined by $\hat{\beta}(\phi, \Psi)(g) := \beta(\phi(gP_{\underline{j}}), \Psi(\imath(gP_{\underline{j}})^{-1} g))$. It remains to observe that $\Psi_{\mathrm{o}}(h) := h^{-1}. \otimes h^{-1}.$ lies in $C^{\mathrm{an}}(P_{\underline{j}}, \mathrm{End}(V_{\mathrm{fin}}) \otimes \mathrm{Hom}(V_f, V_{\mathrm{sm}}))$ and that $\hat{\beta}(\phi, \Psi_{\mathrm{o}}) = f_{\phi}$.

We now are reduced to show that $C^{\mathrm{an}}(G/P_{\underline{j}}, V)$ is the locally convex inductive limit of a sequence of Banach spaces with compact transition maps. Since $G/P_{\underline{j}}$ is compact this is a special case of [Fea] 2.3.2. $\square$

To finish let us reconsider the bottom filtration step. By definition $\mathrm{St}_1 = K$ is the trivial representation, and $L(\underline{d}) = K^{\times}$ acts on the one dimensional space $M_{\underline{d}}$ through the character $a \mapsto a^d$. Let therefore $K_{\chi}$ denote the one dimensional $P_{\underline{d}}$-representation given by the locally analytic character

$$\chi : P_{\underline{d}} \quad \longrightarrow \quad K^{\times}$$
$$g \quad \longmapsto \quad \frac{\det(g)}{(g_{dd})^{d+1}} \, .$$

By comparing weights one easily checks that the natural map $U(\mathfrak{n}_{\underline{d}}^+) \underset{K}{\otimes} M_{\underline{d}} \longrightarrow \mathfrak{b}_{\underline{d}}/\mathfrak{b}$ is bijective which means that $\mathfrak{d}_{\underline{d}} = 0$. Our theorem therefore specializes in this case to the assertion that the map

$$I^{[d]} : [\Omega^d(\mathcal{X})^d]' \quad \overset{\cong}{\longrightarrow} \quad C^{\mathrm{an}}(G, P_{\underline{d}}; K_{\chi})$$
$$\lambda \quad \longmapsto \quad [g \mapsto -\lambda(g_*(d\Xi_{\beta_0} \wedge \ldots \wedge d\Xi_{\beta_{d-1}}))]$$

is a $G$-equivariant topological isomorphism.

Peter Schneider
Mathematisches Institut
Westfälische Wilhelms-Universität Münster
Einsteinstr. 62
D-48149 Münster, Germany
pschnei@math.uni-muenster.de
http://www.uni-muenster.de/math/u/schneider

Jeremy Teitelbaum
Department of Mathematics, Statistics, and Computer Science (M/C 249)
University of Illinois at Chicago
851 S. Morgan St.
Chicago, IL 60607, USA
jeremy@uic.edu
http://raphael.math.uic.edu/~jeremy